\documentclass[11pt]{article}

\usepackage[plainpages=false, colorlinks=true, 
            linkcolor=black, urlcolor=black, citecolor=black]{hyperref}
\usepackage{geometry}
\usepackage{xcolor,lipsum}
\usepackage{booktabs}
\usepackage{caption} 
\usepackage{subcaption} 
\usepackage{geometry} 
\usepackage[labelformat=simple]{subcaption}
\usepackage{tabularx}
\usepackage{tikz}
\usepackage{overpic}
\usepackage{amsmath, amssymb, amsthm,  amsfonts, mathrsfs, cite}
\usepackage{multirow}
\usepackage{color}
\usepackage{booktabs}
\usepackage{bm}
\usepackage{hyperref}
\usepackage{scrextend}
\usepackage{arydshln}
\usepackage {url}
\usepackage[font=scriptsize, labelfont=bf]{caption}
\usepackage{algorithm}
\usepackage{graphicx} 
\graphicspath{{./fig/}}
\usepackage{epstopdf} 
\usepackage{algpseudocode}
\usepackage{graphicx}
\usepackage{array,  tabularx}
\usepackage{algorithm, algcompatible}
\usepackage{booktabs} 
\usepackage{paralist} 
\usepackage{verbatim} 
\usepackage{xcolor, marginnote, enumitem}
\usepackage[numbers]{natbib}
\numberwithin{equation}{section}
\theoremstyle{theorem}
\newtheorem{lemma}{Lemma}
\newtheorem{theorem}{Theorem}
\newtheorem{proposition}{Proposition}

\theoremstyle{remark}
\newtheorem{remark}{Remark}

\theoremstyle{definition}
\newtheorem{definition}{Definition}

\DeclareMathOperator{\dist}{dist}

\DeclareMathOperator{\dom}{dom}

\DeclareMathOperator*{\argmin}{arg\, min}

\DeclareMathOperator{\sign}{sign}
\newcommand{\mM}{\mathcal{M}}
\newcommand{\mA}{\mathcal{A}}

\newcommand{\mP}{\mathcal{P}}
\newcommand{\mS}{\mathcal{S}}
\newcommand{\mF}{\mathcal{F}}

\newcommand{\mV}{\mathcal{V}}
\newcommand{\mD}{\mathcal{D}}

\newcommand{\bitem}{\begin{itemize}}
	\newcommand{\eitem}{\end{itemize}}

\newcommand{\bpm}{\begin{pmatrix}}

\DeclareMathAlphabet{\mathbfit}{OML}{cmm}{b}{it}

\usepackage{xspace}
\usepackage{bold-extra}
\usepackage[most]{tcolorbox}

\colorlet{texcscolor}{blue!50!black}
\colorlet{texemcolor}{red!70!black}
\colorlet{texpreamble}{red!70!black}
\colorlet{codebackground}{black!25!white!25}

\definecolor{Red}{rgb}{0.8,0,0}
\definecolor{Blue}{rgb}{0,0,0.8}
\definecolor{Green}{rgb}{0,0.4,0.4}
\definecolor{darkblue} {rgb} {0.00, 0.0, 1.0}

\DeclareMathAlphabet{\mathbfit}{OML}{cmm}{b}{it}

\usepackage{xspace}
\usepackage{bold-extra}
\usepackage[most]{tcolorbox}

\colorlet{texcscolor}{blue!50!black}
\colorlet{texemcolor}{red!70!black}
\colorlet{texpreamble}{red!70!black}
\colorlet{codebackground}{black!25!white!25}

\date{}

\begin{document}


\title{An Augmented Lagrangian Method-Based Framework in the Adjoint Space for Sparse Reconstruction of Acoustic Sources}

\author{Nirui Tan\thanks{School of  Mathematics,
		Renmin University of China, 100872 Beijing, China. Email: \href{mailto:niruitan@ruc.edu.cn}{tnr00@ruc.edu.cn}.}
	\and Hongpeng Sun\thanks{School of  Mathematics,
		Renmin University of China, 100872 Beijing, China.
		Email: \href{mailto:hpsun@amss.ac.cn}{hpsun@amss.ac.cn}.} \thanks{Corresponding author. }}

\date{} 

\maketitle

\begin{abstract}
We propose a semismooth Newton-based augmented Lagrangian framework for reconstructing acoustic sources in inverse acoustic scattering problems. Rather than working directly in the unknown source space, our semismooth Newton updates operate in the measurement (adjoint) space, which is especially efficient when the number of measurements is much smaller than the discretized source dimension. The source is then recovered via Fenchel–Rockafellar duality. Our approach substantially accelerates computation and reduces costs. Numerical experiments in two and three dimensions demonstrate the high efficiency of the proposed method.

\end{abstract}
\paragraph{Key words.}{inverse acoustic scattering, inverse source problem,  adjoint space method, semismooth Newton method, augmented Lagrangian method, global convergence}

\section{Introduction}\label{sec:intro}

Inverse acoustic scattering plays a key role in applications such as sonar imaging, oil prospecting, and nondestructive testing \cite{CK}. Here we focus on the inverse source problem: reconstructing acoustic sources from boundary measurements. It is well known that the inverse source problem lacks uniqueness \cite{BC,JS} due to the presence of non-radiating sources, making it highly ill-posed. Appropriate regularization that incorporates prior information can mitigate the ill-posedness (e.g.,  $L^2$ regularization \cite{DEL}). For convenience, let $\Omega \subseteq \mathbb{R}^d$ with $d=2$ or $d=3$ is a bounded and compact domain with boundary $\partial \Omega$ of class $C^3$ and contains the sources.  Sparse regularization is a popular choice in signal processing, including $L^1$ regularization and total variation, which can promote sparsity in the solutions or their gradients.  
 However, for acoustic scattering involving the Helmholtz equation, which is an infinite-dimensional problem, the existence of a reconstructed solution in the $L^1(\Omega)$ space cannot be guaranteed due to the lack of weak completeness in $L^1(\Omega)$ (see \cite[Chapter 4]{Bre}).  
 In \cite{KB, CLA2,CLA1},  a larger space, the Radon measure space $\mM(\Omega)$, is introduced for inverse problems and optimal control. We also refer to \cite{AS} for inverse problems, including electrical impedance tomography with
finite measurements and a linearized and locally optimized strategy and algorithms for sparse point acoustic sources \cite{APS}. 
 
It is known that $\mM(\Omega)$ is a Banach space and can be characterized by its dual space $C_0(\Omega)$ through the Riesz representation theorem (see \cite[Chapter 4]{Bre}),
\begin{equation}\label{eq:measure:dual}
\|\mu\|_{\mM(\Omega)} = \sup \left\{\int_{\Omega} u  d\mu : \ u \in C_{0}(\Omega), \ \|u\|_{C_0(\Omega)} \leq 1 \right\}.
\end{equation}
This is also equivalent to $\mM(\Omega) = C_0(\Omega)'$, which means that $\mM(\Omega)$ is weakly compact by the Banach-Alaoglu theorem since $C_0(\Omega)$ is a separable Banach space \cite{Bre}. The existence of a reconstructed sparse solution in the Radon measure space $\mM(\Omega)$ can be guaranteed. 
Moreover, since $ L^1(\Omega)$ can be embedded in $\mM(\Omega)$, regularization in $\mM(\Omega)$ can also promote sparsity. Inspired by  developments in the Radon measure space regularization \cite{CLA2,CLA1}, \cite{XS} considered the sparse reconstruction of an acoustic source in the Radon measure space $\mM(\Omega)$. The complete scattering data in $\Omega$ is employed, and a semismooth Newton method is developed for the reconstruction algorithms \cite{XS}. 

There are several iterative methods for inverse scattering with smooth functionals such as $L^2$ regularization, including the gradient descent method \cite{TA}, recursive linearization \cite{Bao1,Bao2}, Newton-type methods, and distorted Born iteration \cite{RB,chew1,HC,MCM,hoh1,hoh2}. For nonsmooth regularization, a primal-dual method for the total variation regularized inverse medium scattering problem is developed in \cite{BKL1,BKL2,CP}. 
Let us turn to a discrete setting to discuss the efficiency. Suppose that $N$ is the number of unknowns modeling the scatterer and $M$ is the number of measurements. As noted in \cite{chew1}, ``For most practical situations, the number of unknowns is much larger than the number of measurements, i.e., $N \gg M$".  In \cite{MCM}, it was also pointed out that $M \propto \sqrt{N}$ generally holds for the two-dimensional case. Currently, these Newton-type or distorted Born-type methods are iterated directly in the domain of unknowns, which means that a linear system such as $Ax=b$ with $A\in \mathbb{C}^{N\times N}$ must be solved during each Newton update. It can be much more efficient if each Newton update can be performed in the measurement domain by solving a linear system such as $A_b z=c$ with $A_b\in \mathbb{C}^{M\times M}$. 
The $N$ unknowns are then recovered from the $M$ computed dual variables directly via the optimality conditions. This technique can significantly reduce the computational cost. We also extend the proposed framework to the multi-frequency case, since there is some uniqueness with multiple frequencies \cite[Theorem 3.1]{Bao2} and \cite{Bao1}.
 
In this paper, we build a framework in which the Newton or semismooth Newton operates in the measurement domain. We call the corresponding solution (with the same dimension as the measurements) a ``dual solution" or ``adjoint field" \cite{DBP1}.  We then recover the unknown acoustic source through Fenchel–Rockafellar duality theory. Our contributions can be summarized as follows. First, we develop a semismooth Newton-based augmented Lagrangian method in which semismooth Newton updates are performed in the adjoint space. It enables us to solve linear systems with the same dimension as the boundary measurements.  Such ideas for inverse electromagnetic scattering can be found in \cite{DBP1,DBP2}, where incremental iterative methods are developed that involve a linear equation with the same dimensions as the measurements. Compared to \cite{DBP1, DBP2}, our framework can handle nonsmooth regularization, whereas $L^2$ Tikhonov regularization is employed in \cite{DBP1, DBP2}. The semismooth Newton method can thus be highly efficient when the number of boundary measurements is much smaller than the number of unknowns in the acoustic source to be reconstructed.  We would like to point out that a similar dual framework is studied in statistical or signal data processing \cite{LST, TS}. Second, we also design a semismooth Newton method with Moreau-Yosida regularization (see \cite[section 3.1]{CLA1}, \cite[Chapter 9]{KK}, or \cite{XS}) to directly reconstruct the acoustic source in cases where the scale of the measurement data is comparable to that of the unknowns.

The paper is organized as follows. In Section \ref{sec:direct}, we introduce the direct problem. In Section \ref{sec:alm}, we present the semismooth Newton-based augmented Lagrangian method in the adjoint space. In Section \ref{sec:ssn}, we introduce a semismooth Newton method and a first-order primal-dual method. In Section \ref{sec:num}, we present detailed numerical experiments both in two-dimensional and three-dimensional spaces. Finally, in Section \ref{sec:con}, we give a conclusion.  

\section{Inverse source problem with boundary measurements}\label{sec:direct}

In this paper, we will focus on the analysis and reconstruction of  the following inverse scattering problem:

\emph{Reconstructing an acoustic source $\mu$ in the Radon measure space $\mathcal{M}(\Omega)$ from a given scattered field on the boundary $\partial \Omega$}.

\begin{figure}
\centering
\usetikzlibrary {shapes.geometric}
\usetikzlibrary {decorations.shapes}
\tikzset{
  paint/.style={draw=black, fill=red},
  spreading/.style={
    decorate,decoration={shape backgrounds, shape=rectangle,
    shape start size=1.5mm,shape end size=1mm,shape evenly spread={#1}}}
}
\tikzset{my node/.style={trapezium, fill=green, draw=#1!75, text=black}}

\begin{tikzpicture}

    \foreach \i in {0,40,...,240} {
        \draw[fill=red] ({1.935},{\i/100.0 - 0.10}) rectangle ++(0.15,0.15);
    }

    \foreach \i in {0,40,...,240} {
        \draw[fill=red] ({-1.98},{\i/100.0 - 0.10}) rectangle ++(0.15,0.15);
    }

  \fill [paint=red,    spreading=11]            (2.0,2.8)   -- (-1.9,2.8);
 \fill [paint=red,    spreading=11]            (2.0,-0.5)   -- (-1.9,-0.5);


  \node [my node=red]                      at (0.1,1.1)  {$A$};

\end{tikzpicture}
\caption{Illustration of the acoustic source scattering: the red squares represent the receivers for recording the scattered wave, and the 
trapezoid labeled $A$ represents the acoustic sources.}\label{fig:illus}
\end{figure}
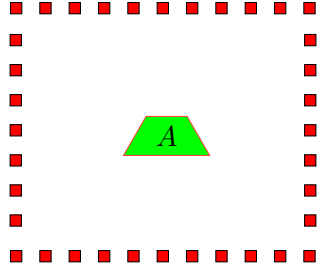

The direct scattering problem in the frequency domain with an inhomogeneous medium in $\mathbb{R}^{d}$ with $d=2$ or $d=3$ is as follows:
\begin{equation}\label{eq:helm}
\begin{cases}
-\Delta u  - k^2 n(x) u = \mu, \quad x \in \mathbb{R}^d, \\
\displaystyle{\lim_{|x| \rightarrow \infty} |x|^{\frac{d-1}{2}} (\frac{\partial u}{\partial |x|} - ik u) = 0,}
\end{cases}
\end{equation}
where $\mu \in \mM(\Omega)$ is a Radon measure and $n(x)$ is the refractive index. We refer to Figure \ref{fig:illus} for an illustration. Henceforth, we assume that $n(x)$ is real, smooth, and compactly supported, i.e., $\Im n = 0$ and there exists a bounded domain $ \Omega\subseteq \mathbb{R}^d$  such that $\text{supp}(n(x)-1) \Subset  \Omega$. 
 Throughout this paper, we assume that $\mu$ is a real, compactly supported Radon measure, which is reasonable in physics. We assume that the bounded domain $\Omega$ is large enough such that $\mu$ is also compactly supported in $\Omega$, i.e., 
$\text{supp}(\mu) \Subset \Omega$.
The Helmholtz equation \eqref{eq:helm} can also be considered for the case  $n(x) \equiv 1$.  

The solution of \eqref{eq:helm} is defined in a ``very weak" sense \cite[Definition 2.1]{XS} using test functions in $C^{2,\alpha}$. Here, for convenience, we directly work with the following volume potential representations, 
\begin{equation}\label{eq:weak:repre}
    (\mathcal{V} \mu)(x) :=\int_{\Omega}G(x,y)d\mu(y),
\end{equation}
where the inhomogeneous background Green's function $G(x,y)$ is defined as the radiating solution \cite{CM}
\begin{equation}\label{eq:green:nonhomo}
\begin{cases}
-\Delta_{x}G(x,y) - k^2n(x) G(x,y)=\delta(x-y), \quad x, y \in \mathbb{R}^d, \\
\displaystyle{\lim_{|x| \rightarrow \infty} |x|^{\frac{d-1}{2}} \Big(\frac{\partial G(x,y)}{\partial |x|} - ik G(x,y)\Big) = 0.}
\end{cases}
\end{equation}
 Denoting $m(x)=1-n(x)$,
we can thus construct $G(x,y)$ by the Lippmann-Schwinger integral equation \cite{CM}
\begin{equation}\label{eq:lip:inte}
G(x,y) =\Phi(x,y) -k^2 \int_{\Omega} \Phi(x,z)m(z)G(z,y)dz
\end{equation}
where $\Phi(x,y)$ is the fundamental solution of the Helmholtz equation, i.e.,  $ \Phi(x,y) =  \frac{i}{4}H_{0}^{(1)}(k|x-y|) $ in $\mathbb{R}^2$ and $\Phi(x,y) =  \frac{e^{ik|x-y|}}{4 \pi|x-y|}$ in $\mathbb{R}^3$.
For the $W^{1,p}$ estimate of the problem \eqref{eq:helm}, we have the following proposition \cite[Theorem 2.8]{XS}.
\begin{proposition}\label{prop:w1p}
    	For the solution of \eqref{eq:helm} with representation in \eqref{eq:weak:repre}, we have the following regularity estimate,
	\begin{equation}\label{eq:w1p:estimate}
	\|u\|_{W^{1,p}(\Omega)} \leq C \|\mu\|_{\mM(\Omega)}, \quad 1 \leq p < \frac{d}{d-1}, \quad d = 2 \ \text{or} \ d=3.
	\end{equation}
	Here $C$ is a positive constant that does not depend on $\mu$.
\end{proposition}



According to the trace theorem of $W^{1,p}(\Omega)$ in \cite[Theorem 5.36]{adam}, we have
\begin{equation}\label{eq:trace}
    W^{1,p}(\Omega) \hookrightarrow L^{s}(\partial \Omega),\quad  p \leq s \leq s^*=(d-1)p/(d-p)=(d-1)/(\frac{d}{p}-1).
\end{equation}
Since $p<d/(d-1)$ as in  \eqref{eq:w1p:estimate}, we have the estimate on $s$
\begin{equation}
    s = \begin{cases}
    s_0, \quad \forall s_0, \ \ \text{such that } \ p \leq s_0 < \infty, \ \ d=2, \\
    s_0, \quad \forall s_0, \ \ \text{such that } p \leq s_0 < 2, \ \ d=3.
    \end{cases}
\end{equation}
With the above trace property \eqref{eq:trace}, to reconstruct the acoustic source $\mu \in \mM(\Omega)$, we will make use of the following sparse regularization functional,
\begin{equation}\label{eq:sparse:functional}
\min_{\mu} J(\mu):=\frac{1}{2}\|\mV\mu-u_b\|_{L^s(\partial\Omega)}^2+\alpha\|  \mu\|_{\mM(\Omega)}, \quad    \forall s  \begin{cases}
    \text{such that } \ p \leq s < \infty, \ \ &d=2, \\
 \text{such that } p \leq s < 2, \ \ &d=3,
    \end{cases}
\end{equation}
where $\alpha$ is the regularization parameter and $u_b$ is the measured scattered fields on $\partial \Omega$. $\mV \mu$ satisfies the equation \eqref{eq:helm} as discussed.

\begin{theorem}\label{thm:existence:ori}
	There exists a solution $\mu \in \mM(\Omega)$ of the regularization functional \eqref{eq:sparse:functional}.
\end{theorem}

\begin{proof}
	The proof is similar to \cite[Proposition 2.2]{CLA1} and \cite{KB}.  Since the energy in \eqref{eq:sparse:functional} is $\frac{1}{2}\|u_b\|_{L^s(\partial \Omega)}^2$ when $\mu=0$, so we can find a minimizing sequence $\{\mu_n\}_n$ in $\mM(\Omega)$ that is bounded by $\frac{1}{2\alpha}\|u_b\|_{L^s(\partial \Omega)}^2$.
	Since $\mM(\Omega)$ is a weakly sequentially compact \cite{Bre} (see Chapter 4), there exists a weakly convergent subsequence $\mu_{n,k}$ converging weakly to $\mu^* \in \mM(\Omega)$.
	
	Denoting $u_{n,k} = \mV(\mu_{n,k})$, we see $u_{n,k} \in W^{1,p}(\Omega)$ with $1\leq p < \frac{d}{d-1}$ with \eqref{eq:w1p:estimate}. By Proposition \ref{prop:w1p}, 
    $u_{n,k}$ weakly converges to $\mV(\mu^*)$, and $\mV(\mu^*)$ is a solution \eqref{eq:helm} in the sense of the very weak solution \cite[Definition 2.1]{XS}, defined as the dual of $C^{2,\alpha}$. 
    By the trace theorem \cite[Theorem 5.36]{adam} and \eqref{eq:trace}, the trace operator is continuous from \(W^{1,p}(\Omega)\) to \(L^s(\partial \Omega)\), so \(u_{n,k}|_{\partial \Omega} \rightharpoonup V(\mu^*)|_{\partial \Omega}\) weakly in \(L^s(\partial \Omega)\). 
    The data term \(\frac{1}{2} \|\cdot - u_b\|^2_{L^s(\partial \Omega)}\) is convex and continuous, hence weakly lower semicontinuous.
    The regularization term \(\alpha \|\cdot\|_{M(\Omega)}\) is a norm, hence weak lower semicontinuous (see \cite[Chapter 4]{Bre}). Therefore,  
\[
J(\mu^*) \leq \liminf_{n \to \infty} J(\mu_n) = \inf J,  
\]  
so \(\mu^*\) is a minimizer  of \eqref{eq:sparse:functional} and the existence follows.
\end{proof}
For the nonsmooth minimization problem \eqref{eq:sparse:functional},  it is convenient to consider the predual problem under the powerful Fenchel–Rockafellar duality theory; see \cite{KB, KK, CLA2, CLA1} for its various applications in inverse problems and optimal control problems. The semismooth Newton method can be employed to compute dual problems efficiently. However, the problem \eqref{eq:sparse:functional} involves a complex-valued function. For the application of Fenchel–Rockafellar duality theory, we need to reformulate complex-valued operators and functions in terms of real matrix operators and real vector functions. 
Let us denote 
\begin{equation}\label{eq:vectorize:matrixization}
\mV = \mV_{R} + i\mV_{I} \ \text{with} \ 
V   = \begin{pmatrix}
\mV_{R} & -\mV_{I} \\
\mV_{I} & \mV_{R}
\end{pmatrix}, \ \ 
\mu := \begin{pmatrix}
\mu_{R}  \\\mu_{I}
\end{pmatrix}, \ \text{and} \   
u_b := \begin{pmatrix}
\Re u_{b}  \\ \Im u_{b}
\end{pmatrix},
\end{equation}
where $\mV_{R} = \Re (\mV)$, $\mV_{I} = \Im (\mV)$. We continue to use $\mu$ and $u_b$ as their complex realizations. 

\section{Semismooth Newton-based augmented Lagrangian method}\label{sec:alm}

Now, we consider the discrete setting of the functional in \eqref{eq:sparse:functional} and denote $\mV_b = V|_{\partial \Omega}$, i.e., the trace of the volume potential after complex realization. For the discretization of $\mV_{b}$, we employ the framework in \cite{vanikko}, which is implemented in \cite{BKL1,BKL2} with publicly available code. The Lippmann-Schwinger equation can be used to solve for the inhomogeneous medium; we provide some comments in Section \ref{sec:num}. Since $\mV_{b}$ may lack injectivity, for the convenience of the semismooth Newton methods and the associated optimization framework, we introduce an extra $L^2$ regularization, and work in $L^2(\partial \Omega)$ instead of $L^s(\partial \Omega)$  in discrete settings,
\begin{equation}\label{eq:sparse:functional:dis}
\min_{\mu } \mP(\mu),  \quad  \mP(\mu) :=\frac{1}{2}\|\mV_{b} \mu-u_b\|_{L^2(\partial\Omega)}^2+\frac{\alpha_{0}}{2}\|\mu\|_{L^2(\Omega)}^2+\alpha\|  \mu\|_{L^1(\Omega)}.
\end{equation}
In \eqref{eq:sparse:functional:dis}, we also replace $\mM(\Omega)$ by $L^1(\Omega)$ since these two norms are equal in a discrete setting. The $L^2$ term $\frac{\alpha_{0}}{2}\|\mu\|_{L^2(\Omega)}^2$ not only improves the regularity of $\mu$ but also can guarantee the uniqueness of the minimizer of the functional \eqref{eq:sparse:functional:dis} due to its strong convexity. The combination of the $L^1$ and $L^2$ \eqref{eq:sparse:functional:dis} can help choose a solution for the inverse source problem, without overcoming its inherent nonuniqueness.
Now, let us turn to the augmented Lagrangian method for solving \eqref{eq:sparse:functional:dis}.
First, set
\begin{equation}\label{eq:h-and-p}
    h(y)=\frac{1}{2}\|y-u_b\|_{L^2(\partial\Omega)}^2,\quad p(\mu)=\frac{\alpha_0}{2}\|\mu\|_{L^2(\Omega)}^2+\alpha\|  \mu\|_{L^1(\Omega)}.
\end{equation}
Let us first present the general Fenchel-Rockafellar duality theory \cite[Theorem 4.34]{KK} or \cite[Theorem 15.23]{HBPL}. Let $f$, $g$ be convex, proper, lower-semicontinuous functions in Hilbert spaces $Y$ and $X$, and $\Lambda :X\rightarrow Y$ is a linear, bounded, and continuous operator. Then, under certain conditions (such as $0 \in \text{sri} (\dom f -\Lambda (\dom g ))$ \cite[Theorem 15.23]{HBPL}, where ``sri" denotes strong relative interior and $\dom$ denotes the definition domain), the dual problem of the following primal minimization problem
\begin{equation}\label{eq:primal}
\min_{x} f(\Lambda x) + g(x)
\end{equation}
is 
\begin{equation}\label{eq:dual}
\sup_y(- f^*(y) - g^*(-\Lambda^* y)) \Leftrightarrow \inf_y f^*(y) + g^*(-\Lambda^* y)
\end{equation}
and we have the following properties on the optimal solutions $(x^*,y^*)$
\begin{align}
& f(\Lambda x^*) + g(x^*)=- f^*(y^*) - g^*(-\Lambda^* y^*), \notag \\
&\Lambda x^*  \in \partial f^*(y^*), \quad -\Lambda^*y^* \in \partial g(x^*). \label{eq:optimal:pd}
\end{align}
Here $x^*$ and $y^*$ are the solutions of the primal problem \eqref{eq:primal} and dual problem \eqref{eq:dual} respectively. $\Lambda^*:X \rightarrow Y$ is the adjoint operator of $\Lambda$. The optimality conditions \eqref{eq:optimal:pd} can be used to recover solutions between the primal and the dual solutions. 

We will first present the semismooth Newton-based augmented Lagrangian method, followed by the convergence analysis.

\subsection{Semismooth Newton-based augmented Lagrangian framework in the dual space}
Now let us turn to \eqref{eq:sparse:functional:dis}. It  can be reformulated as 
\begin{equation}\label{eq:primal:2}
    \min_{\mu } \   h(\mV_{b}\mu)+p(\mu).
\end{equation}
By direct calculations, we see that the dual functions of $h$ and $p$  are
\begin{equation}\label{eq:h:star}
    h^*(y) = \frac{1}{2}\|y\|_{L^2(\partial \Omega)}^2 + \langle y, u_b\rangle_{L^{2}(\partial \Omega)},
\end{equation}
and
\begin{equation}
    p^*(z)=  \langle z', z \rangle_{L^2(\Omega)} -\frac{\alpha_0}{2}\|z'\|_{L^2(\Omega)}^2-\alpha\|  z'\|_{L^1(\Omega)}, \quad z' = \mS_{\frac{\alpha_1}{\alpha_0}}(\frac{z}{\alpha_0}), 
\end{equation}
where $\mS$ is the soft thresholding operator
\begin{equation}
    [\mS_{\sigma}(z)]_i:= \begin{cases}
        0, & |z_i| \leq \sigma \\
        z_i-\sigma \text{sign}(z_i), & |z_i| > \sigma
    \end{cases}
\end{equation}
and $\text{sign}(\cdot)$ is the sign function. The explicit formulation of $p^*$ is not needed in what follows. 
By Fenchel–Rockafellar  duality, we obtain the dual problem of \eqref{eq:primal:2} 
\begin{equation}\label{eq:primaldual}
  \underset{y \in L^2(\partial \Omega)}{\min} \  \mD(y):=p^*(-\mV_b^{T}y)+h^*(y),
\end{equation}
which is equivalent to $ \max_{y \in L^2(\partial \Omega)}   ( - p^*(-\mV_b^{T}y)-h^*(y))$ by \eqref{eq:dual}.
Furthermore, letting $-\mV_b^{T}y=z$, the above minimization problem becomes the following linearly constrained optimization problem
\begin{equation}\label{eq:dual:constr}
\underset{y \in L^2(\partial \Omega), \ z \in L^2(\Omega)}{\min} \ p^*(z)+h^*(y),\quad \text{s.t.},  \quad \mV_b^{T}y+z=0.
\end{equation}
We now introduce the augmented Lagrangian function 
\begin{equation}\label{eq:alg:fuc}
     \quad L_{\sigma}(y,z;\lambda):=p^*(z)+h^*(y)+\langle \lambda,{\mV_{b}}^*y+z\rangle_{L^2(\Omega)}+\frac{\sigma}{2}\|\mV_b^{T}y+z\|_{L^2(\Omega)}^2.
\end{equation}
We arrive at the following equivalent inf-sup problem for \eqref{eq:primal:2} under certain conditions,
\begin{equation}
    {\inf_{y \in L^2(\partial \Omega), \ z \in L^2(\Omega)}} \  \ {\sup_{\lambda \in L^2(\Omega)}} \ L_{\sigma}(y,z;\lambda).
\end{equation}
The augmented Lagrangian method thus follows, with nondecreasing update of $\sigma_k \rightarrow c_{\infty} < +\infty$, e.g., $\sigma_{k+1}= c_0 \sigma_k$, where $c_0 \in [2,10]$ is a constant,
\begin{subequations}\label{eq:alm:y}
\begin{align}
    (y^{k+1}, z^{k+1}) &= \argmin_{y \in L^2(\partial \Omega), \ z \in L^2(\Omega)} L_{\sigma^k}(y,z;\lambda^k), \\
\lambda^{k+1} &= \lambda^k + \sigma_k(\mV_b^{T}y^{k+1}+z^{k+1}).\label{eq:update:lambda}
\end{align}
\end{subequations}

The augmented Lagrangian method (abbreviated as ALM) differs from ADMM (the alternating direction method of multipliers), as ALM solves the $(y,z)$  variables together, whereas ADMM solves $y$ and $z$ separately.
ALM can be derived from the proximal point algorithm \cite{Roc1} to the dual problem \cite{Roc2}, and increasing the step sizes $\sigma_k$ can lead to superlinear convergence for certain functions. ALM is also employed for parameter identifications for elliptic systems \cite{CZ}. For developments of ALM, we refer to \cite{HE,POW} and \cite{FG,GL1}. 

Semismooth Newton methods \cite{KK,MU} are widely employed for the nonlinear and nonsmooth updates of ALM \cite{LST,ZST} and \cite{Gor,HK1, KK1}. For its application in total variation (or total generalized variation) regularized image processing problems, we refer to \cite{sun1, sun2}. We now turn to semismooth Newton methods for solving the coupled nonlinear system of $(y,z)$ in \eqref{eq:alm:y}. The first-order optimality conditions on $y^{k+1}$ and $z^{k+1}$, respectively, in \eqref{eq:alm:y} are
\begin{align}  
&\nabla h^*(y^{k+1})+{\mV_{b}}\lambda^k+\sigma \mV_{b}(\mV_b^{T}y^{k+1}+z^{k+1})  = 0, \label{eq:coup:z}\\ 
&\partial p^*(z^{k+1})+\lambda^k+\sigma_k(z^{k+1}+\mV_b^{T}y^{k+1}) \ni 0. \label{eq:coup:p}
\end{align}

Usually, for discretizations, the dimensions of $L^{2}(\partial \Omega)$ are much smaller than the dimensions of $L^2(\Omega)$ \cite{chew1}. It can be much more efficient to solve only for the Newton update of $y^{k+1}$. Let us first solve for $z^{k+1}$ in terms of $y^{k+1}$.
By $ (\sigma_k I+\partial p^*)(z^{k+1})\ni-\sigma \mV_b^{T}y^{k+1}-\lambda^k $, we now solve $z^{k+1}$ from \eqref{eq:coup:p}
\begin{equation}\label{eq:solve:z}
 z^{k+1}=(I+\frac{1}{\sigma_k}\partial p^*)^{-1}(-\mV_b^{T}y^{k+1}-\frac{\lambda^k}{\sigma_k}).
\end{equation}
Substituting $z^{k+1}$ in \eqref{eq:solve:z} into \eqref{eq:coup:z}, we obtain
\[
\nabla h^*(y^{k+1})+\mV_{b}\lambda^k+\sigma \mV_{b}\mV_b^{T}y^{k+1}+\sigma_k \mV_{b}(I+\frac{1}{\sigma_k}\partial p^*)^{-1}(-\mV_b^{T}y^{k+1}-\frac{\lambda^k}{\sigma_k})\ni 0.
\]
With Moreau’s identity
$$x=(I+\sigma \partial G)^{-1}(x)+\sigma(I+\frac{1}{\sigma}\partial G^*)^{-1}(\frac{x}{\sigma}),$$
and by letting $x= -\sigma_k \mV_b^{T}y^{k+1} -\lambda^k$, and $G(x)=p(x)$, we get $\sigma_k z^{k+1}$ as follows
\begin{equation}\label{eq:update:zk+1}
\sigma_k(I+\frac{1}{\sigma_k}\partial p^*)^{-1}(-\mV_b^{T}y^{k+1}-\frac{\lambda^k}{\sigma_k})=-\sigma_k \mV_b^{T}y^{k+1} -\lambda^k - (I +\sigma_k \partial p)^{-1}(-\sigma_k\mV_b^{T}y^{k+1} -\lambda^k).
\end{equation}
Substituting the above formula (which is actually $\sigma_k z^{k+1}$) into  \eqref{eq:coup:z}, we arrive at
\[
\nabla h^*(y^{k+1})+\mV_{b}\lambda^k+\sigma_k \mV_{b}\mV_b^{T}y^{k+1}+ \mV_{b}\underbrace{(-\sigma_k \mV_b^{T}y^{k+1} -\lambda^k - (I +\sigma_k \partial p)^{-1}(-\sigma_k \mV_b^{T}y^{k+1} -\lambda^k))}_{\sigma_k z^{k+1}}=0.
\]
This finally leads to the following nonlinear equation for $y^{k+1}$ in $L^{2}(\partial \Omega)$ 
\begin{equation}\label{eq:h:ssn1}
  \partial h^*(y)-\mV_{b}(I+\sigma\partial p)^{-1}(-\lambda^k-\sigma_k \mV_b^{T}y)=0.   
\end{equation}
Since $\nabla h^*(y)=y+u_b$ with \eqref{eq:h:star}, the nonlinear equation \eqref{eq:h:ssn1} on $y^{k+1}$ becomes
\begin{equation} \label{eq:nonlinear:F}
  \mF_k(y)=0, \quad \mF_k(y):=  y+u_b-\mV_{b}(I+\sigma_k\partial p)^{-1}(-\lambda^k-\sigma_k \mV_b^{T}y)=0.
\end{equation}
Note that $y \in L^{2}(\partial \Omega)$ and the nonlinear equation  \eqref{eq:nonlinear:F} is on the boundary $\partial \Omega$, which can be of much smaller dimension compared to the source $\mu$ to be reconstructed.
Since $(I+\sigma\partial p)^{-1}(\cdot)$ is a semismooth function for any positive $\sigma$, we have the following lemma for the semismooth Newton derivative of $\mF$.
\begin{lemma}
    The semismooth Newton derivative of $\mF_k(y)$ is 
    \begin{equation}\label{eq:newton:deri}
        (\partial F_k(y))(z) = (I +  \frac{\sigma_k}{1+\sigma_k\alpha_0} \mV_{b} \mathcal{X}^k_{y} \mV_b^{T})(z)
    \end{equation}
    where $\mathcal{X}_{y}^k$ can be chosen as follows
    \begin{equation}\label{eq:kakka}
        \mathcal{X}^k_{y}: = \begin{cases}
            \{1\}, \quad &|\lambda^k + \sigma_k\mV_b^{T}y| > \sigma_k \alpha \\
            \{0\}, \quad &|\lambda^k + \sigma_k\mV_b^{T}y| < \sigma_k  \alpha \\
            [0,1], \quad &|\lambda^k + \sigma_k\mV_b^{T}y| = \sigma_k  \alpha
        \end{cases}.
    \end{equation}
\end{lemma}
\begin{proof}
For the resolvent $(I+\sigma_k\partial p)^{-1}$ in \eqref{eq:nonlinear:F}, let us introduce
\[
w:= (I+\sigma_k\partial p)^{-1}( \mu )={\argmin_{\tilde w}}\ \frac{\|\tilde w-  \mu  \|_2^2}{2}+ \frac{\sigma_k\alpha_0}{2}\|\tilde w\|_{L^2(\Omega)}^2+\sigma_k\alpha\| \tilde w\|_{L^1(\Omega)}.
\]
By the optimality condition, we have
\begin{align}
 & 0 \in w -  \mu  + \sigma_k\alpha_0 w+\sigma_k\alpha\partial\|w\|_1 \Leftrightarrow   \mu  \in(1+\sigma_k\alpha_0)w+\sigma_k\alpha\partial\|w\|_1   \\
  &\Rightarrow w=(I+\frac{\sigma_k\alpha}{1+\sigma_k\alpha_0}\partial\|.\|_1)^{-1}(\frac{ \mu }{1+\sigma_k\alpha_0})=\mS_{\frac{\sigma_k\alpha}{1+\sigma_k\alpha_0}}(\frac{ \mu }{1+\sigma_k\alpha_0}).
\end{align}
Now, introducing $l(y): =(I+\sigma_k\partial p)^{-1}(-\lambda-\sigma_k\mV_b^{T}y)$,  we can write its expression componentwise
\begin{equation}
[l(y)]_i = \begin{cases}
\frac{(-\lambda-\sigma_k\mV_b^{T}y)_i}{1+\sigma\alpha_0} - \frac{\sigma\alpha}{1+\sigma\alpha_0}, \quad &\text{if} \ \ 
(-\lambda-\sigma_k\mV_b^{T}y)_i > \sigma_k\alpha, \\
0 \quad &\text{if} \ \ |(-\lambda-\sigma_k\mV_b^{T}y)_i| \leq \sigma_k\alpha  \\
\frac{(-\lambda-\sigma_k \mV_b^{T}y)_i}{1+\sigma_k\alpha_0} + \frac{\sigma_k\alpha}{1+\sigma_k\alpha_0}, \quad &\text{if} \ \ 
(-\lambda-\sigma_k \mV_b^{T}y)_i<- \sigma_k\alpha.
\end{cases}
\end{equation}
It can be checked that $l(y)$ is a $PC^{1}$ function (piecewise differentiable function) \cite[Proposition 4.3.1]{Sch}.  Since $l(y)$ is continuously differentiable on $y_i$ when $|(-\lambda-\sigma_k\mV_b^{T}y)_i| > \sigma_k\alpha $, we have 
$[(\partial l)(z)]_i=   (-\frac{\sigma}{1+\sigma_k \alpha_0}\mV_{b} (z))_i$. Whereas  $|(-\lambda-\sigma_k \mV_b^{T}y)_i| < \sigma_k\alpha $, we obtain $[(\partial l)(z)]_i=  0$.
For the case $|(-\lambda-\sigma_k \mV_b^{T}y)_i| = \sigma_k\alpha$, by \cite[Theorem 4.3.1]{Sch},  we have $ s  (\frac{-\sigma_k}{1+\sigma_k \alpha_0}\mV_{b} (z))_i \in [(\partial l)(z)]_i $ for any $s\in [0,1]$. Since the soft thresholding operator $\mS_{\frac{\sigma_k\alpha}{1+\sigma_k\alpha_0}}$ is  anisotropic on $y$. We finally arrive at the semismooth Newton derivative $\partial F$ as in \eqref{eq:newton:deri}.
\end{proof}

 Now, let us give the following remark for the additional quadratic term in \eqref{eq:sparse:functional:dis}.
 \begin{remark}\label{rem:1}
 With \cite[Proposition 12.60]{Roc3}, $p^*$ is differentiable with Lipschitz gradient. From the Newton derivative \eqref{eq:newton:deri}, it can be seen that with nonzero $\alpha_0$, the Newton derivative can be better conditioned compared to $\alpha_0 \equiv 0$ without the $L^2$ term $\frac{\alpha_{0}}{2}\|\mu\|_{L^2(\Omega)}^2$,  especially while $\sigma_k$ tends to a very large $c_{\infty}$ as required by the ALM in \eqref{eq:alm:y}.
 \end{remark}

The Newton update for solving $y^{k+1}$ are as follows, starting with $y^0=y^k$ for $l=0,2,\ldots, $ 
\begin{equation}\label{eq:newtonupdate}
    \mathcal{N}(y^l) (y^{l+1} -y^l)= - F(y^l), \quad \mathcal{N}(y^l)\in \partial F(y^l).
\end{equation}
For the Newton derivative as in \eqref{eq:newton:deri}, we have the following remark.
\begin{remark}\label{rem:2}
    As seen from \eqref{eq:newton:deri}, the semismooth Newton derivative is positive definite and bounded. Thus for the Newton update in \eqref{eq:newtonupdate} and any $\mathcal{N}(y^l)\in \partial F(y^l)$, we have $\mathcal{N}(y^l) \succ I$ and $\mathcal{N}(y^l)^{-1}$ is uniformly bounded. The Newton update \eqref{eq:newtonupdate} is well-defined. The convergence of the semismooth Newton methods does not depend on the choice of the semismooth Newton derivative \cite[Theorem 8.16]{KK}. For the Newton derivative in \eqref{eq:newton:deri}, we simply choose $     \mathcal{X}^k_{y}=1$ for the case $|\lambda^k + \sigma_k\mV_b^{T}y| = \sigma_k  \alpha$ as in \eqref{eq:kakka}. 
\end{remark}
After obtaining $y^{l+1}$, we can calculate $z^{k+1}$ by \eqref{eq:update:zk+1} as follows
\begin{equation}\label{eq:z:up}
    z^{l+1} = \mM^k({y^{l+1}}), \quad \mM^k({y^{l+1}}):=\frac{1}{\sigma_k}[- \sigma_k\mV_b^{T}y^{l+1} -\lambda^k - (I +\sigma_k \partial p)^{-1}(-\sigma_k \mV_b^{T}y^{l+1} -\lambda^k)].
\end{equation}
Since Newton methods, including semismooth Newton methods, are locally convergent, line search is usually employed for globalization and to ensure convergence. We use the following Armijo line search strategy. By defining $d^l = y^{l+1}-y^{l}$, we proceed with an $\emph{aggressive}$ Armijo-type line search involving parameters $\beta^t>0$, $t = 0,1,2, \ldots$, $\beta \in (0,1)$, and $ c>0$. The objective is to find the smallest integer $t$ satisfying:
\begin{equation}\label{eq:line:search}
L_{\sigma_k}(y^l+\beta^t  d^l, \mM^k(y^l+\beta^t  d^l);\lambda^k) \leq L_{\sigma_k}(y^l, \mM^k(y^l);\lambda^k)  - c\beta^t  \|d^l\|_2^2. 
\end{equation}
The final step size $t_l := \beta^t $ obtained from the line search in \eqref{eq:line:search} is used to update of $y^{k+1}$ as follows:
\begin{equation}\label{eq:LSstep}
   y^{l+1} = y^l + t_l d^l.
\end{equation}
Once a stopping criterion for inner Newton iteration is satisfied, e.g.,  for some $l=l_0$, we set $y^{k+1} = y^{l_0}$ and obtain $z^{k+1}$ with $z^{k+1}=\mM^k(y^{k+1})$. 
After finishing the inner semismooth Newton iterations, we then update $\lambda^{k+1}$ as in \eqref{eq:update:lambda} for the outer ALM updates.  
Once the stopping criterion of the outer ALM iterations is satisfied with $(y^K,z^K,\lambda^K)$, 
the solution $\mu$ corresponding to the primal problem can be approximated by 
\begin{equation}\label{eq:mu}
    \mu^K =(I+\frac{\alpha}{\alpha_0}\partial \|.\|_1)^{-1}(\frac{-\mV_b^{T}y^K}{\alpha_0}),
\end{equation}
which has an explicit solution. Since, by the primal-dual optimality conditions of \eqref{eq:primaldual} at saddle-points $(y^*,\mu^*)$, we have 
\[
-\mV_b^{T}y^*\in \partial p(\mu^*) \Leftrightarrow -\mV_b^{T}y^* \in \alpha \partial \|.\|_1+\alpha_0 \mu^* \Leftrightarrow \mu^* =(I+\frac{\alpha}{\alpha_0}\partial \|.\|_1)^{-1}(\frac{-\mV_b^{T}y^*}{\alpha_0}).
\]
Finally, let us conclude this section with the proposed semismooth Newton-based ALM as the following Algorithm \ref{alg:alm}. We then discuss its convergence in the next section.

\begin{algorithm}[h]
	\caption{Boundary Measurements-Based Dual Augmented Lagrangian Method with Semismooth Newton Solver for  \eqref{eq:sparse:functional} (ALM-bd-SSN)}\label{alg:alm}
\begin{algorithmic}
\REQUIRE ~
     Given linear operator $\mV_{b}$ and $u_b$, regularization parameters $(\alpha_0, \alpha)$, parameter settings for ALM $\sigma_0$, $c_0$, parameters for line search $\beta$, $c$, and initial vectors $(y_0, z_0, \lambda_0)$, maximum outer ALM iteration number $K_{\text{max}}$, maximum SSN iteration number $L$.
    \STATE set $y=y_0,\sigma=\sigma_0,\lambda=\lambda_0$
    \WHILE{$0\leq i\leq K_{\text{max}}$  and some stopping criterion for outer ALM iterations is satisfied}
    \WHILE{$l\leq L$ and some stopping criterion for inner SSN is satisfied}
    \STATE Calculate $\mathcal{X}^k_{y}$ in \eqref{eq:kakka} and Newton derivative \eqref{eq:newton:deri}, and solve  \eqref{eq:newtonupdate} for SSN update.
    \STATE Do Armijo line search \eqref{eq:line:search} and update $y^{l+1}$ with \eqref{eq:LSstep}. 
    \ENDWHILE
    \STATE  Compute $z^{k+1}=\mM^k(y^{k+1})$ by \eqref{eq:z:up}. 
    \STATE Update $\lambda^{k+1}$ by \eqref{eq:update:lambda}.
    \STATE Update $\sigma_{k+1}=c_0 \sigma_k$.
    \ENDWHILE
      \STATE Obtain the reconstructed $\mu$ with \eqref{eq:mu} finally.
       \end{algorithmic}
\end{algorithm}

\subsection{Convergence of the SSN-based ALM on the dual space}
Henceforth, we assume all the spaces including $L^2(\partial \Omega)$, $L^2( \Omega)$, all the variables including $y$, $\mu$, $z$, and $\lambda$, and all the linear operators including $\mV_b$ and $\mV_b^{T}$ are in finite-dimensional settings. Now we introduce some basic definitions and properties of multivalued mappings from convex analysis \cite{DR, LST}. Let $F: X \rightrightarrows Y$ be a multivalued mapping. The graph of $F$ is defined as the set
\[
\text{gph} F: = \{ (x,y) \in X\times Y| y\in F(x)\}.
\]
The inverse of $F$, i.e., $F^{-1}:  Y \rightrightarrows X$ is defined as the multivalued mapping whose graph is $\{(y,x)| (x,y) \in \text{gph} F\}$. The distance from $x$  to the set $C\subset X$ is defined by
\[
\text{dist}(x,C): = \inf\{\|x-x'\|\ | \ x' \in C\}.
\]
For the local convergence rate of ALM, we need the metrical subregularity  \cite{DR}.
\begin{definition}[Metric Subregularity \cite{DR}]\label{def:metricregular}
	A mapping $F: X \rightrightarrows Y$ is called metrically subregular at $\bar x$ for $\bar y $ if $(\bar x, \bar y) \in \text{gph} F$ and there exists a modulus $\kappa \geq 0$  along with a neighborhood $U$ of $\bar x$ and $V$ of $\bar y$ such that
	\begin{equation}\label{eq:metricregular}
	\dist(x, F^{-1}(\bar y)) \leq \kappa  \dist(\bar y, F(x) \cap V ) \quad \text{for all} \ \ x \in U.
	\end{equation}
\end{definition}
Now, let us turn to the primal problem \eqref{eq:sparse:functional:dis} and its dual problem \eqref{eq:primaldual}.
For the $L^1$ norm, we see
\begin{equation}\label{eq:ani:li}
\|\mu\|_{1} =\sum_{i=1}^M |\mu_i|
\end{equation}
which is a polyhedral function. 
Introduce the Lagrangian function
\begin{equation}\label{eq:lag1}
l(y,z,\lambda) = p^*(z)+h^*(y)+\langle \lambda,{\mV_{b}}^*y+z\rangle_{L^2(\Omega)}.
\end{equation}
It is well known that $l$ is a convex-concave function on $(y,z,\lambda)$. Define the maximal monotone operator $T_{l}$ by 
\begin{equation}
T_{l}(y,z,\lambda) =\{(y',z',\lambda')|(y',z',-\lambda')\in \partial l(y,z,\lambda)\},
\end{equation}
and the corresponding inverse is given by
\begin{equation}
T_{l}^{-1}(y',z',\lambda') =\{(y,z,\lambda)|(y',z',-\lambda')\in \partial l(y,z,\lambda)\}.
\end{equation}
Here, adding the negative sign as in $-\lambda'$  is because $l$ is concave on $\lambda$.
\begin{theorem}\label{thm:metric:regular:lag}
	For the dual problem \eqref{eq:primaldual}, assuming the KKT system has at least one solution, then $T_{l}$ is metrically subregular at $( y^*, z^*, \lambda^*)^T$ for the origin. Similarly, assuming $(\partial \mP)^{-1}(0) \neq \emptyset$ as $\mP$ in \eqref{eq:sparse:functional:dis},  $\partial \mP$ is metrically subregular at $\lambda^*$ for the origin.
\end{theorem}
\begin{proof}
	With direct calculation, we have 
	\[
	T_{l}(y,z,\lambda) = (\partial h^*(y)+ \mV_b \lambda, \partial p^*(z)+ \lambda, -z-\mV_b^{T}y)^{T}: = \mathcal{A}(x) + \mathcal{B}(x),
	\]
	where
	\begin{equation}
	\mathcal{A}\begin{pmatrix}
	y\\z\\ \lambda
	\end{pmatrix}
	:=\begin{pmatrix}
	\partial h^* & 0& 0\\
	0 & \partial p^* &0 \\
	0 &0&0
	\end{pmatrix}\begin{pmatrix}
	y\\z\\ \lambda
	\end{pmatrix},
	\quad 
	\mathcal{B}\begin{pmatrix}
	y\\z\\ \lambda
	\end{pmatrix} := 
	\begin{pmatrix}
	0 & 0 &\mV_b  \\
	0&0&I\\
	-\mV_b^{T} & -I &0
	\end{pmatrix}
	\begin{pmatrix}
	y\\z\\ \lambda
	\end{pmatrix}
	.
	\end{equation}
	It is known that the  $\|\cdot\|_{1}$ \eqref{eq:ani:li} is a polyhedral convex function. We can see that $\partial h^*$ is a polyhedral function of $y$ and $\partial p^*$ is a polyhedral function of $z$. We thus conclude that the monotone operator $\mathcal{A}$ is polyhedral. Besides, the operator $\mathcal{B}$ is maximal monotone and linear. Thus $T_{l}$ is a polyhedral mapping \cite{Rob}. By the corollary in \cite{Rob}, we see $T_{l}$ is metrically subregular at $(y^*, z^*, \lambda^*)^{T}$ for the origin.

Let us now turn to the metric subregularity of $\partial \mP$, the subdifferential of the dual problem \eqref{eq:primaldual}. Assuming $(\partial \mP)^{-1}(0) \neq \emptyset$, since
\begin{equation}\label{eq:subgradient:dual}
(\partial \mP)(\mu) = \mV_b^{T}(\mV_b \mu-u_b) + \alpha_0 \mu + \alpha \partial \|\mu\|_1.
\end{equation}
We conclude that it is a polyhedral mapping of $\mu$. With the corollary in \cite{Rob},  $\partial \mP$ is metrically subregular.
\end{proof}
Here, we follow the standard stopping criteria for the inexact augmented Lagrangian method originated \cite{Roc1, Roc2}.
\begin{align}
& L_{\sigma^k}(y^{k+1},z^{k+1};\lambda^k) - \inf_{y,z} L_{\sigma^k}(y,z;\lambda^k) \leq \epsilon_k^2/(2\sigma_k), \quad \sum_{k=0}^{\infty}\epsilon_k < \infty, \label{stop:a}        \tag{A} \\
& L_{\sigma^k}(y^{k+1},z^{k+1};\lambda^k) - \inf_{y,z} L_{\sigma^k}(y,z;\lambda^k)  \leq  (\delta_k^2/(2\sigma_k))\|\lambda^{k+1}-\lambda^k\|^2, \quad \sum_{k=0}^{\infty}\delta_k < +\infty, \label{stop:b1} \tag{B1}\\
&\text{dist}(0, \partial L_{\sigma^k}(y,z;\lambda^k)|_{y,z}  ) \leq  (\delta_k'/\sigma_k)\|\lambda^{k+1} - \lambda^k\|, \quad 0 \leq \delta_k' \rightarrow 0. \label{stop:b2}\tag{B2}
\end{align}

Before discussing the local convergence rate, let us turn to the relation between $\mu$ and the Lagrangian multiplier $\lambda$. It is known that the saddle-points of the augmented Lagrangian functional \eqref{eq:alg:fuc} and the Lagrangian functional \eqref{eq:lag1} are the same. Let us focus on  \eqref{eq:lag1}. With direct calculation, we have
\begin{align}
  &\inf_{y,z}   \sup_{\lambda} l(y,z,\lambda)=p^*(z)+h^*(y)+\langle \mV_{b}^Ty +z, \lambda \rangle \Leftrightarrow \\
  & \sup_{\lambda} -[\sup_{y}\langle y, -\mV_{b} \lambda \rangle -h^*(y) + \sup_{z}\langle -\lambda,z \rangle -p^*(z)] \Leftrightarrow \\
  & \sup_{\lambda} - (h(-\mV_b \lambda) + p(-\lambda)) \Leftrightarrow  - \inf_{\lambda}(h(-\mV_b \lambda) + p(-\lambda)). \label{eq"lambda}
\end{align}
Compared to \eqref{eq:primal:2}, we conclude that optimal solutions $\mu^*$ of \eqref{eq:primal:2} and $\lambda^*$ of \eqref{eq"lambda} have the relation $\mu^*=-\lambda^*$. The convergence of $\mu^k$ is determined by $\lambda^k$. We can also obtain the convergence rate of $\mu^k$ from the convergence rate of $y^k$ by \eqref{eq:mu}. Since $(I+\frac{\alpha}{\alpha_0}\partial \|.\|_1)^{-1}$ is a firmly nonexpansive operator, we have
\[
\|\mu^k- \mu^{k+1}\| \leq c_0  \|y^k-y^{k+1}\|, \quad \|\mu^k- \mu^*\| \leq c_0 \|y^k-y^{*}\|, \quad c_0:={\|\mV_b^{T}\|}/{\alpha_0}.
\]
Henceforth, we only focus on the convergence rate of $(y^k,z^k,\lambda^k)$.

We denote $\mathcal{X}^P$  as the solution set of the problem \eqref{eq:sparse:functional:dis}.  With these preparations, we have the following global convergence and local convergence rate for the ALM method \eqref{eq:alm:y}, under the mild condition that the KKT system of \eqref{eq:lag1} has at least one solution, as in Theorem \ref{thm:metric:regular:lag}.
\begin{theorem}\label{thm:ani:KKT}
	For the dual problem \eqref{eq:primaldual} and corresponding ALM \eqref{eq:alm:y}, denote the iteration sequence $(y^k, z^k,\lambda^k)$ generated by ALM-bd-SSN with stopping criteria \eqref{stop:a}. The sequence $(y^k, z^k,\lambda^k)$  is bounded and  convergences to $(y^*, z^*, \lambda^*)$ globally. $T_d=\partial \mP$ is  metrically subregular for the origin with modulus $\kappa_d$ and  with the additional stopping criteria \eqref{stop:b1},  the sequence $\{\lambda^k\}_k$ converges to  $\lambda^* \in \mathcal{X}^P$. For  sufficiently large $k$, we have the following local linear convergence
	\begin{equation}\label{eq:convergence:rate:dual:ani}
	\dist(\lambda^{k+1}, \mathcal{X}^{P}) \leq \theta_k \dist(\lambda^k, \mathcal{X}^P),
	\end{equation}
	where
	\[
	\theta_k = [\kappa_d(\kappa_d^2 + \sigma_k^2)^{-1/2} + \delta_k](1-\delta_k)^{-1}, \ \emph{as} \ k\rightarrow \infty, \  \theta_k \rightarrow  \theta_{\infty} = \kappa_d(\kappa_d^2+ \sigma_{\infty}^2)^{-1/2} < 1.
	\]
	Furthermore, $T_l$ is metrically subregular at $(y^*, z^*, \lambda^*)$ for the origin with modulus $\kappa_l$. When the additional stopping criteria \eqref{stop:b2} is employed, for sufficiently large $k$, we have
	\begin{equation}\label{eq:convergence:rate:up:ani}
	\|(y^{k+1}, z^{k+1}) - (y^k, z^k)\| \leq \theta_k'\|\lambda^{k+1}-\lambda^k\|,
	\end{equation}
	where $\theta_k'=\kappa_l(1+\delta_k')/\sigma_k$ with $\displaystyle{\lim_{k\rightarrow \infty}\theta_k' = \kappa_l/\sigma_{\infty}}$.
\end{theorem} 
\begin{proof}
	Since discrete $L^2(\partial \Omega)$ is a finite-dimensional reflexive space and the dual function \eqref{eq:primaldual} is lower semicontinuous, proper, and strongly convex due to the strong convexity of $h^*$, it is coercive. Thus, the existence of the solution can be guaranteed {\cite[Theorem 4.25]{KK}}. Furthermore, since $\dom \mD = L^2(\partial \Omega)$, by Fenchel-Rockafellar theory {\cite[Chapter 4.3]{KK}}, the solution to the dual problem \eqref{eq:primaldual} is not empty and $\min_{\mu} \mP(\mu)=\max_{y}-\mD(y)$.
	By {\cite[Theorem 4]{Roc2}} (or {\cite[Theorem 1]{Roc1}}, where the augmented Lagrangian method \eqref{eq:alm:y} essentially comes from proximal point method applying to its dual problem, i.e., the primal problem \eqref{eq:sparse:functional:dis}), with criterion \eqref{stop:a}, we obtain the boundedness of $\{\lambda^k\}_k$. The uniqueness of $(y^*,z^*)$ follows from the strong convexity of $h^*$ on $y$ and the $z^*=-\mV_b^{T} y^*$, which is one of the KKT conditions. The boundedness of $(y^k,z^k)$ and convergence of $(y^k, z^k,  \lambda^k)$ then follow by {\cite[Theorem 4]{Roc2}}. 
	
	By Theorem \ref{thm:metric:regular:lag}, we have metrical subregularity of $T_d= \partial \mP$. With the stopping criteria  \eqref{stop:a} and  \eqref{stop:b1},  the local convergence rate  \eqref{eq:convergence:rate:dual:ani}   can thus be obtained from {\cite[Theorem 5]{Roc2}} (or {\cite[Theorem 2]{Roc1}}). Now we turn to the local convergence rate of $(y^k,z^k)$. By the metrical subregularity of $T_l$ as in Theorem \ref{thm:metric:regular:lag}, for sufficiently large $k$, we have
	\[
	\| (y^{k+1}, z^{k+1}) - (y^*, z^*)\| + \dist (\lambda^k, \mathcal{X}^P) \leq \kappa_l \text{dist} (0, T_l(y^{k+1}, z^{k+1}, \lambda^{k+1})).
	\]
	Together with the stopping criteria \eqref{stop:b2} and \cite{Roc2} (Theorem 5 and Corollary in Section 4), we arrive at 
	\begin{align*}
	\| (y^{k+1}, z^{k+1}) - (y^*, z^*)\| &\leq \kappa_l \sqrt{{\delta_k'^2}{\sigma_k^{-2}} \|\lambda^{k+1} -\lambda^k\|^2 + {\sigma_k^{-2}} \|\lambda^{k+1} -\lambda^k\|^2 } \\
	& = \kappa_l\sqrt{\delta_k'^2+1}\sigma_k^{-1}\|\lambda^{k+1} -\lambda^k\| \leq \theta_k' \|\lambda^{k+1} -\lambda^k\|,
	\end{align*}
	which leads to \eqref{eq:convergence:rate:up:ani}.
	\qed
\end{proof}
\section{Semismooth Newton method and first-order primal-dual method}\label{sec:ssn}
In this section, we will introduce a semismooth Newton method and a first-order primal-dual method for the inverse source scattering problem. Let us first discuss the semismooth Newton method.
For the primal formulation \eqref{eq:sparse:functional:dis}, by choosing 
\begin{equation}
    f(z) :=  \frac{1}{2}\langle  B^{-1}z, z \rangle_{L^2(\Omega)} -\langle z, B^{-1}\mV_{b}^T u_b \rangle_{L^2(\Omega)} + d_0, \quad g(z): = \alpha \|z\|_{L^1(\Omega)},
\end{equation}
where $B:=\mV_{b}^T \mV_{b}+\alpha_0 I$.
The primal formulation can be written as $\min_{\mu} f(B\mu) +g(\mu)$.
By the Fenchel–Rockafellar  duality, the predual problem is $\min_{y}f^*(y) + g^*(-B^Ty)$, which is
\begin{align}
    \underset{y}{\min}&\langle y,By+\mV_{b}^T u_b\rangle_{L^2(\Omega)}-\frac{1}{2}\langle y+B^{-1}\mV_{b}^Tu_b,By+\mV_{b}^T u_b\rangle _{L^2(\Omega)}+\langle y+B^{-1}\mV_{b}^T u_b,\mV_{b}^Tu_b\rangle_{L^2(\Omega)}  \notag  \\
    &-d_0+\mathcal{I}_{\{\|B^T y\|_\infty\leq \alpha\}}(y)
\end{align}
We then use Moreau–Yosida regularization for dealing with the constraint $\|B^T y\|_\infty\leq \alpha$ as follows
\begin{align}
   & \underset{y}{\min} \ E(y), \ \ E(y):=\langle y,By+\mV_{b}^Tu_b\rangle_{L^2(\Omega)}-\frac{1}{2}\langle y+B^{-1}\mV_{b}^Tu_b,Bz+\mV_{b}^Tu_b\rangle _{L^2(\Omega)} 
   \label{eq:func:ssn2}\\
   &+\langle y+B^{-1}\mV_{b}^Tu_b,\mV_{b}^Tu_b\rangle_{L^2(\Omega) } -d_0+\frac{1}{2\gamma}\|\max(0,\gamma(B^T y-\alpha)\|_2^2+\frac{1}{2\gamma}\|\min(0,\gamma(B^T y+\alpha)\|_2^2. \notag
\end{align}
Although $B^{-1}$ appears in the  functional \eqref{eq:func:ssn2}, direct calculation of the first-order optimality condition shows that no $B^{-1}$ is involved in the following nonlinear equation, fortunately, 
\begin{equation}
    F(y)=\nabla E(y)=By +\mV_{b}^Tu_b+\gamma B\max(0,B^T y-\alpha)+\gamma B\min(0,B^T y+\alpha)=0.
\end{equation}
The formulation $\max(0,B^Ty-\alpha)$ or $\min(0,B^Ty+\alpha)$ is understood componentwise. For their subgradients, we introduce 
 $\mathcal{A}_k^+=\{x\in \Omega:B^Ty^k(x) \geq \alpha\}$, $\mathcal{A}_k^-=\{x\in \Omega:B^Ty^k(x)\leq -\alpha\}$, and $\mathcal{A}_k=\mathcal{A}_k^+\cup \mathcal{A}_k^-$ along with $\chi_{\mA^{+}}=\text{Diag}([\chi_{\mA^+}]_1,[\chi_{\mA^+}]_2, \ldots, )$, $\chi_{\mA^{-}}=\text{Diag}(\chi_{\mA^1}]_1,[\chi_{\mA^1}]_2, \ldots, )$ and $\chi_{\mA}=\text{Diag}([\chi_{\mA}]_1,[\chi_{\mA}]_2, \ldots, )$,  which depend on  $y$ and their diagonal elements are defined by
\begin{equation}\label{eq:katta:def2}
[\chi_{\mA^+}]_i = \begin{cases}
1, \quad [B^Ty^k]_i \geq \beta , \\
0, \quad  [B^Ty^k]_i < \beta ,
\end{cases}
\ 
[\chi_{\mA^-}]_i = \begin{cases}
1, \quad   [B^T y^k]_i \leq -\beta, \\
0, \quad  [B^Ty^k]_i > -\beta,
\end{cases}
\  [\chi_{\mA}]_i =  [\chi_{\mA^+}]_i + [\chi_{\mA^-}]_i.
\end{equation}
With these preparations, the semismooth Newton method for solving the nonlinear system $F(y) = 0$ reads as
 \begin{equation}
    \mathcal{N}(y^k)y^{k+1}=\mathcal{N}(y^k)y^k-F(y^k), 
 \end{equation}
 where $B+\gamma B\mathcal{X}_{\mathcal{A}_k}B^T=\mathcal{N}(y^k)\in \partial F(y^k)$, so we get
 \begin{equation}\label{eq:newton:up2}
    (B+\gamma B\mathcal{X}_{\mathcal{A}_k}B^T)y^{k+1}=-\mV_{b}^Tu_b + \gamma\alpha B(\mathcal{X}_{\mathcal{A}_k^+}-\mathcal{X}_{\mathcal{A}_k^-})\vec{1}.
 \end{equation}

We let $\gamma  \rightarrow \infty$ using a path-following strategy \cite{HK1} (or \cite[Section 9.1]{KK}). Finally, we can recover $\mu$ with $y$ by \eqref{eq:optimal:pd}
\begin{equation}
    B\mu^* = \nabla f^*(y) \rightarrow \mu = B^{-1} \nabla f^*(y^*)=y^*+B^{-1}(\mV_b^{T}u_b).
\end{equation}
We conclude the above discussions with the following semismooth Newton algorithm, i.e., Algorithm \ref{alg:ssn2}.
 \begin{algorithm}
    \caption{Semismooth Newton Method  (abbreviated as SSN)}
    \label{alg:ssn2}
    \begin{algorithmic}[1]
        \REQUIRE $y^0\in L^2(\Omega)$, strictly increasing path-following parameters $[\gamma_0,\gamma_1,\ldots,\gamma_{{K_{\text{max}}}}]$, with $\gamma_i>0$ $i=1,\ldots, K_{\text{max}}$. 
        \ENSURE $y,\mu$
        \STATE Initialization $y_{\gamma_0}^0=y^0$
        \WHILE{$0\leq i\leq I,\gamma=\gamma_i$}
        \WHILE{$k\leq K$}
        \STATE Calculate active sets and the corresponding functions as in \eqref{eq:katta:def2} for $\gamma=\gamma_i$. 
        \STATE Solve for $y^k\in  L^2(\Omega)$ by \eqref{eq:newton:up2} and denote it $y^k_{\gamma_i}$.
        \STATE Update $\mathcal{A}_k^+,\mathcal{A}_k^-,\mathcal{A}_k$
        \STATE Until $\mathcal{A}_k^+=\mathcal{A}_{k-1}^+,\mathcal{A}_k^-=\mathcal{A}_{k-1}^-$, set $y^0_{\gamma_{i+1}}=y^k_{\gamma_i}.$
        \ENDWHILE
        \ENDWHILE
    \STATE $y=y_{\gamma_{K_{\text{max}}}}^k$
    \STATE $\mu=B^{-1} \nabla f^*(y)=y+B^{-1}(\mV_b^{T}u_b)$.
    \end{algorithmic}
\end{algorithm}

\subsection{First-order primal-dual method}
Finally, let us turn to the first-order methods. The Chambolle-Pock first-order primal-dual algorithm \cite{CP} is employed and studied for the inverse medium scattering problem in \cite{BKL1,BKL2}. Here, for comparison, we also present the first-order primal-dual algorithm. With
\[
\frac{1}{2}\|\mV_{b}\mu-u_b\|_{L^2(\partial\Omega)}^2=\underset{p\in L^2(\partial\Omega)}{\max}\langle \mV_{b}\mu-u_b,p\rangle_{L^2(\partial\Omega)}-\frac{1}{2}\|p\|_{L^2(\partial\Omega)}^2,
\]
we get the following primal-dual problem
\begin{align}
&\underset{\mu}{\min}\ \underset{p}{\max}\langle \mV_{b}\mu-u_b,p\rangle_{L^2(\partial\Omega)}-\frac{1}{2}\|p\|_{L^2(\partial\Omega)}^2+\frac{\alpha_0}{2}\|\mu\|_{L^2(\Omega)}^2+\alpha\|\mu\|_{L^1(\Omega)} \Leftrightarrow\\
&   \underset{\mu}{\min} \ \underset{p}{\max}\langle \mV_{b}\mu,p\rangle_{L^2(\partial\Omega)}+\frac{\alpha_0}{2}\|\mu\|_{L^2(\Omega)}^2+\alpha\|\mu\|_{L^1(\Omega)}-\Big(\frac{1}{2}\|p\|_{L^2(\partial\Omega)}^2+\langle u_b,p\rangle_{L^2(\partial\Omega)}^2\Big).\label{eq:pda:sd}
\end{align}
Setting $F^*(p)=\frac{1}{2}\|p\|_{L^2(\partial\Omega)}^2+\langle u_b,p\rangle_{L^2(\partial\Omega)}$ and  $G(\mu)=\frac{\alpha_0}{2}\|\mu\|_{L^2(\Omega)}^2+\alpha\|\mu\|_{L^1(\Omega)}$, we have 
\[
(I+\sigma \partial F^*)^{-1}(\bar p)=\frac{\bar p-\sigma u_b}{1+\sigma}, \quad (I+\tau \partial G)^{-1}(\bar \mu)=\sign(\bar \mu)\max\Big\{\frac{|\bar \mu|-\tau\alpha}{1+\tau\alpha_0},0\Big\}.
\]
We now can employ the first-order primal-dual algorithm \cite[Algorithm 1]{CP} to solve \eqref{eq:sparse:functional} in its saddle-point formulation \eqref{eq:pda:sd},  as in the following  Algorithm \ref{alg:pda}.

\begin{algorithm}[H]
\caption{First-Order Primal-Dual Algorithm (abbreviated as PDA)}
\label{alg:pda}
\begin{algorithmic}[1] 
\REQUIRE ~~ 
    initial vectors $(p_0,\mu_0)\in L^2(\partial\Omega)\times L^2(\Omega)$, step sizes $(\sigma, \tau)$,
\ENSURE ~~ 
   $p$, $\mu$
    \STATE set $p=p_0$, $\mu=\mu_0$, $\bar{\mu}=\mu_0$
    \FOR{$0\leq k\leq  K_{\text{max}}$}
    \STATE $ p_{k+1} = (I  +\sigma \partial F^*)^{-1}(p_k+\sigma \mV_b \bar \mu_k)= \frac{p_k+\sigma \mV_{b} \bar{\mu}_k-\sigma u_b}{1+\sigma}$
    \STATE $\mu_{k+1} =(I+\tau \partial G)^{-1}(\mu_k-\tau \mV_b^{T}p_{k+1})=\sign(\mu_k-\tau \mV_b^{T}p_{k+1})\max\big\{\frac{|\mu_k-\tau \mV_b^{T}p_{k+1}|-\tau\alpha}{1+\tau\alpha_0},0\big\}$
    \STATE $ \bar{\mu}_{k+1} = \mu_{k+1} + \theta (\mu_{k+1} - \mu_k)$
    \ENDFOR
   \end{algorithmic}
\end{algorithm}

\section{Numerical experiments}\label{sec:num}
In this section, we will first give a brief discussion of the discretization of the Lippmann-Schwinger volume potential, followed by the presentation of the numerical examples. We will first give the linear mapping from the acoustic source $\mu$ to the boundary measurement $u_b^s$ for the inhomogeneous medium.
 Since we have
\begin{equation}
-\Delta u^s- k^2n(x) u^s=\mu \Rightarrow
-\Delta u^s - k^2 u^s = k^2(\frac{1}{k^2}\mu+(n-1)u^s),
\end{equation}
which leads to the following Lippmann-Schwinger integral equation
\begin{equation}\label{eq:sac:lip}
    u^s=\mV_k (\frac{1}{k^2}\mu+(n(x)-1)u^s).
\end{equation}
Now let us introduce $q(x)=n(x)-1$ and

\[
\mV_k  h: = k^2\int_{\Omega} \Phi(x,y)h(y)dy.
\]
We  have
\begin{equation}\label{eq:lipman1}
u^s(x) = \frac{1}{k^2}(I-\mV_k q)^{-1} \mV_k \mu. 
\end{equation}
Substituting \eqref{eq:lipman1} into \eqref{eq:sac:lip}, we have
\begin{equation}
        u^s(x)=\frac{1}{k^2}\mV_k (I+q(I-\mV_k q)^{-1} \mV_k)\mu.
\end{equation}
We finally obtain the linear mapping from the source $\mu$ to the measurement 
\begin{equation}
    u_b^s(x) = T_r \circ \big[\frac{1}{k^2}\mV_k (I+q(I-\mV_k q)^{-1} \mV_k)\big]\mu,  \quad T_r u:=u|_{ \partial \Omega},
\end{equation}
where $T_r$ denotes the trace mapping. The source to measurement mapping can be 
\begin{equation}
    \mathcal{T} : =  T_r \circ \big[\frac{1}{k^2}\mV_k (I+q(I-\mV_k q)^{-1} \mV_k)\big]: \quad \mathcal{M}(\Omega) \rightarrow T_r(W^{1,p}(\partial \Omega)).
\end{equation}
The volume potential $\mV_k$ and $(I-\mV_kq)^{-1}$ can be computed using a Fourier-based collocation method via periodization of the Lippmann-Schwinger volume potential \cite{vanikko}. We adapted the code developed in \cite{BKL1,BKL2} where the discretization of $  \mathcal{T}$ is given, and we refer to \cite{BKL1,BKL2} for more details. Integral equations are widely used for inverse scattering \cite{CK}; see \cite{Qzz} for the corresponding developments of rough surfaces scattering.

All algorithms are implemented in MATLAB R2024b on an Ubuntu 24.04 workstation with dual Intel Xeon E5-2697A v4 CPUs (2.60 GHz) and 288 GB of memory. 
The parameters  of the compared algorithms  are as follows

\begin{itemize}
\item ALM: ALM-bd-SSN as in Algorithm \ref{alg:alm}. We choose $\sigma_0=1$ and $\sigma_{k+1}=6\sigma_k$ with $c_0=6$. For the Armijo line search, we choose $\beta=0.3$ and $c=10^{-4}$ as in \eqref{eq:line:search}.
\item SSN: SSN as Algorithm \ref{alg:ssn2}. We choose $\gamma_0=1$ and $\gamma^i=10^i$,  $i=1,\ldots, 8$.
\item PDA: First-order primal-dual algorithm as in Algorithm \ref{alg:pda} \cite{CP}. We choose $\sigma=0.5$, $\tau= 1/((\|\mV_b^2\|+10^{-6})\sigma)$, and $\theta\equiv 1$.
\end{itemize}
Gaussian noise is added to the simulated data $u$ to obtain noisy boundary measurements $u_b^{\delta}$,
which is essentially the same as in \cite[(73)]{BKL1} 
\[
 (u_b^{\delta} - u_b)/\|u_b\|_2 = \delta (N_{\text{re}} +iN_{\text{im}})
\]
where $N_{\text{re}}$ and $N_{\text{im}}$ are two real vectors sampled from 
standard normal distribution. We measure the reconstruction error 
 by the following relative error 
\[
\text{N-Error}:=\|\mu^{K_{\text{max}}}-\mu_{\text{exa}}\|_2/\|\mu_{\text{exa}}\|_2,
\]
where $\mu_{\text{exa}}$ is the exact acoustic source and $K_{\text{max}}$ is the corresponding maximum outer iteration number, and it varies for different algorithms.

To avoid ``inverse crime", we use the implementation from  \cite{BKL1,BKL2}. The strategy is to use different discretizations (especially grid sizes) for generating the synthetic data in the direct problem and for reconstructing the unknown sources \cite{BKL1} (see also \cite[Chapter 2.3.6]{TA}).

For the homogeneous medium, we choose the velocity  $c\equiv 1$ both in $\mathbb{R}^2$ and $\mathbb{R}^3$. For the inhomogeneous medium,  we employ the one from \cite{hoh1}
\begin{equation}\label{eq:medium:in}
    q(x)=\chi(x-\bold{0}), \quad \chi(t)=
\begin{cases}
\exp(\frac{-1}{1-|t|^2}) & |t| \leq  1,\\
  0                   & |t| > 1,
\end{cases}
\end{equation}
where $\bold{0}=(0,0)^T$ in  $\mathbb{R}^2$ and $\bold{0}=(0,0,0)^T$ in $\mathbb{R}^3$.

\begin{figure}[htbp]
    \centering

    \begin{subfigure}[b]{0.3\textwidth}
        \includegraphics[width=\textwidth]{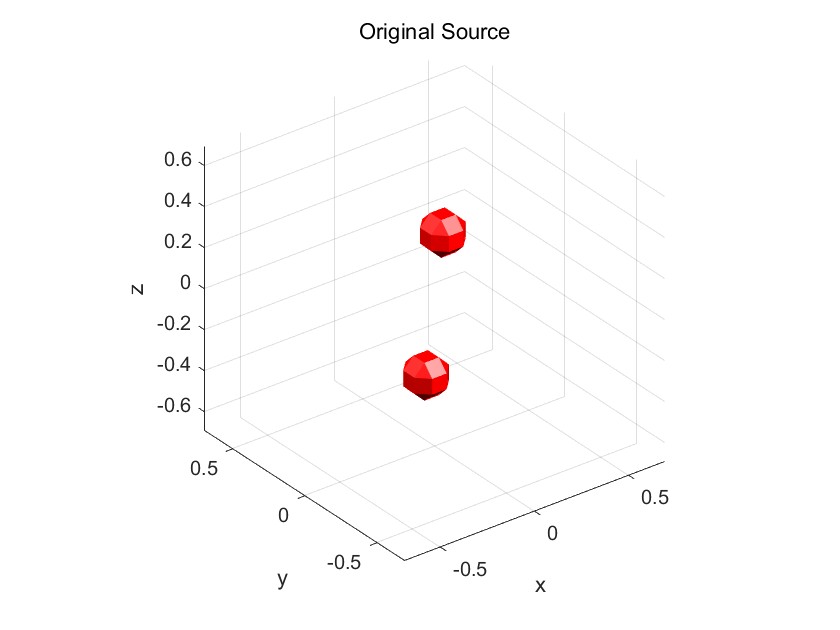}
        \caption{}
        \label{fig:1d}
    \end{subfigure}
    \hfill
    \begin{subfigure}[b]{0.3\textwidth}
        \includegraphics[width=\textwidth]{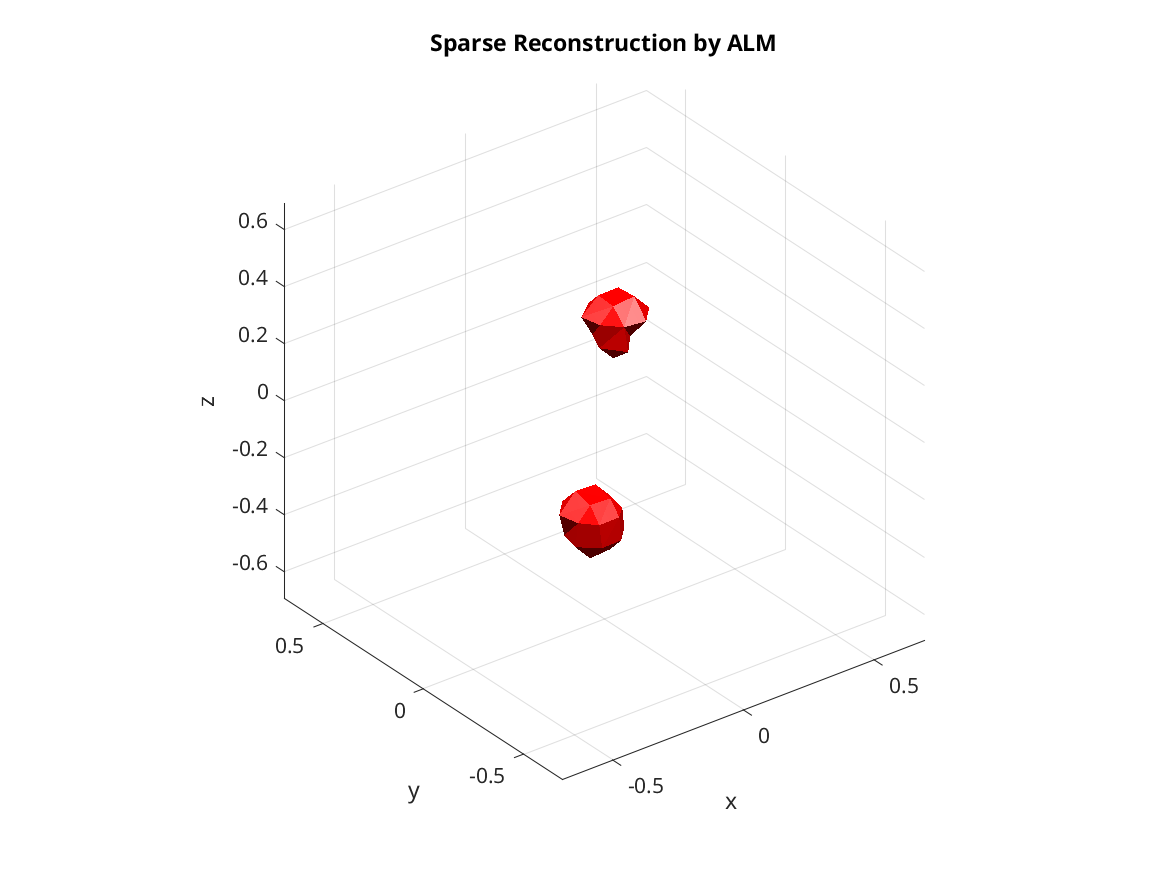}
        \caption{}
        \label{fig:1e}
    \end{subfigure}
    \hfill
    \begin{subfigure}[b]{0.3\textwidth}
        \includegraphics[width=\textwidth]{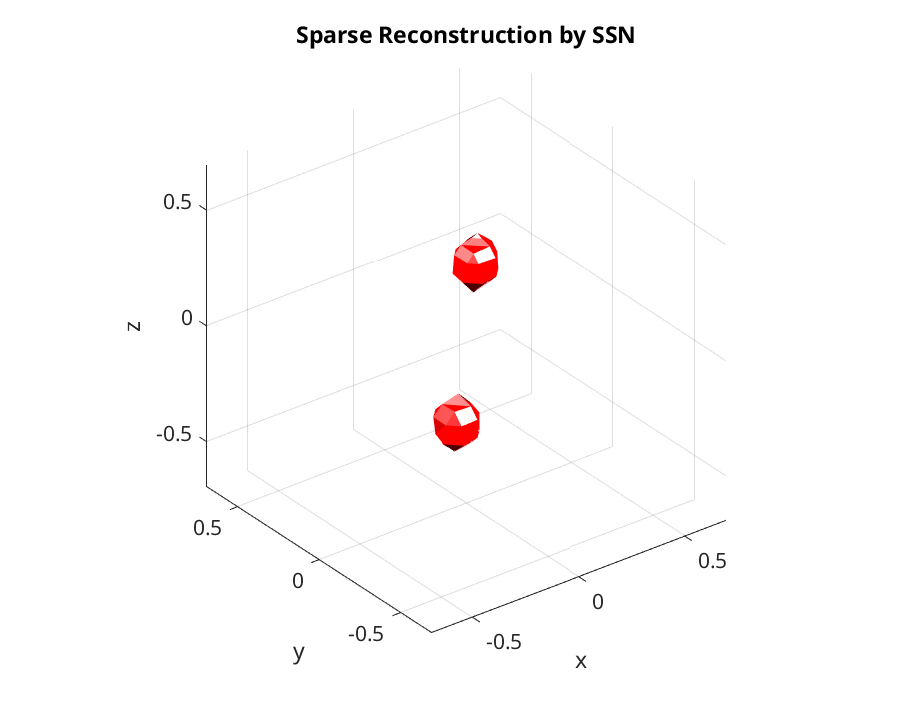}
        \caption{}
        \label{fig:1f}
    \end{subfigure}
    
    \begin{subfigure}[b]{0.3\textwidth}
        \includegraphics[width=\textwidth]{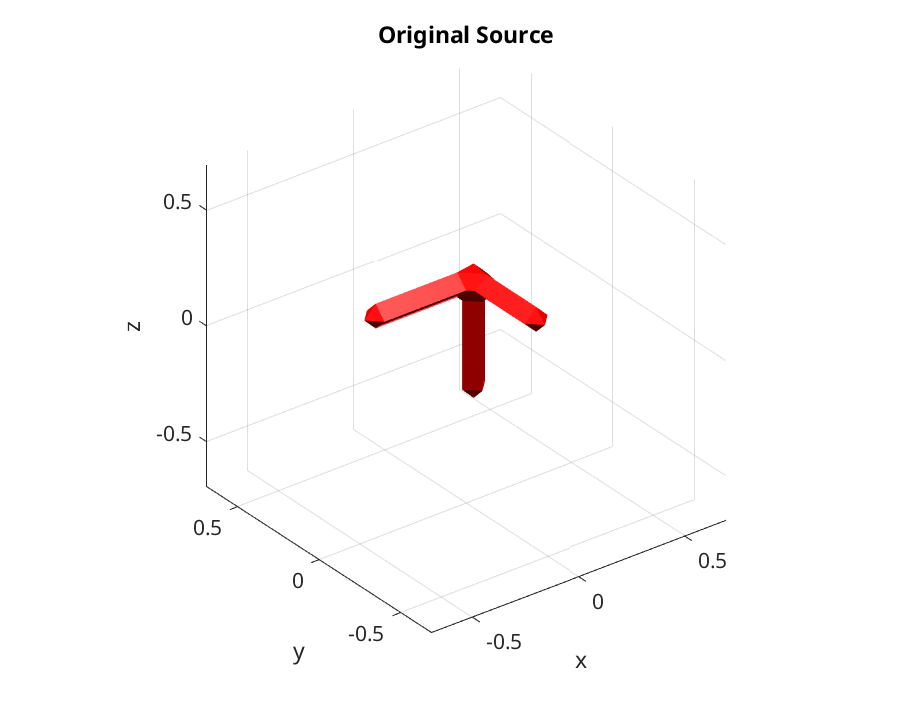}
        \caption{}
        \label{fig:1g}
    \end{subfigure}
    \hfill
    \begin{subfigure}[b]{0.3\textwidth}
        \includegraphics[width=\textwidth]{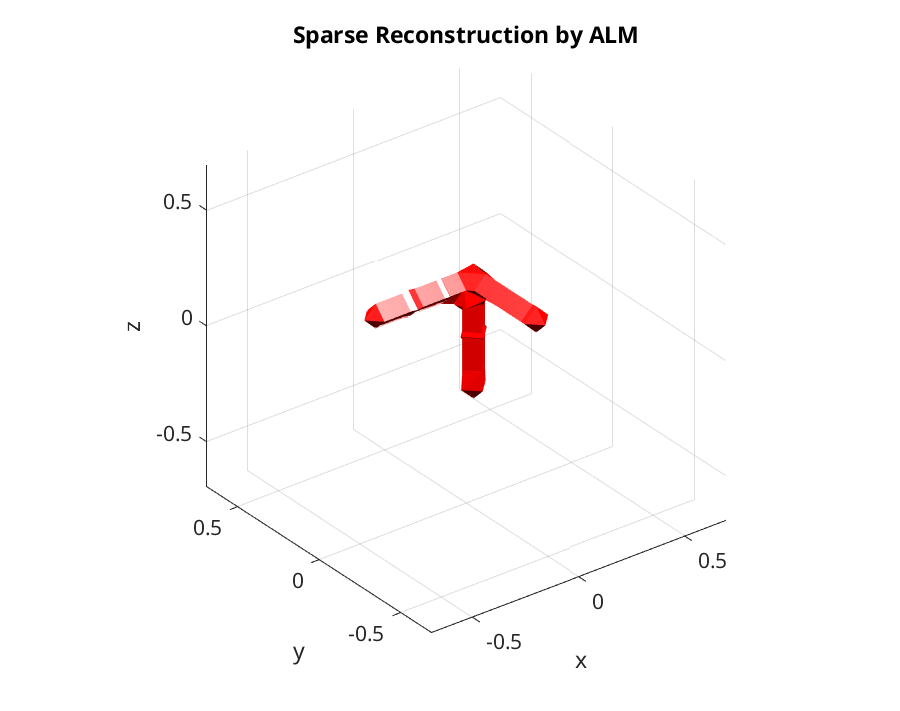}
        \caption{}
        \label{fig:1h}
    \end{subfigure}
    \hfill
    \begin{subfigure}[b]{0.3\textwidth}
        \includegraphics[width=\textwidth]{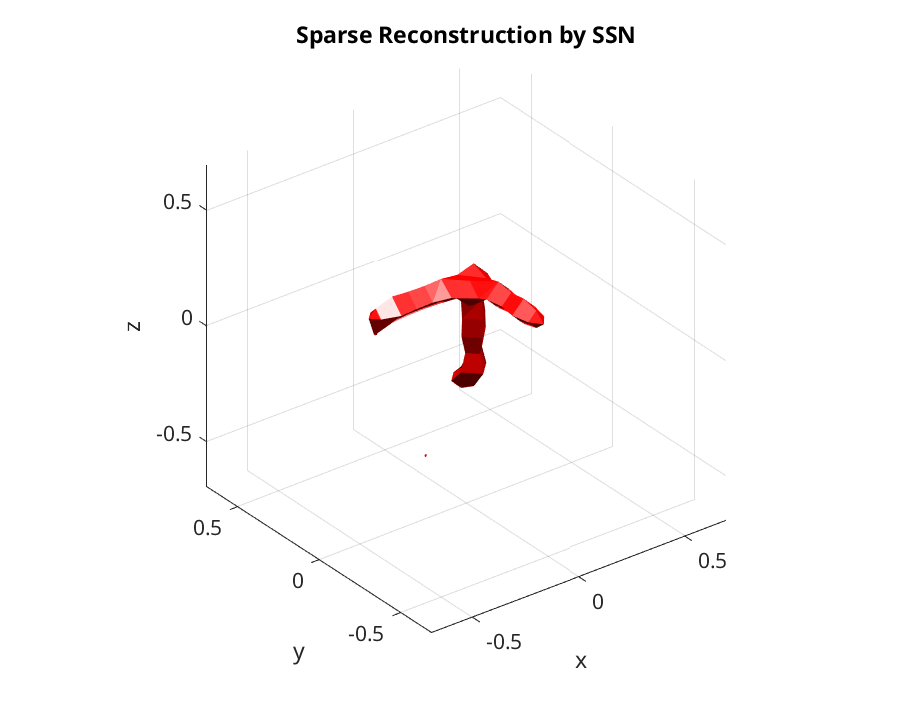}
        \caption{}
        \label{fig:1i}
    \end{subfigure}
    
    \begin{subfigure}[b]{0.3\textwidth}
        \includegraphics[width=\textwidth]{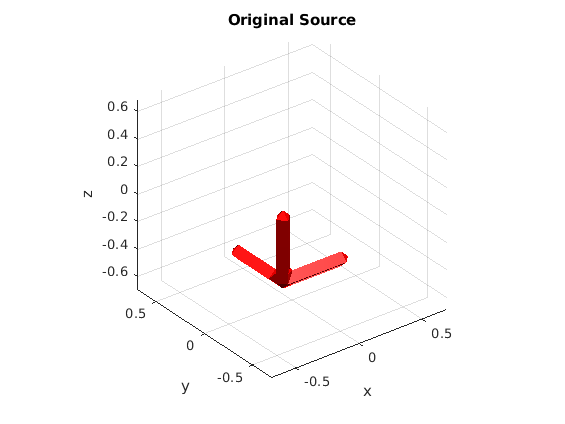}
        \caption{}
        \label{fig:1j}
    \end{subfigure}
    \hfill
    \begin{subfigure}[b]{0.3\textwidth}
        \includegraphics[width=\textwidth]{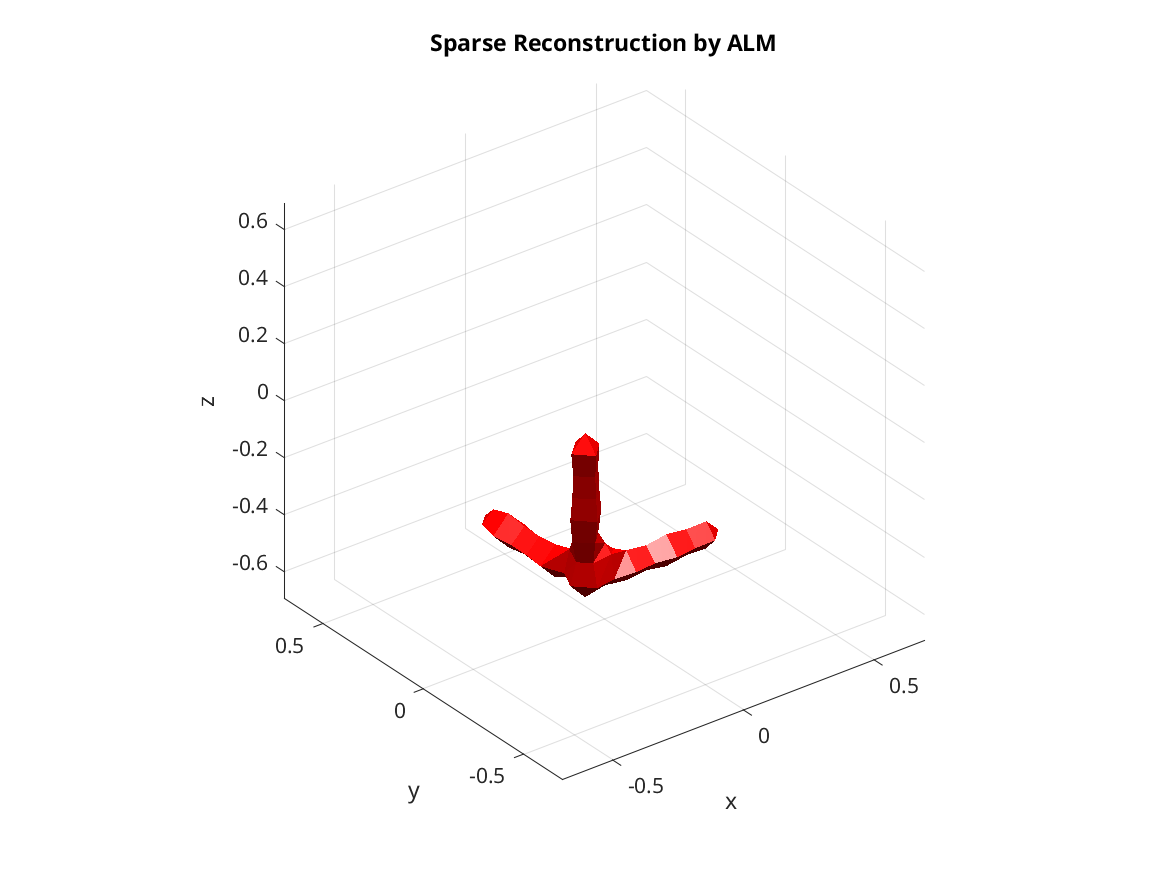}
        \caption{}
        \label{fig:1k}
    \end{subfigure}
    \hfill
    \begin{subfigure}[b]{0.3\textwidth}
        \includegraphics[width=\textwidth]{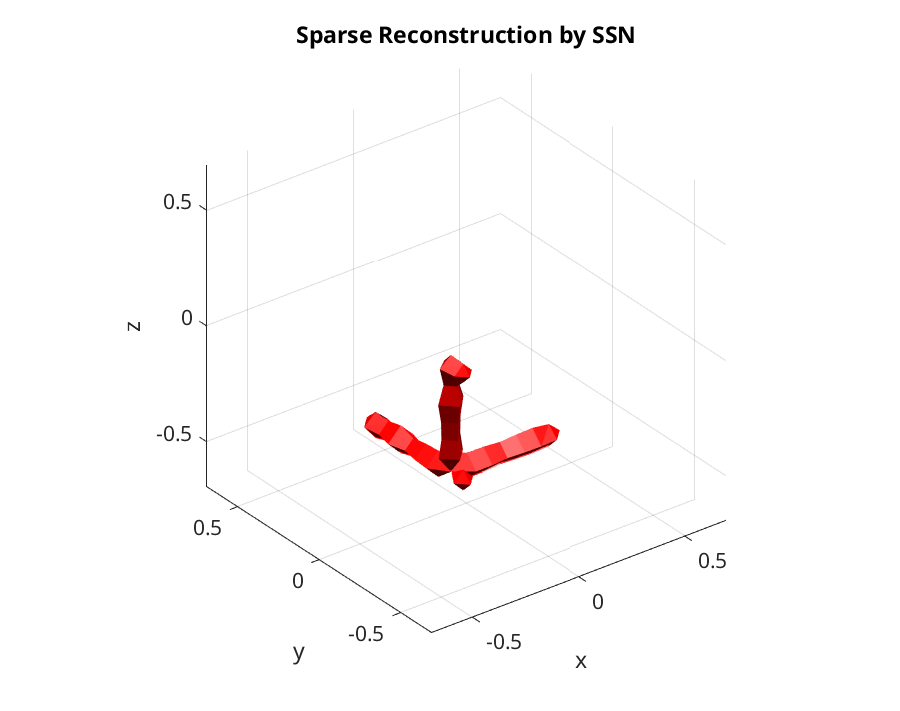}
        \caption{}
        \label{fig:1l}
    \end{subfigure}

    \begin{subfigure}[b]{0.3\textwidth}
        \includegraphics[width=\textwidth]{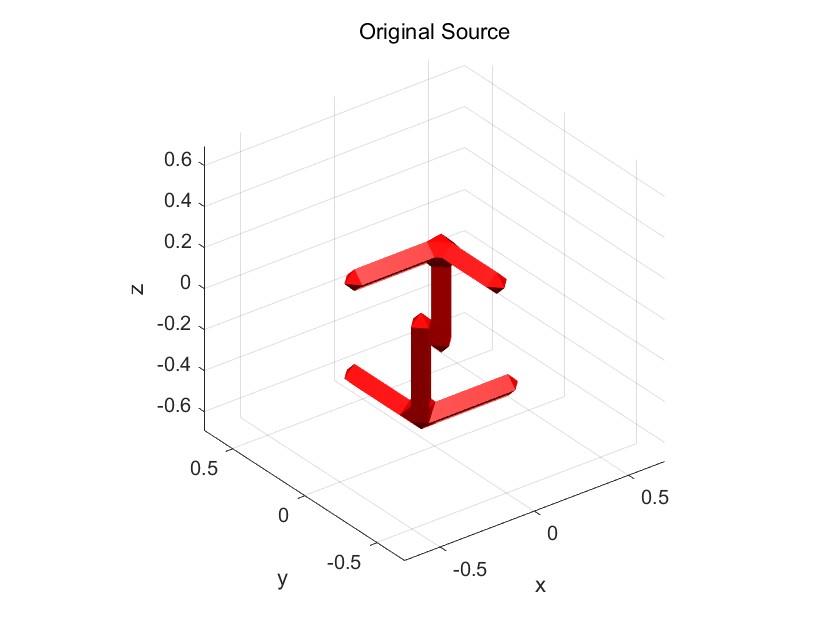}
        \caption{}
        \label{fig:1j}
    \end{subfigure}
    \hfill
    \begin{subfigure}[b]{0.3\textwidth}
        \includegraphics[width=\textwidth]{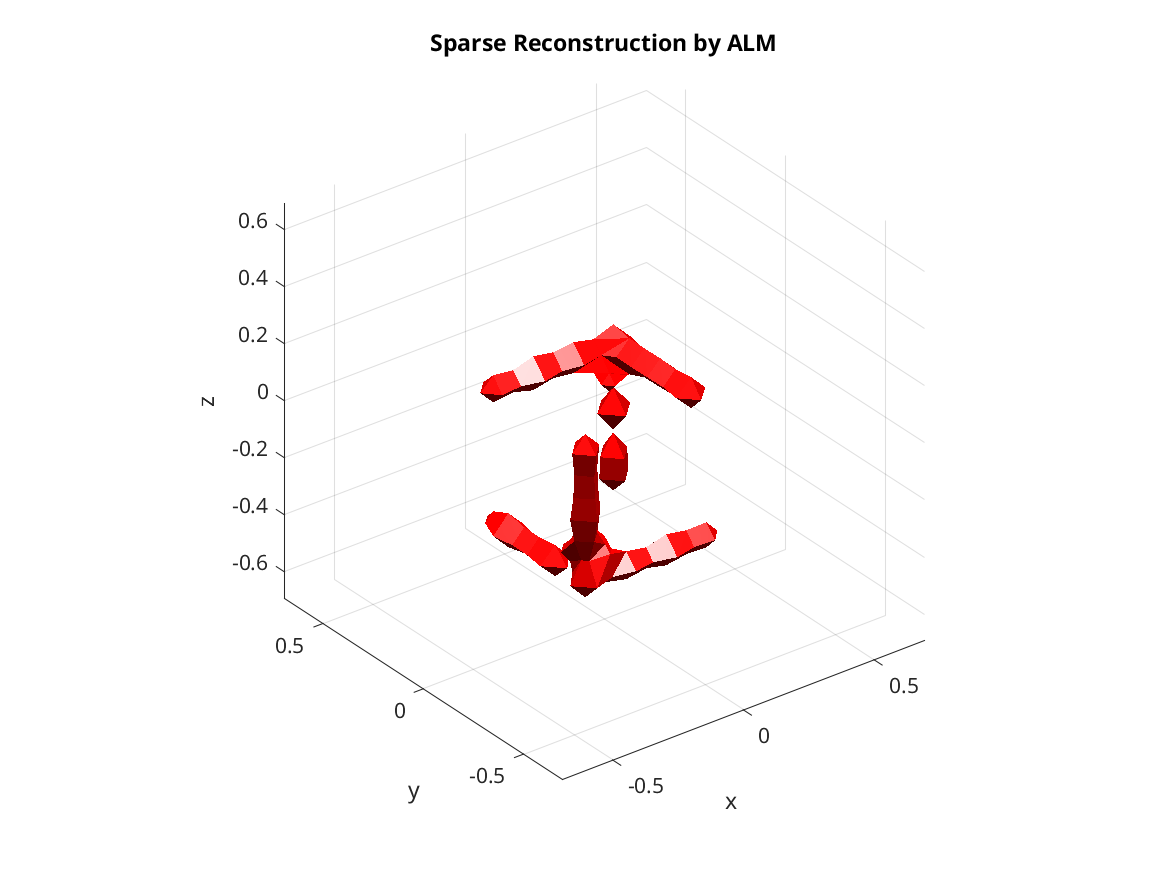}
        \caption{}
        \label{fig:1k}
    \end{subfigure}
    \hfill
    \begin{subfigure}[b]{0.3\textwidth}
        \includegraphics[width=\textwidth]{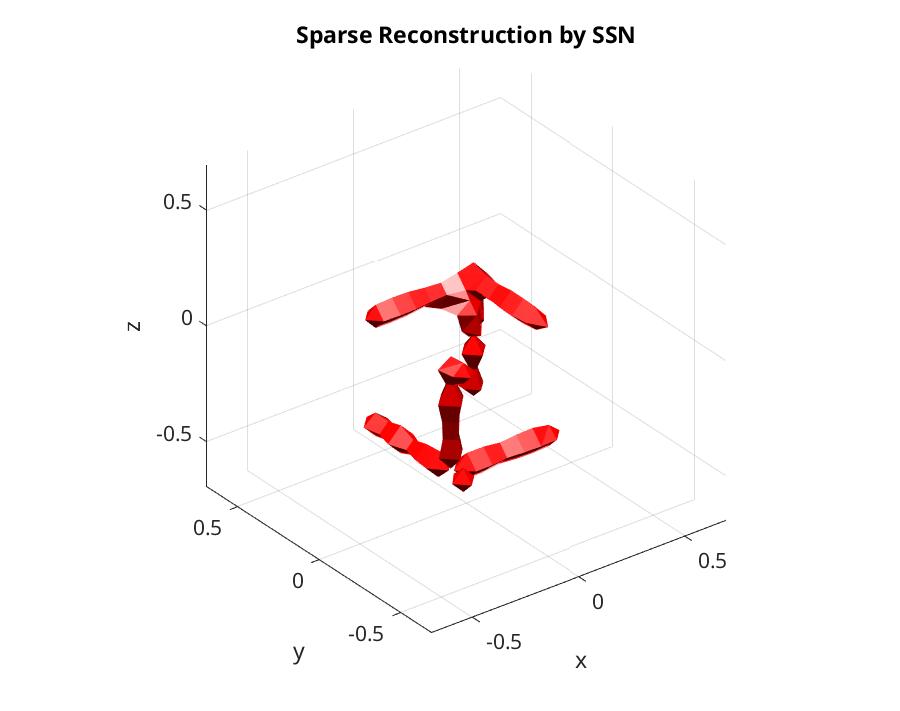}
        \caption{}
        \label{fig:1l}
    \end{subfigure}
    \caption{Reconstruction of 3-dimensional acoustic sources in homogeneous media with $k=6$, noise level 1\%, and $\alpha$ = 5e-7 together with
$\alpha_0$= 2e-8. The solution is computed by $-\lambda^{K_{\text{max}}}$.  The figures in the leftmost column are the original acoustic sources. The images in the middle and the rightmost columns are reconstructed results of ALM and SSN, respectively.  The images in the first,  second, third, and fourth rows are the sources with the two balls, the right-up tripod, the left-down tripod, and the two tripods, respectively. }
    \label{fig:big_figure:noise001:3d:homo}
\end{figure}

\begin{figure}[htbp]
    \centering

    \begin{subfigure}[b]{0.3\textwidth}
        \includegraphics[width=\textwidth]{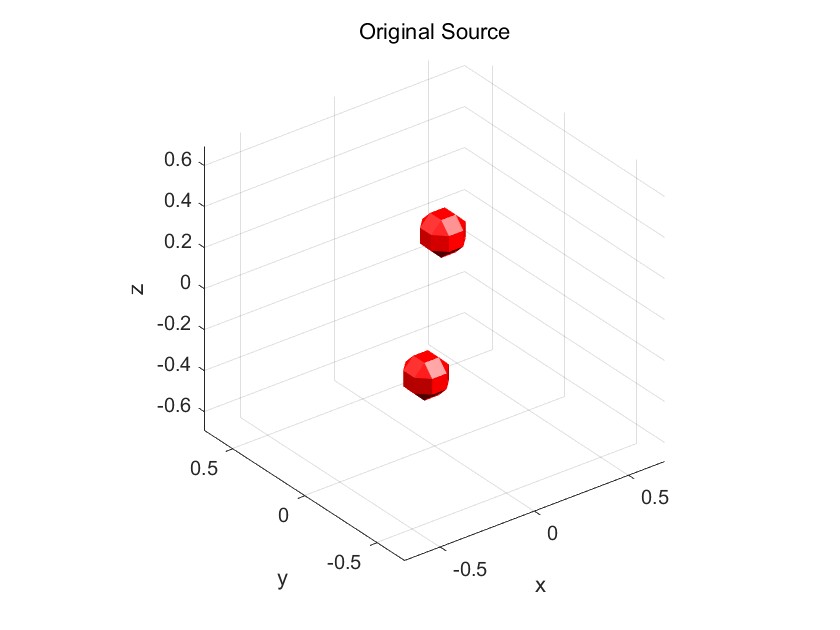}
        \caption{}
        \label{fig:1d}
    \end{subfigure}
    \hfill
    \begin{subfigure}[b]{0.3\textwidth}
        \includegraphics[width=\textwidth]{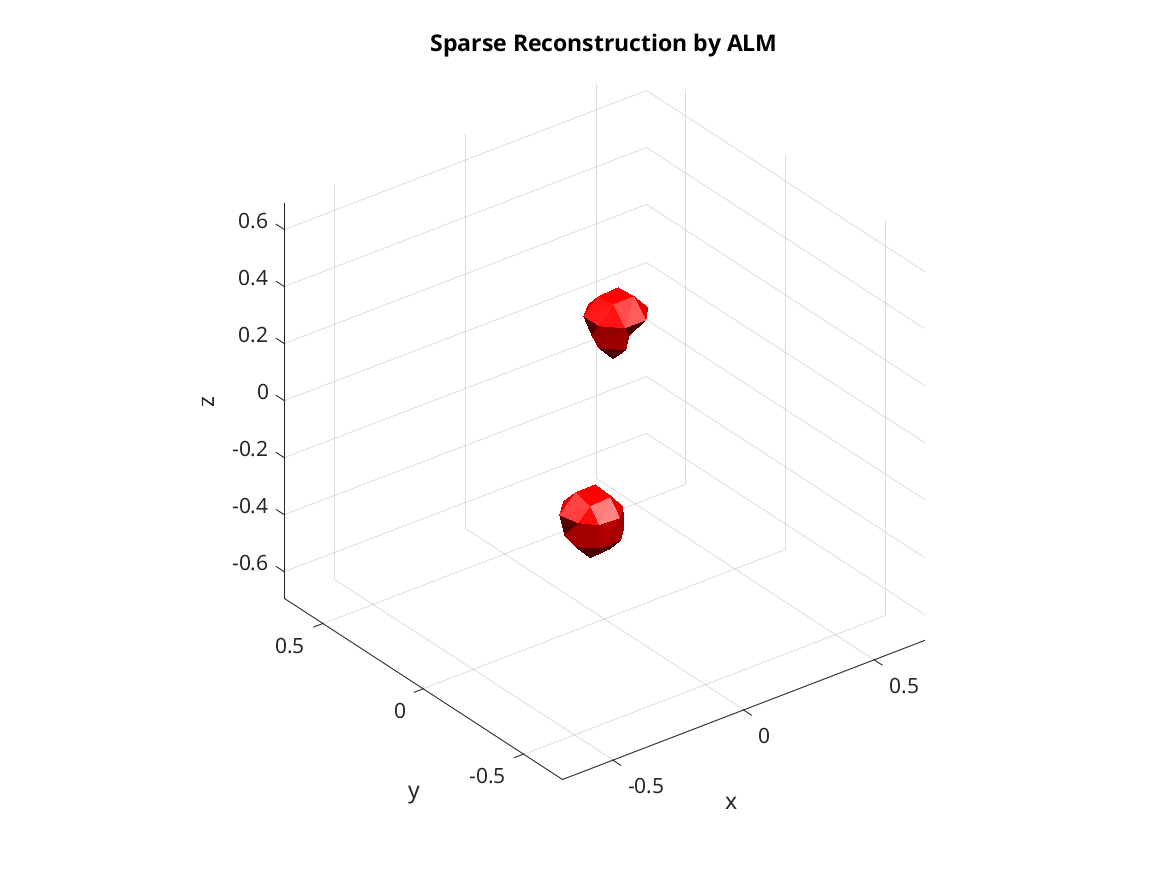}
        \caption{}
        \label{fig:1e}
    \end{subfigure}
    \hfill
    \begin{subfigure}[b]{0.3\textwidth}
        \includegraphics[width=\textwidth]{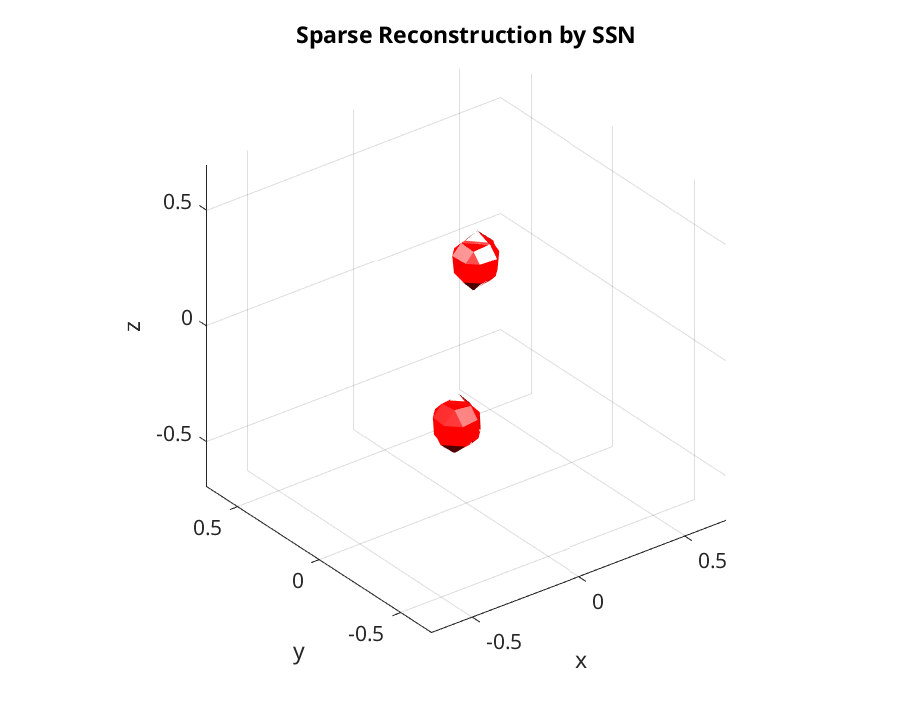}
        \caption{}
        \label{fig:1f}
    \end{subfigure}
    
    \begin{subfigure}[b]{0.3\textwidth}
        \includegraphics[width=\textwidth]{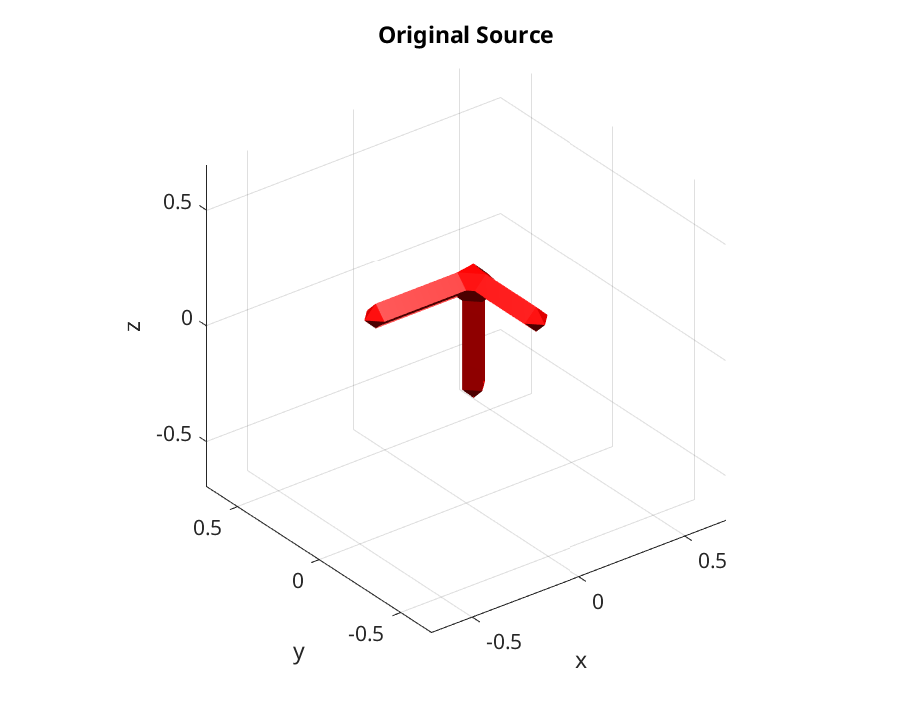}
        \caption{}
        \label{fig:1g}
    \end{subfigure}
    \hfill
    \begin{subfigure}[b]{0.3\textwidth}
        \includegraphics[width=\textwidth]{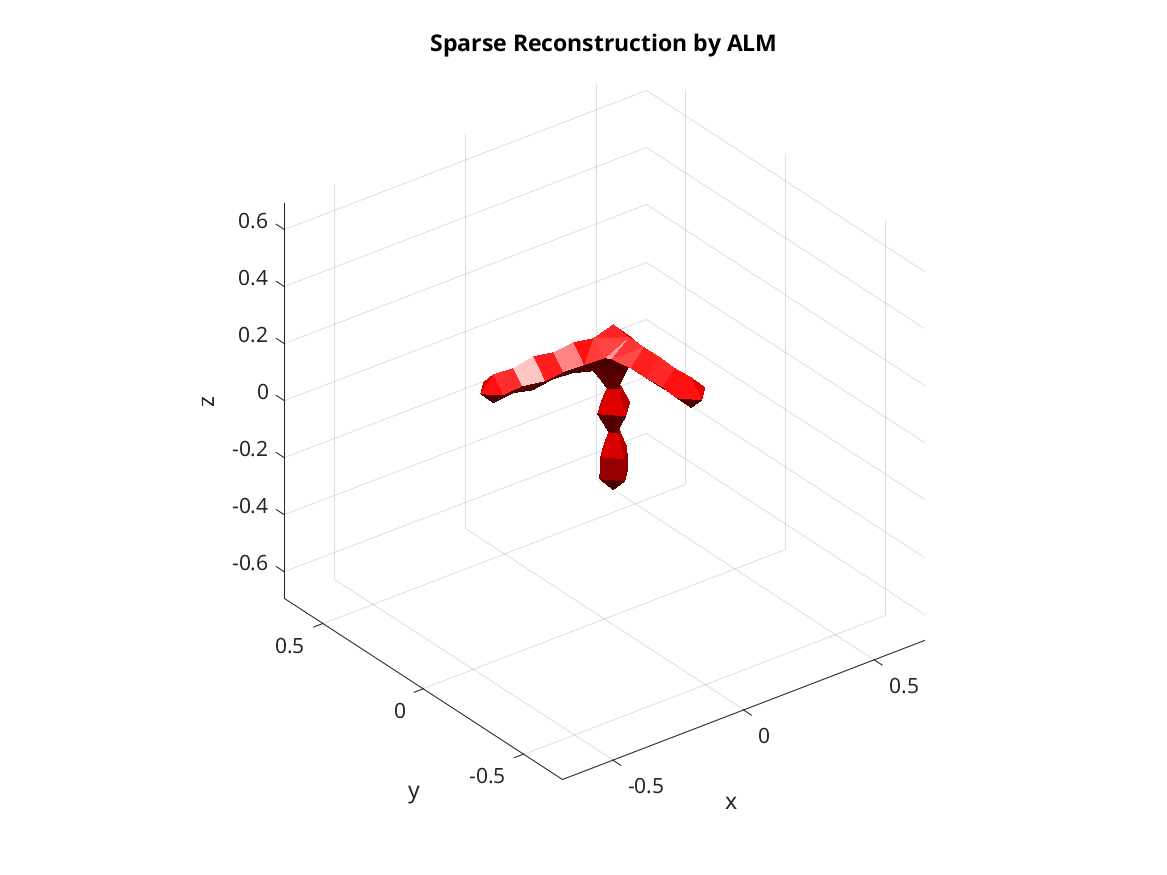}
        \caption{}
        \label{fig:1h}
    \end{subfigure}
    \hfill
    \begin{subfigure}[b]{0.3\textwidth}
        \includegraphics[width=\textwidth]{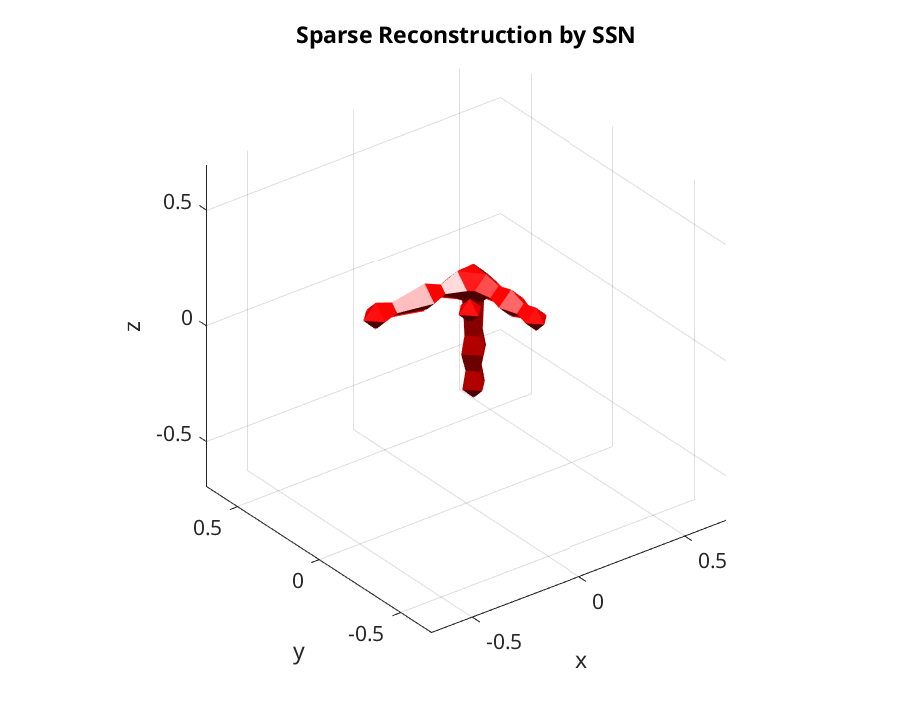}
        \caption{}
        \label{fig:1i}
    \end{subfigure}
    
    \begin{subfigure}[b]{0.3\textwidth}
        \includegraphics[width=\textwidth]{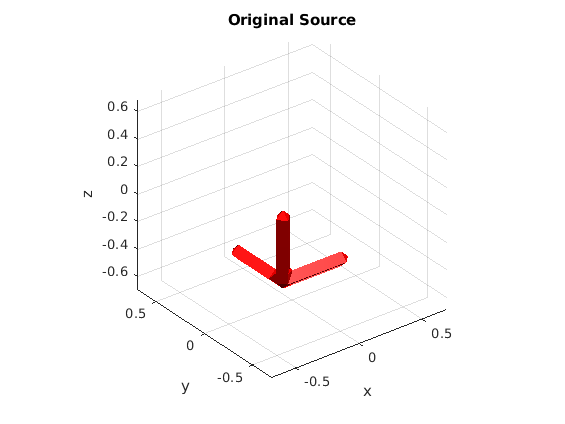}
        \caption{}
        \label{fig:1j}
    \end{subfigure}
    \hfill
    \begin{subfigure}[b]{0.3\textwidth}
        \includegraphics[width=\textwidth]{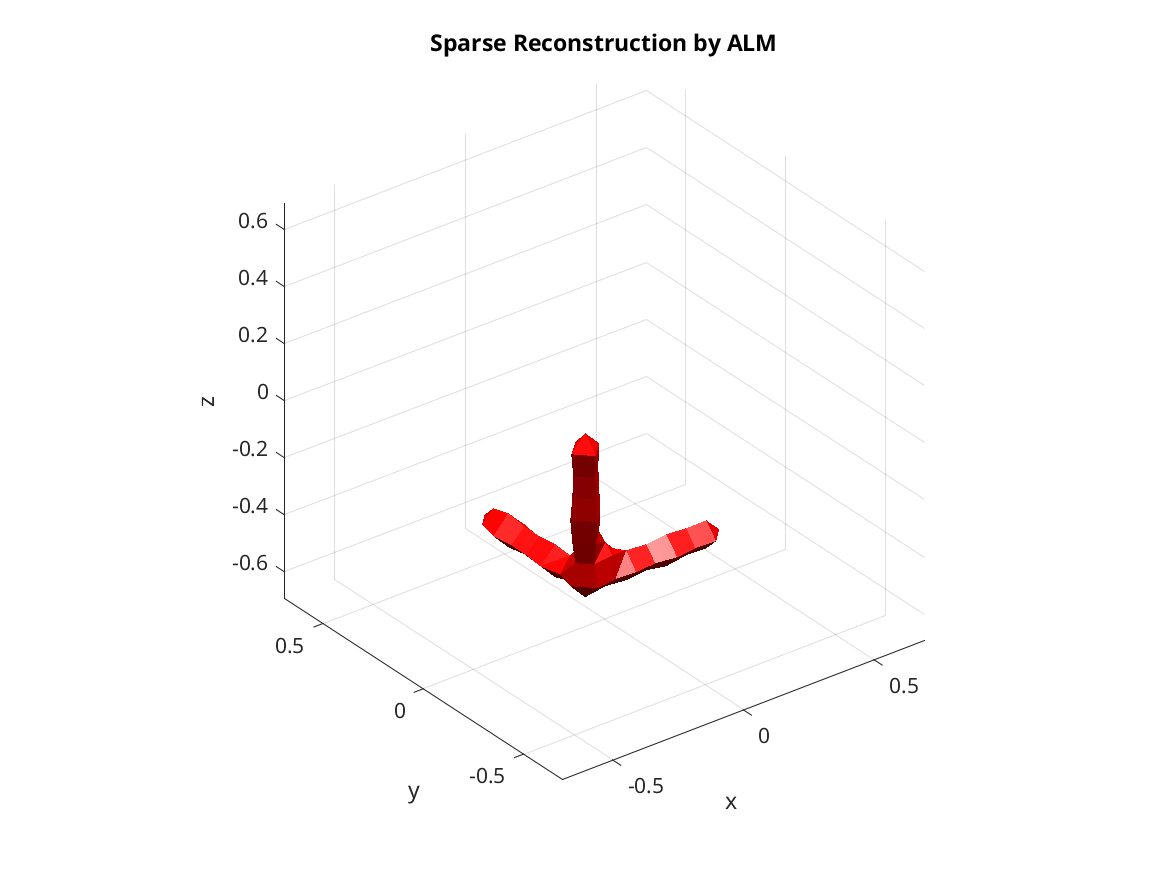}
        \caption{}
        \label{fig:1k}
    \end{subfigure}
    \hfill
    \begin{subfigure}[b]{0.3\textwidth}
        \includegraphics[width=\textwidth]{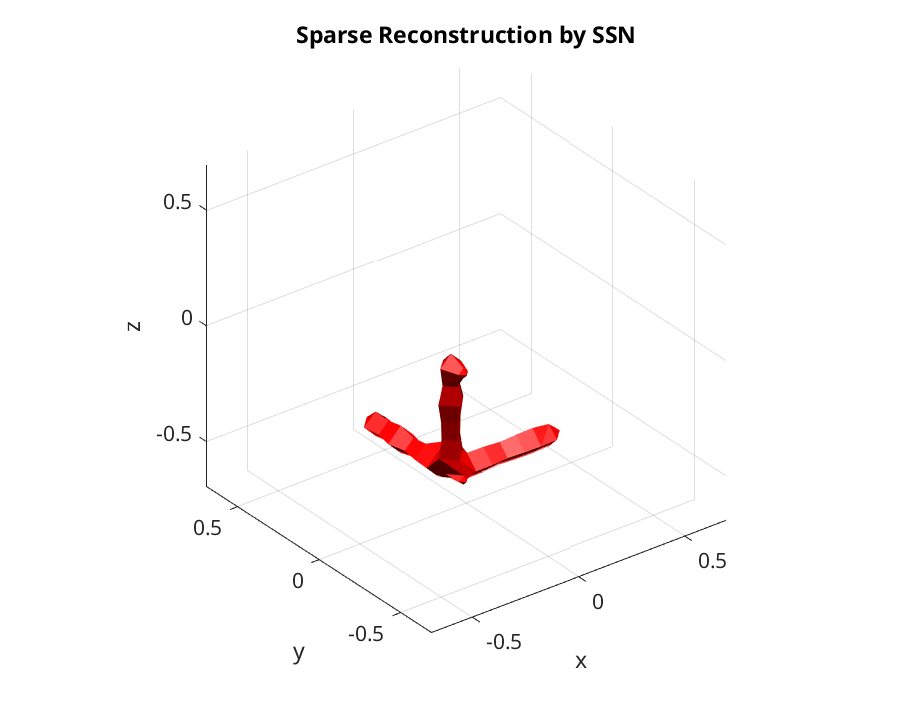}
        \caption{}
        \label{fig:1l}
    \end{subfigure}

    \begin{subfigure}[b]{0.3\textwidth}
        \includegraphics[width=\textwidth]{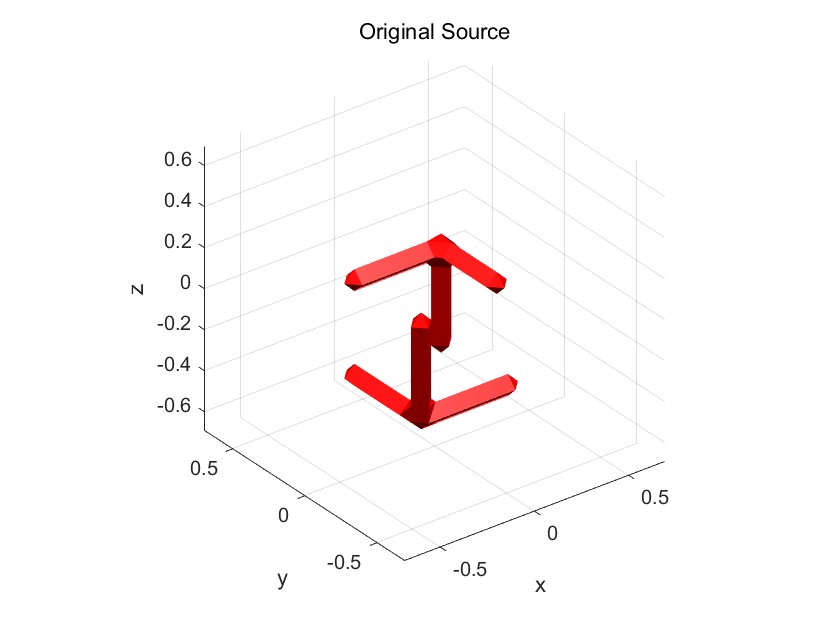}
        \caption{}
        \label{fig:1j}
    \end{subfigure}
    \hfill
    \begin{subfigure}[b]{0.3\textwidth}
        \includegraphics[width=\textwidth]{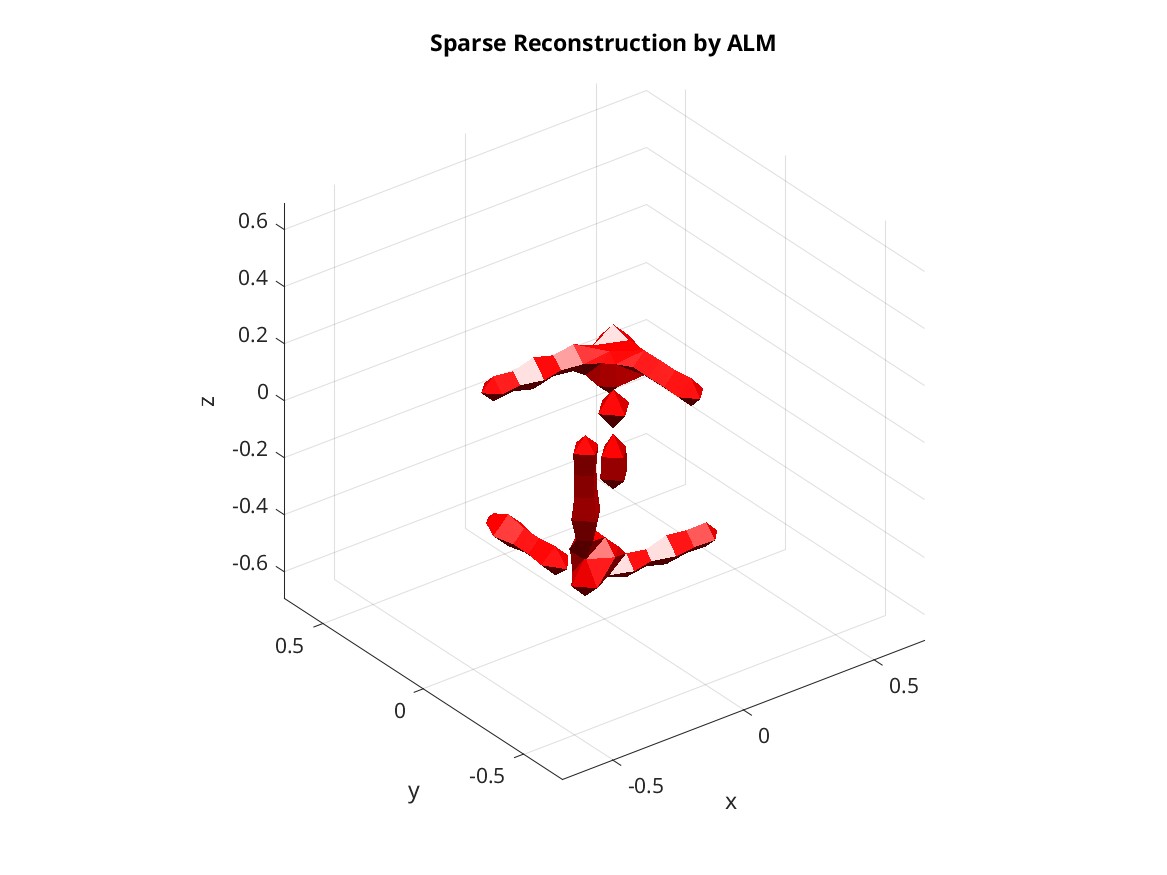}
        \caption{}
        \label{fig:1k}
    \end{subfigure}
    \hfill
    \begin{subfigure}[b]{0.3\textwidth}
        \includegraphics[width=\textwidth]{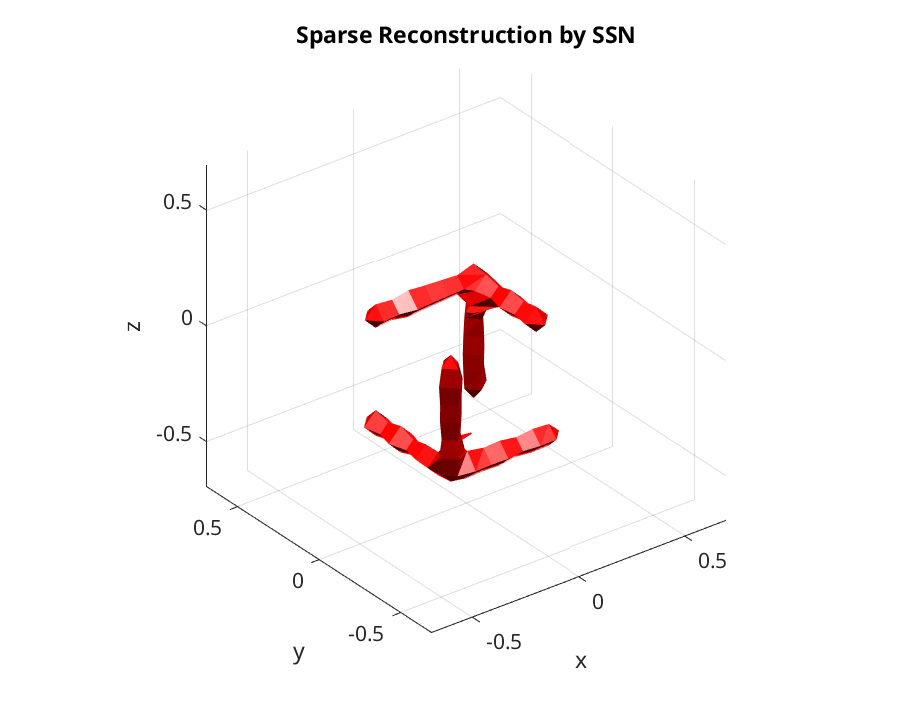}
        \caption{}
        \label{fig:1l}
    \end{subfigure}
    \caption{Reconstruction of 3-dimensional acoustic sources in inhomogeneous media with $k=6$ and  noise level 1\%, and $\alpha$ = 5e-7 together with
$\alpha_0$= 2e-8. The solution is computed by $-\lambda^{K_{\text{max}}}$. The information of the acoustic sources and the corresponding reconstruction algorithms are the same as in Figure \ref{fig:big_figure:noise001:3d:homo}.}
    \label{fig:big_figure:inhomo:3d:001}
\end{figure}

We computed the primal solution $\mu^{K_{\text{max}}}$ by \eqref{eq:mu} except the numerical examples in  Figures \ref{fig:big_figure:noise001:3d:homo} and \ref{fig:big_figure:inhomo:3d:001}, which are computed by $-\lambda^{K_{\text{max}}}$ due to \eqref{eq"lambda} and the explanations there.

Now, let us turn to the numerical examples, and begin with the single-frequency cases.

\begin{figure}
     
    \begin{subfigure}[b]{0.2\textwidth}
        \includegraphics[width=\textwidth]{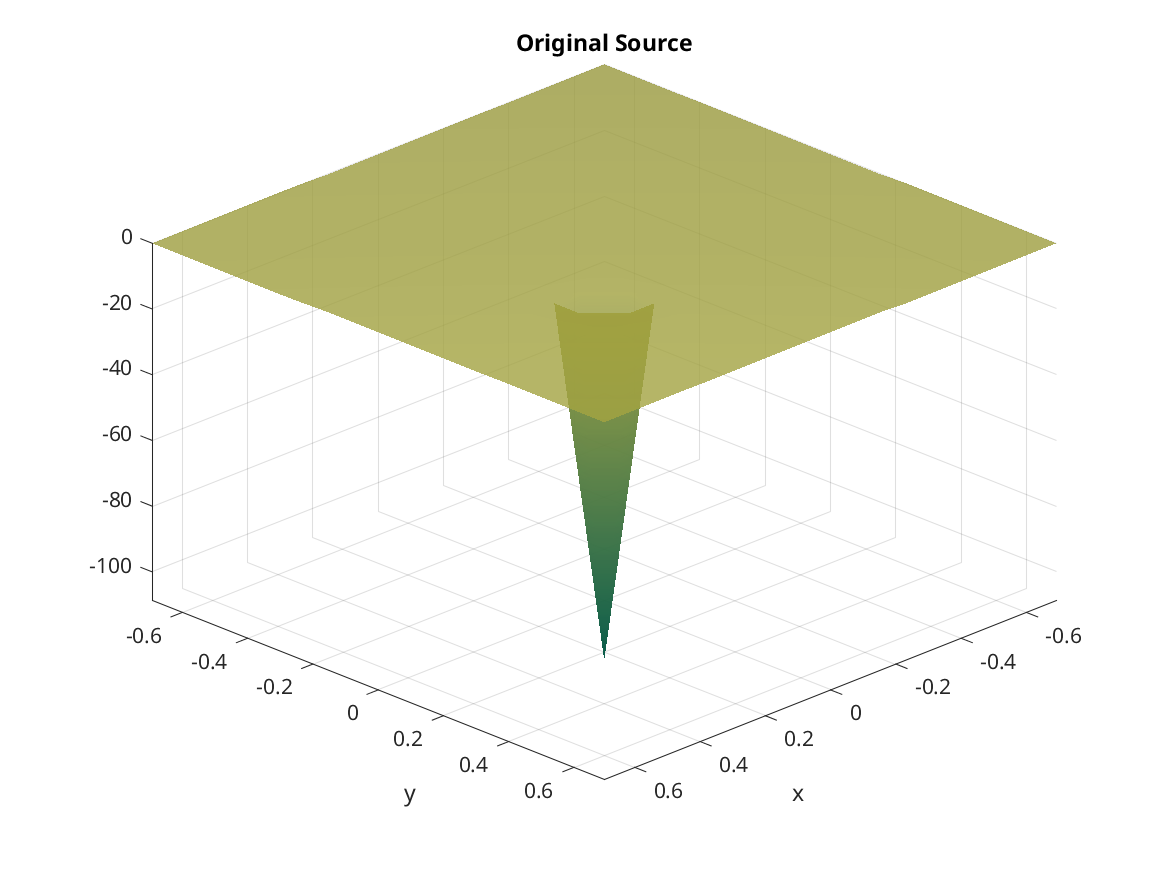}
        \caption{}
        \label{}
    \end{subfigure}
    \hfill
    \begin{subfigure}[b]{0.2\textwidth}
        \includegraphics[width=\textwidth]{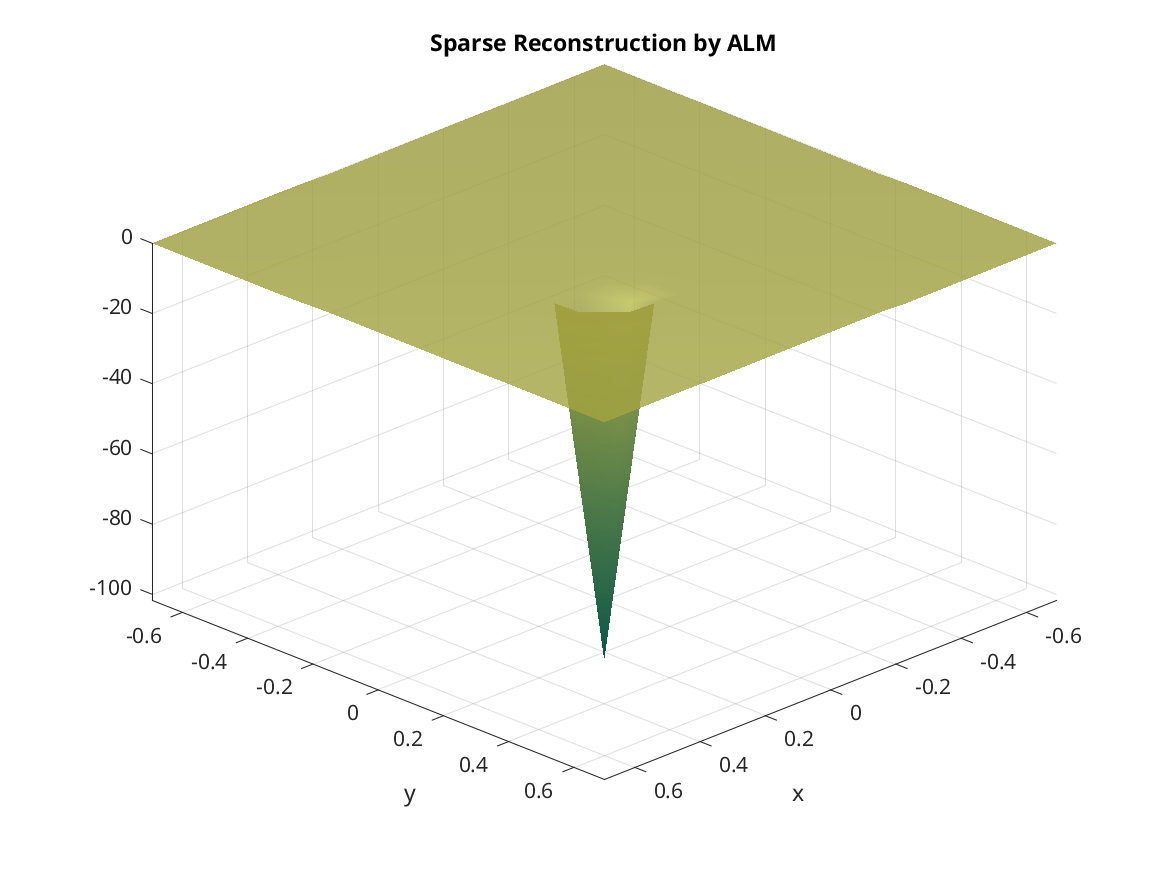}
        \caption{}
        \label{}
    \end{subfigure}
    \hfill
    \begin{subfigure}[b]{0.2\textwidth}
        \includegraphics[width=\textwidth]{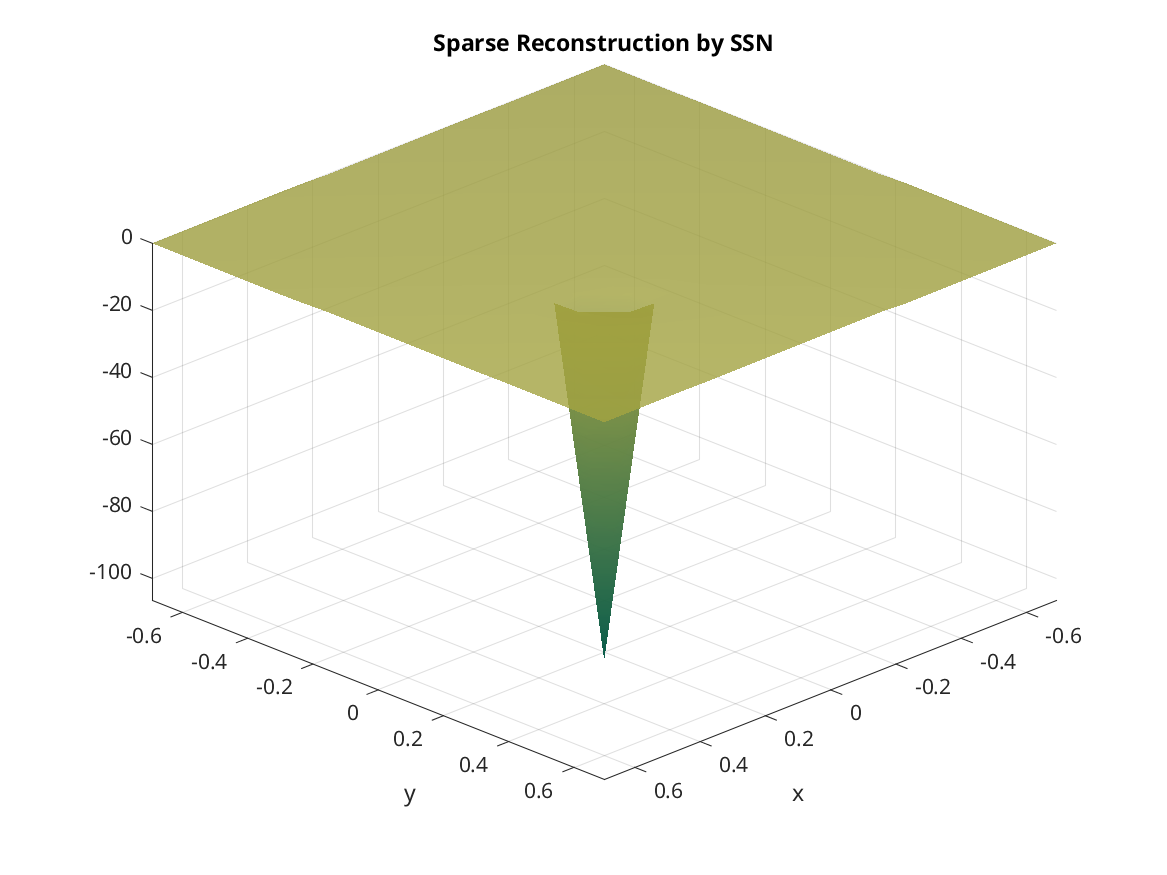}
        \caption{}
        \label{}
    \end{subfigure}
    \hfill
    \begin{subfigure}[b]{0.2\textwidth}
        \includegraphics[width=\textwidth]{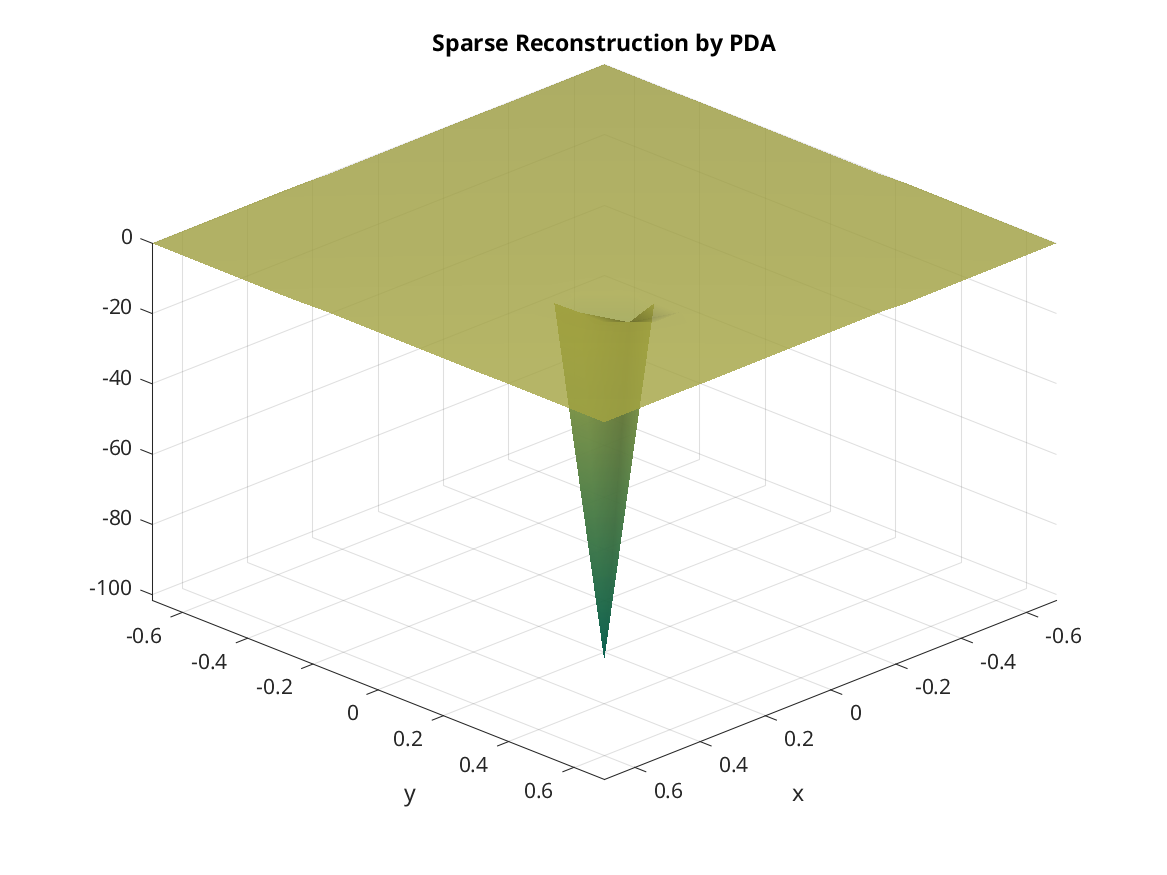}
        \caption{}
        \label{}
    \end{subfigure}
    
    \begin{subfigure}[b]{0.2\textwidth}
        \includegraphics[width=\textwidth]{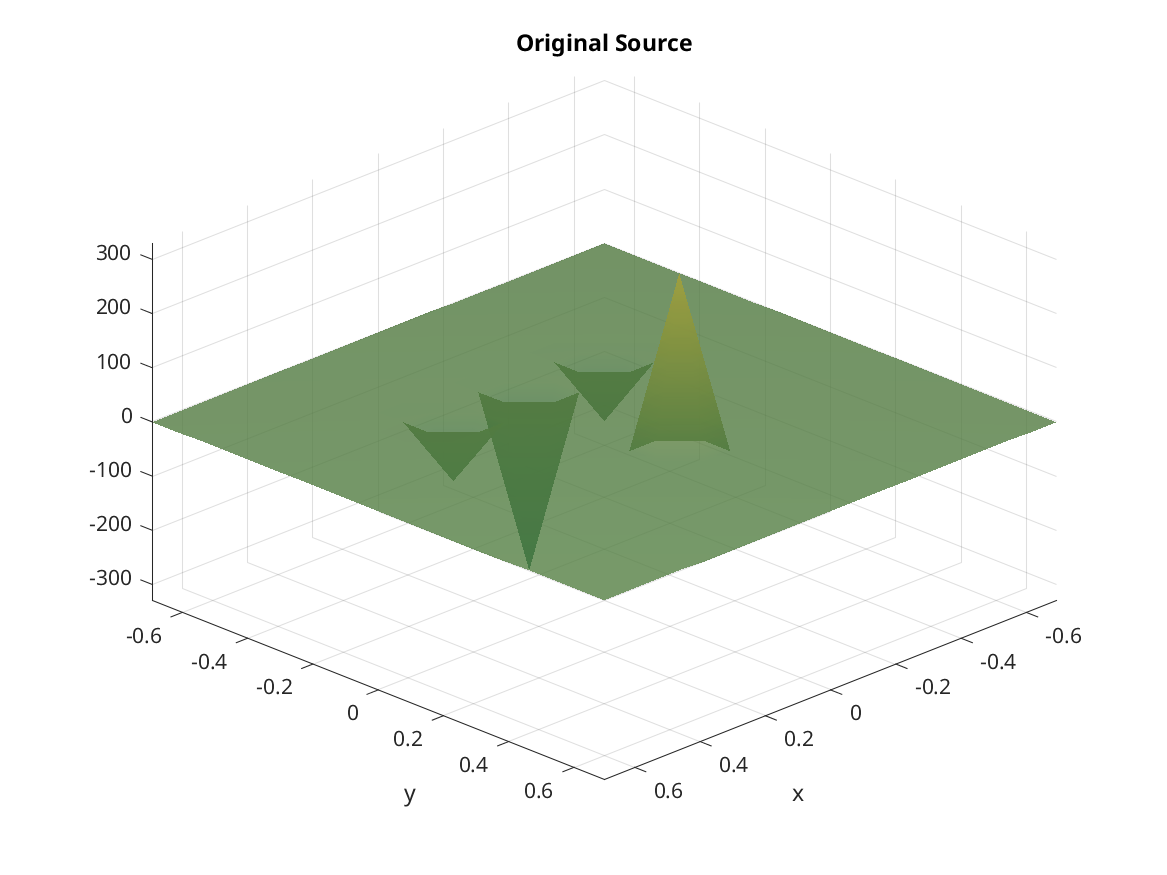}
        \caption{}
        \label{}
    \end{subfigure}
    \hfill
    \begin{subfigure}[b]{0.2\textwidth}
        \includegraphics[width=\textwidth]{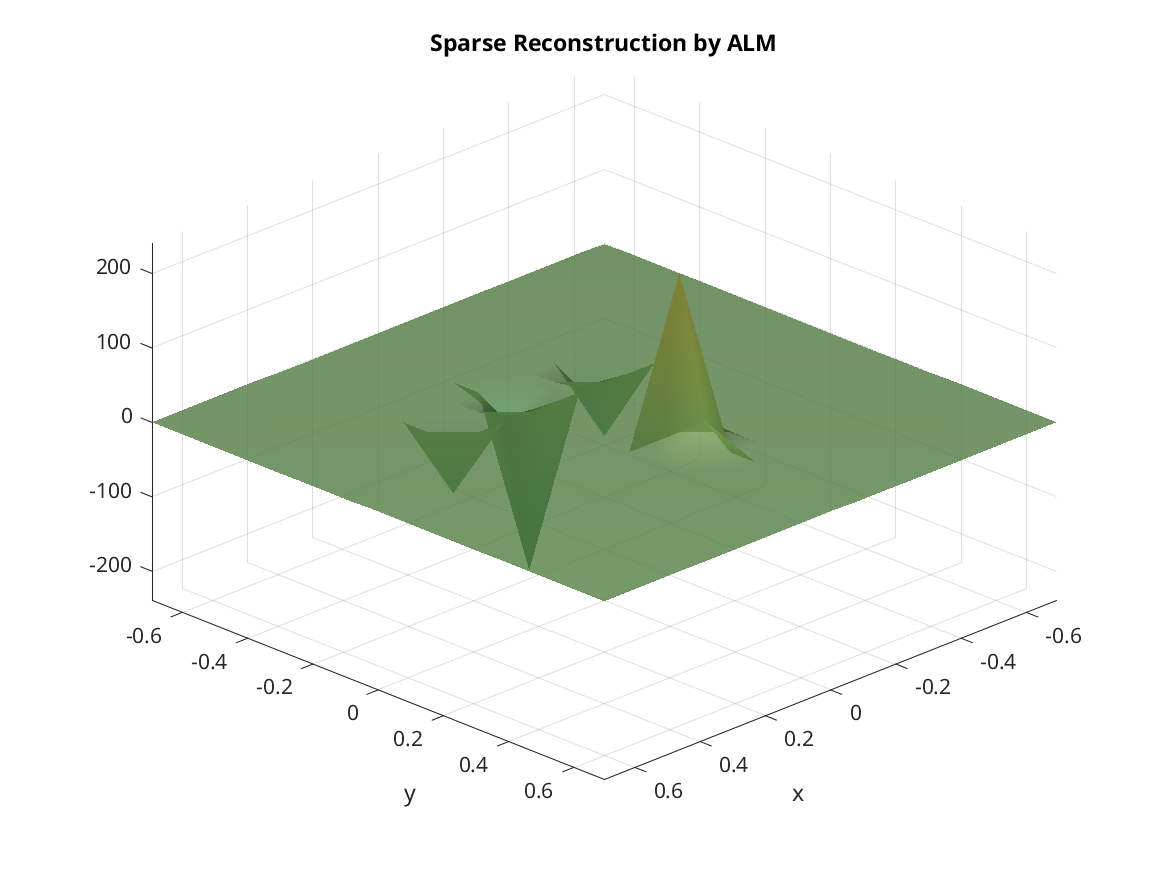}
        \caption{}
        \label{}
    \end{subfigure}
    \hfill
    \begin{subfigure}[b]{0.2\textwidth}
        \includegraphics[width=\textwidth]{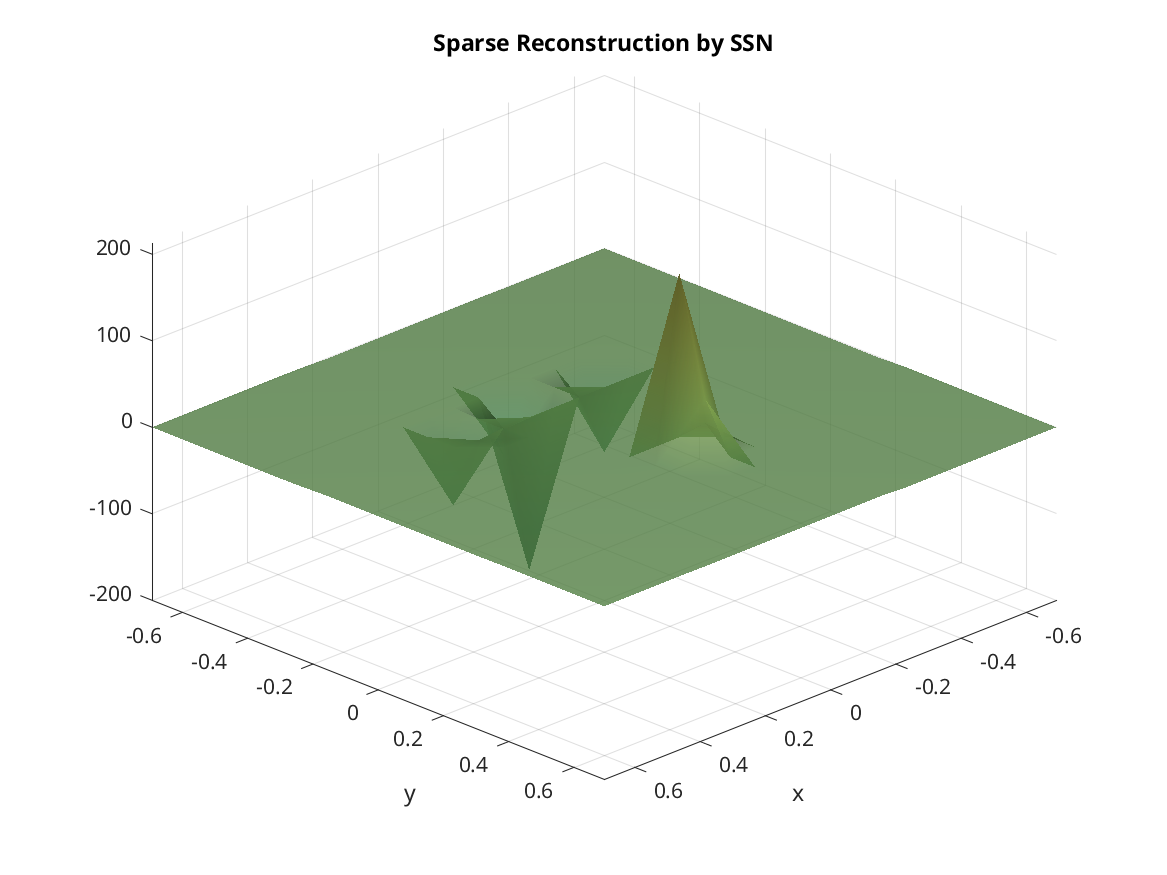}
        \caption{}
        \label{}
    \end{subfigure}
    \hfill
    \begin{subfigure}[b]{0.2\textwidth}
        \includegraphics[width=\textwidth]{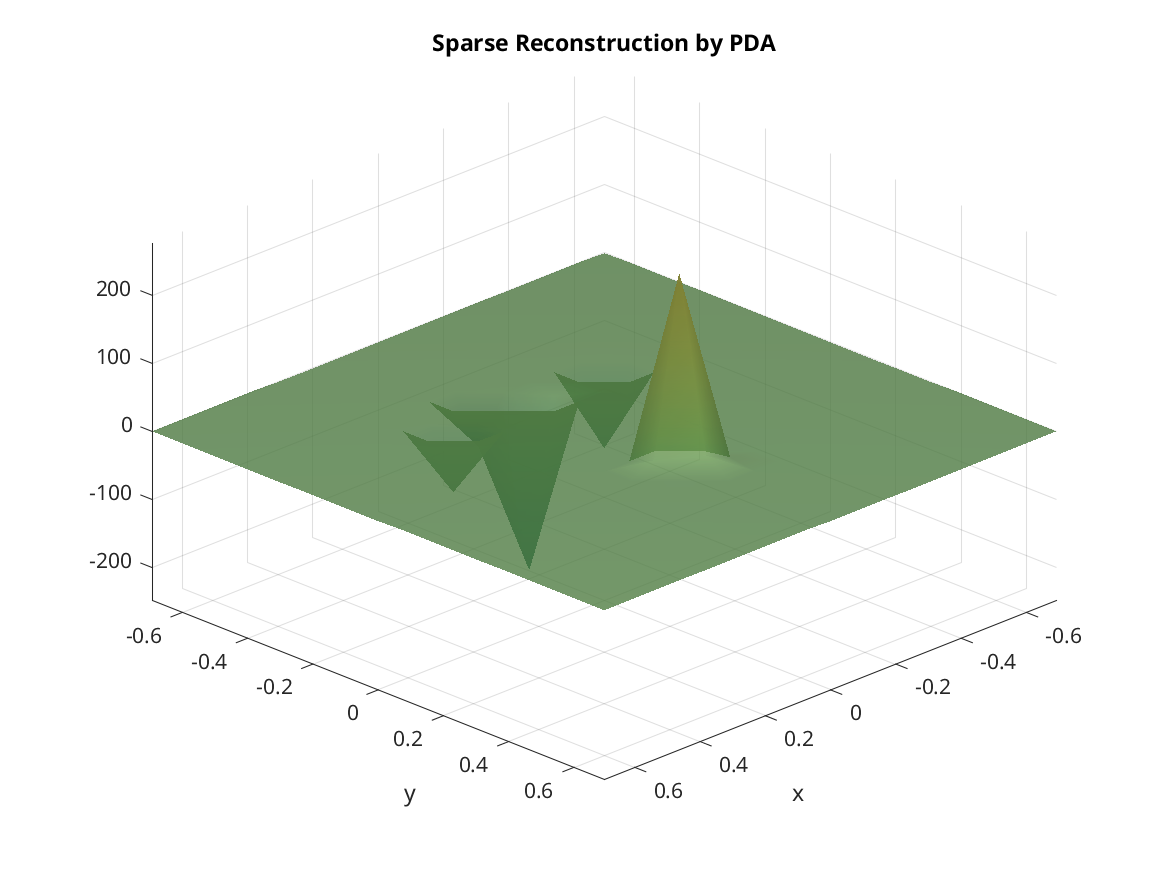}
        \caption{}
        \label{}
    \end{subfigure}
    
    \begin{subfigure}[b]{0.2\textwidth}
        \includegraphics[width=\textwidth]{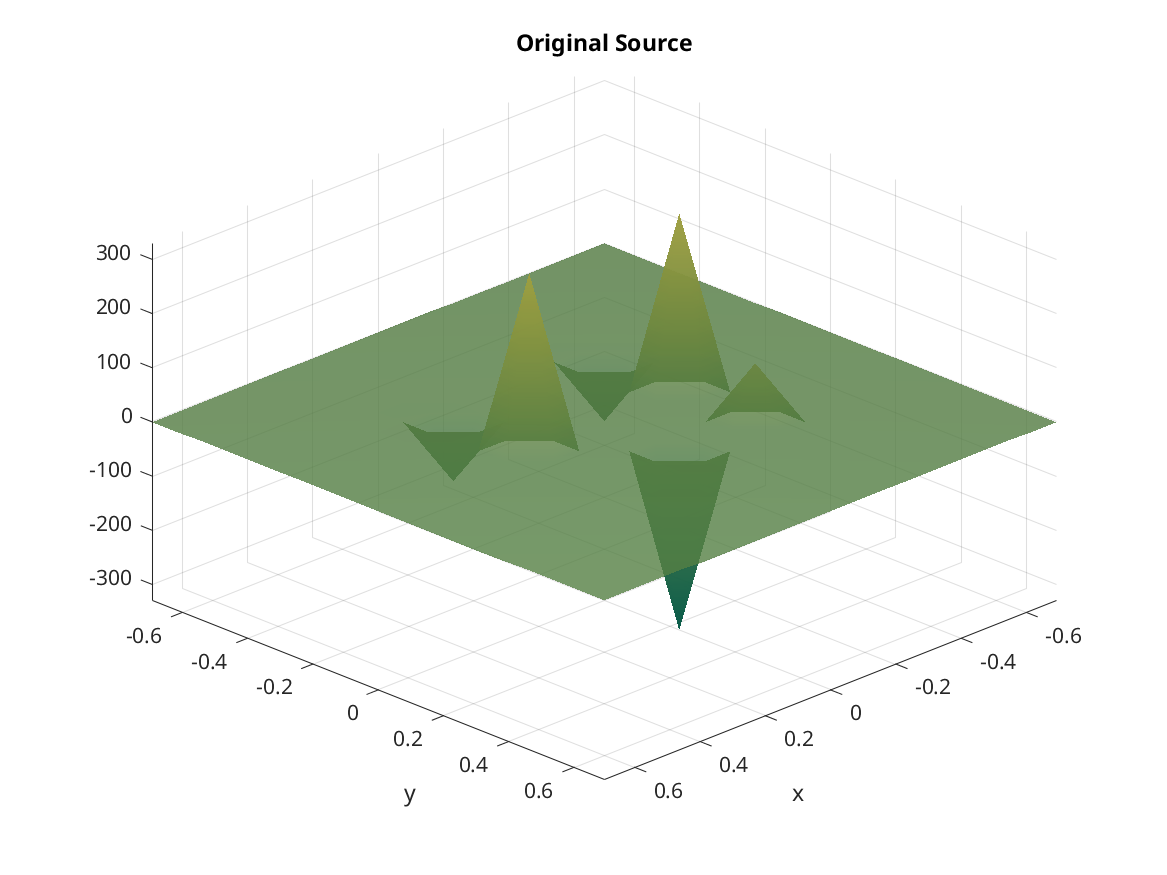}
        \caption{}
        \label{}
    \end{subfigure}
    \hfill
    \begin{subfigure}[b]{0.2\textwidth}
        \includegraphics[width=\textwidth]{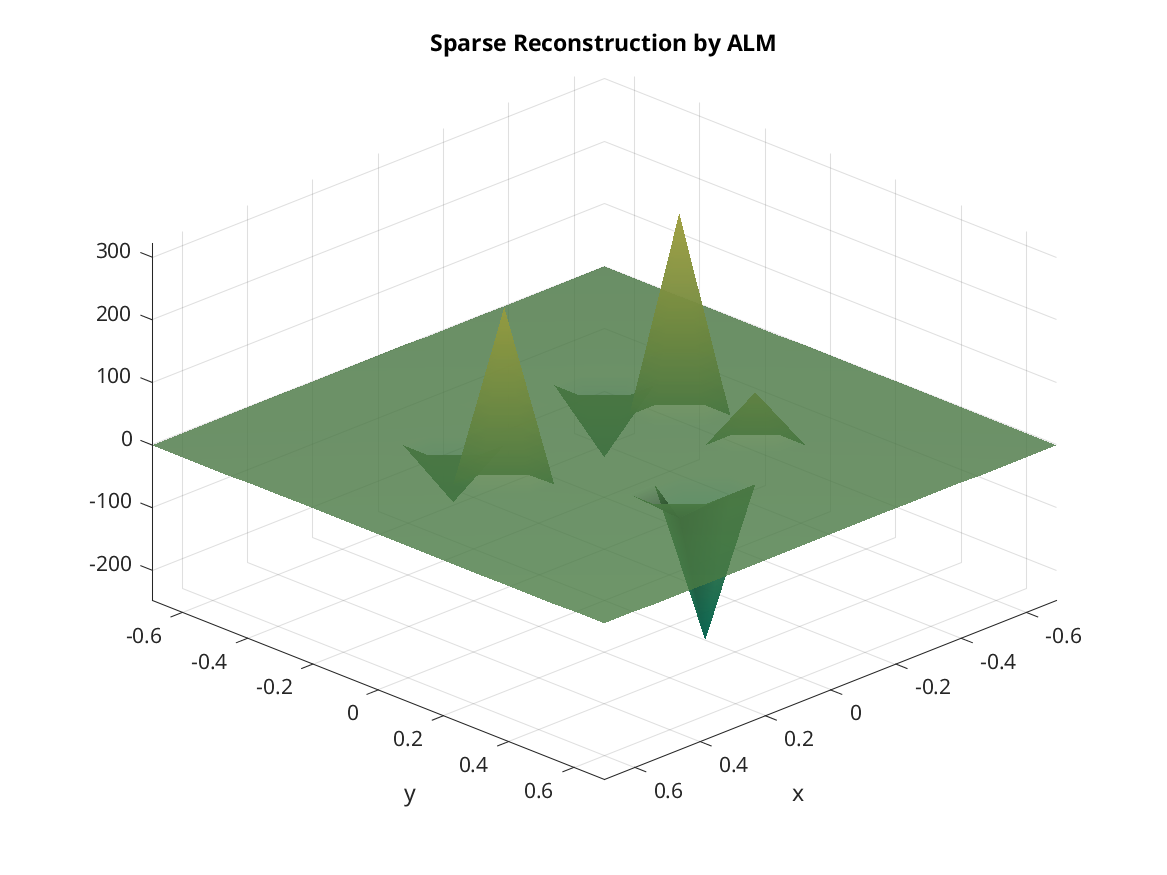}
        \caption{}
        \label{}
    \end{subfigure}
    \hfill
    \begin{subfigure}[b]{0.2\textwidth}
        \includegraphics[width=\textwidth]{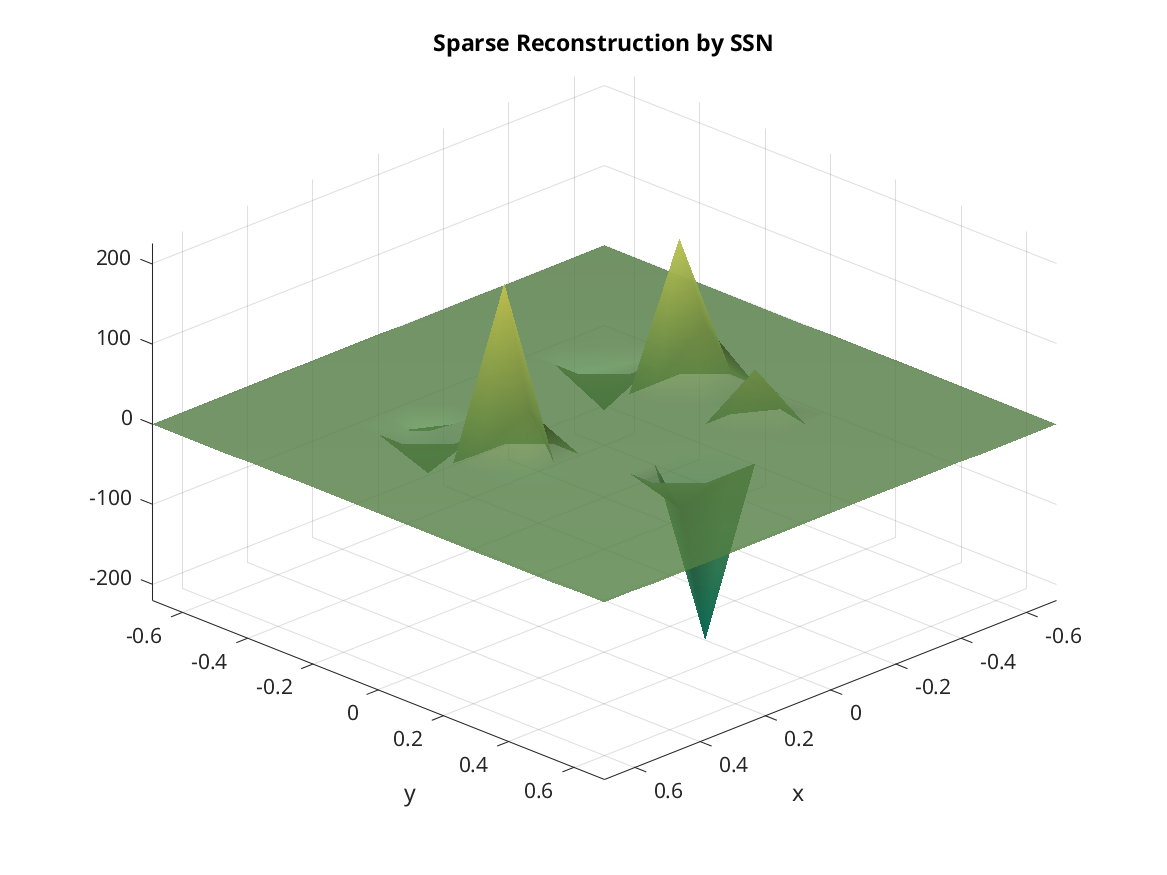}
        \caption{}
        \label{}
    \end{subfigure}
    \hfill
    \begin{subfigure}[b]{0.2\textwidth}
        \includegraphics[width=\textwidth]{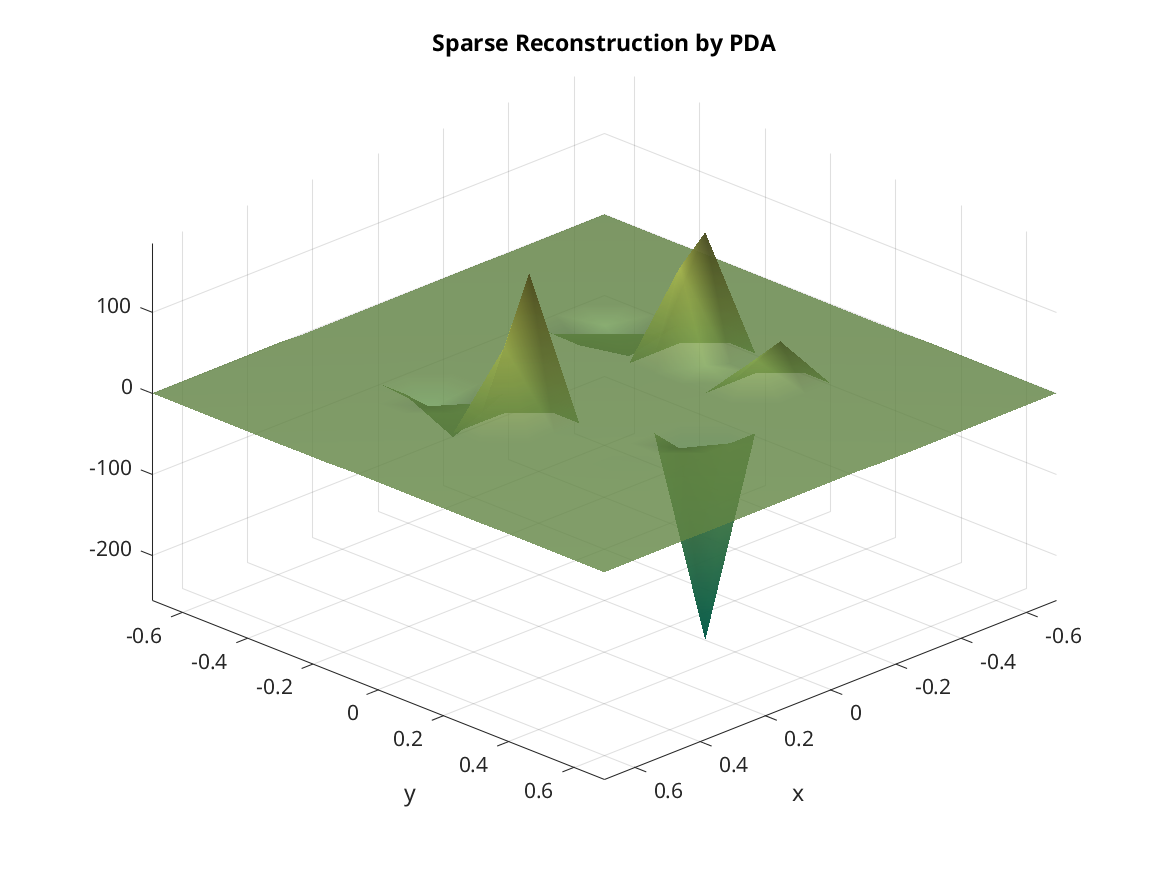}
        \caption{}
        \label{}
    \end{subfigure}
    
    \begin{subfigure}[b]{0.2\textwidth}
        \includegraphics[width=\textwidth]{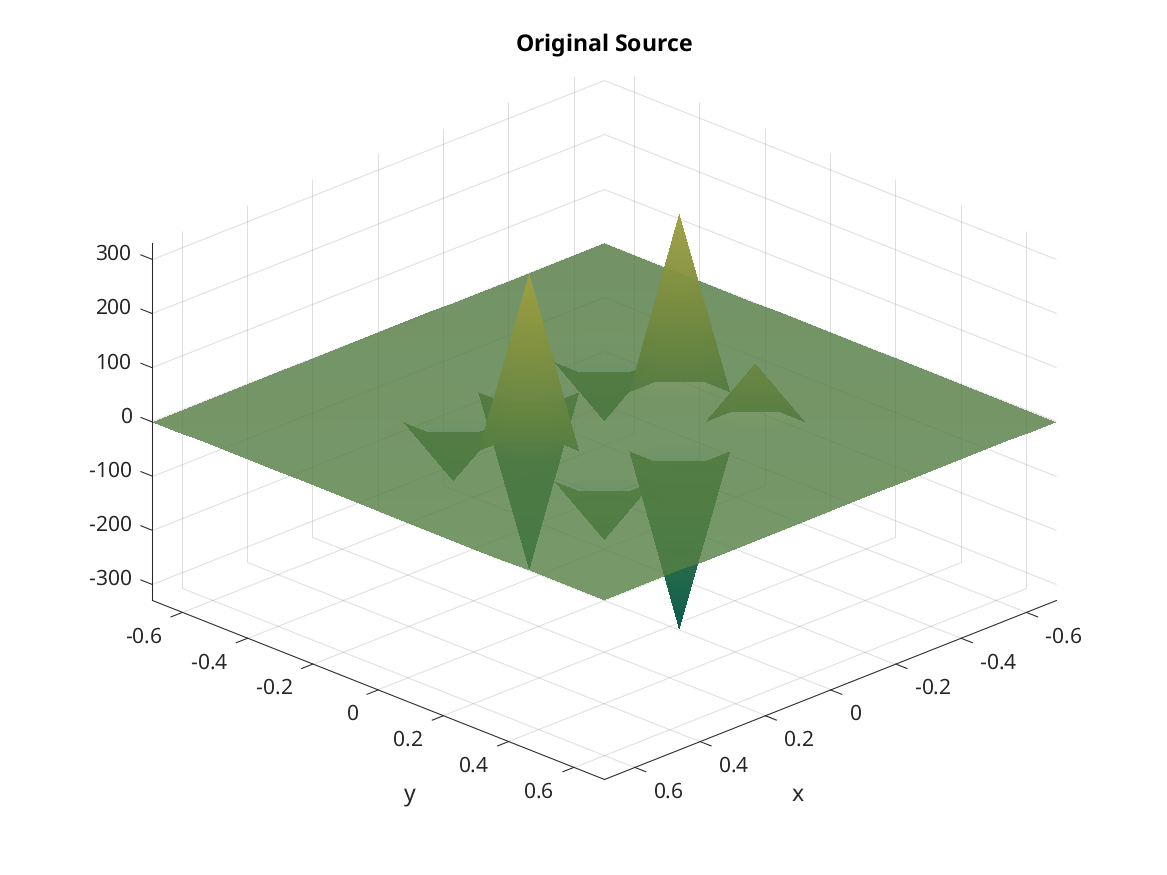}
        \caption{}
        \label{fig:1j}
    \end{subfigure}
    \hfill
    \begin{subfigure}[b]{0.2\textwidth}
        \includegraphics[width=\textwidth]{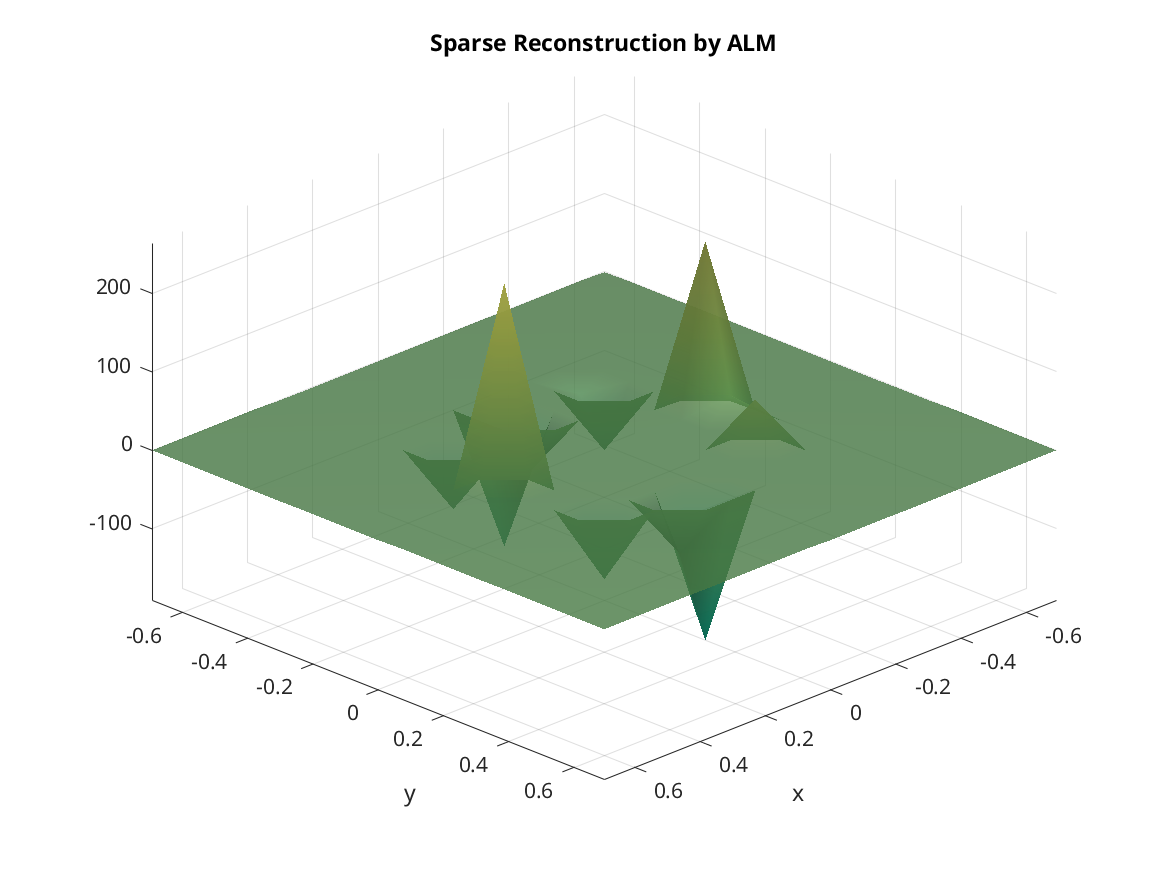}
        \caption{}
        \label{}
    \end{subfigure}
    \hfill
    \begin{subfigure}[b]{0.2\textwidth}
        \includegraphics[width=\textwidth]{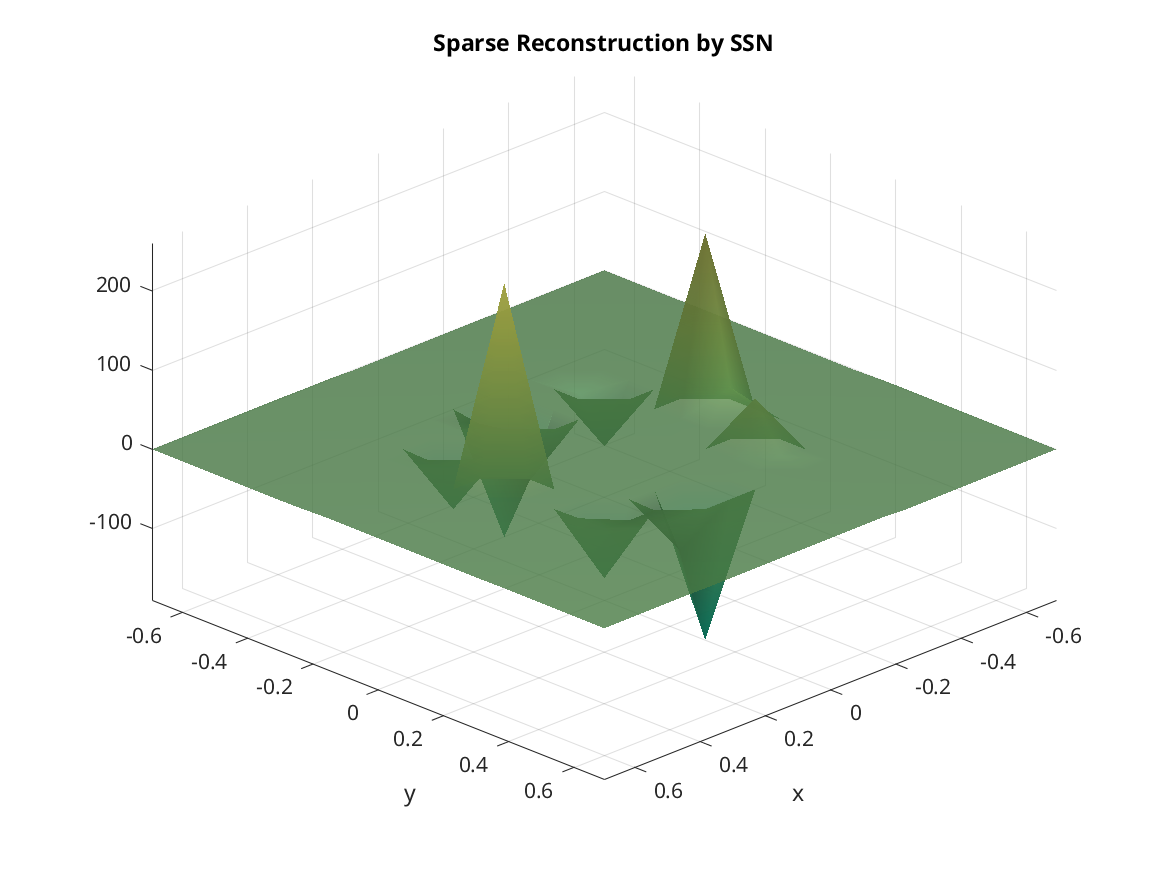}
        \caption{}
        \label{}
    \end{subfigure}
    \hfill
    \begin{subfigure}[b]{0.2\textwidth}
        \includegraphics[width=\textwidth]{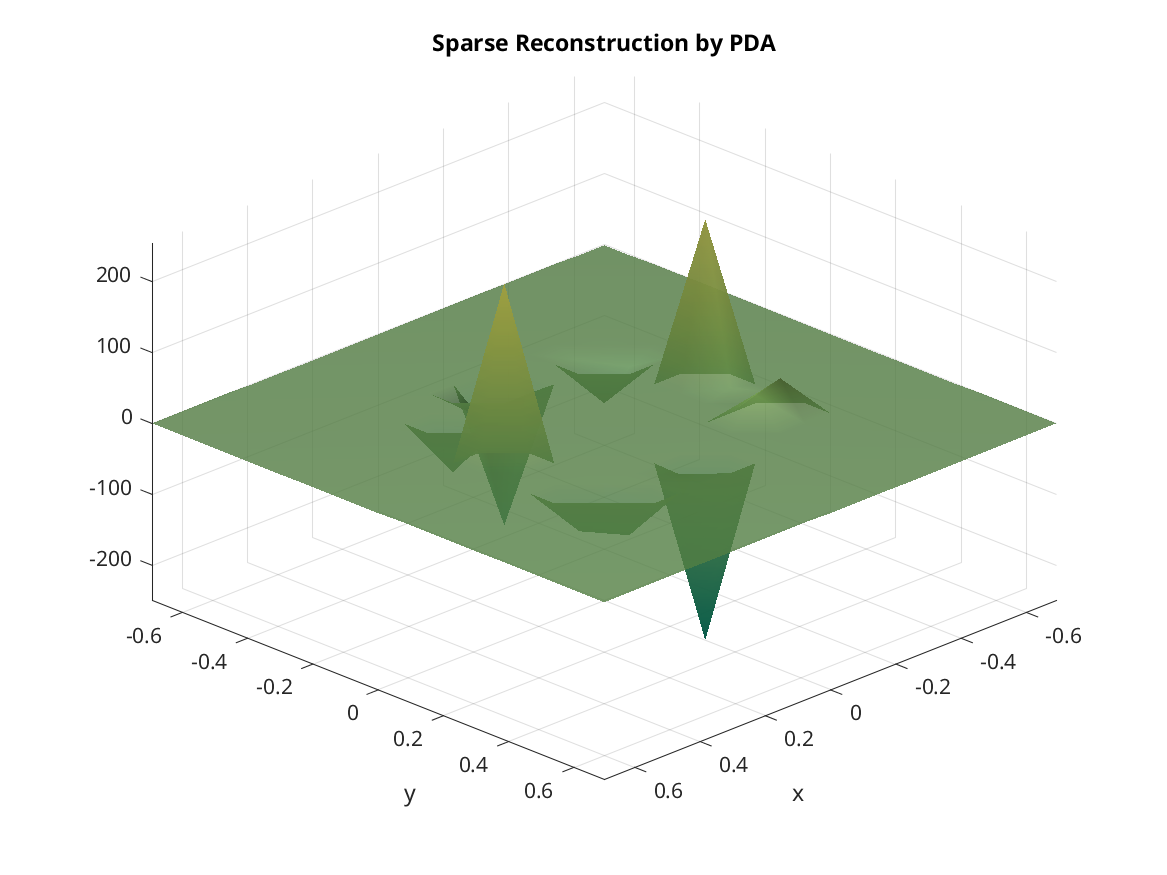}
        \caption{}
        \label{}
    \end{subfigure}
    \caption{Reconstruction of acoustic sources with multiple peaks in homogeneous media with $k=6$. The leftmost column shows the exact acoustic sources with one, four, six, and eight peaks, respectively. The images in the second, the third from the left, and the rightmost columns are reconstructed results of ALM, SSN, and PDA, respectively.  The images in the first, second, third, and fourth rows are the sources with one, four, six, and eight peaks, respectively.  }
    \label{fig:big_figure1}
\end{figure}

\begin{figure}[htbp]
    \centering
    
    \begin{subfigure}[b]{0.2\textwidth}
        \includegraphics[width=\textwidth]{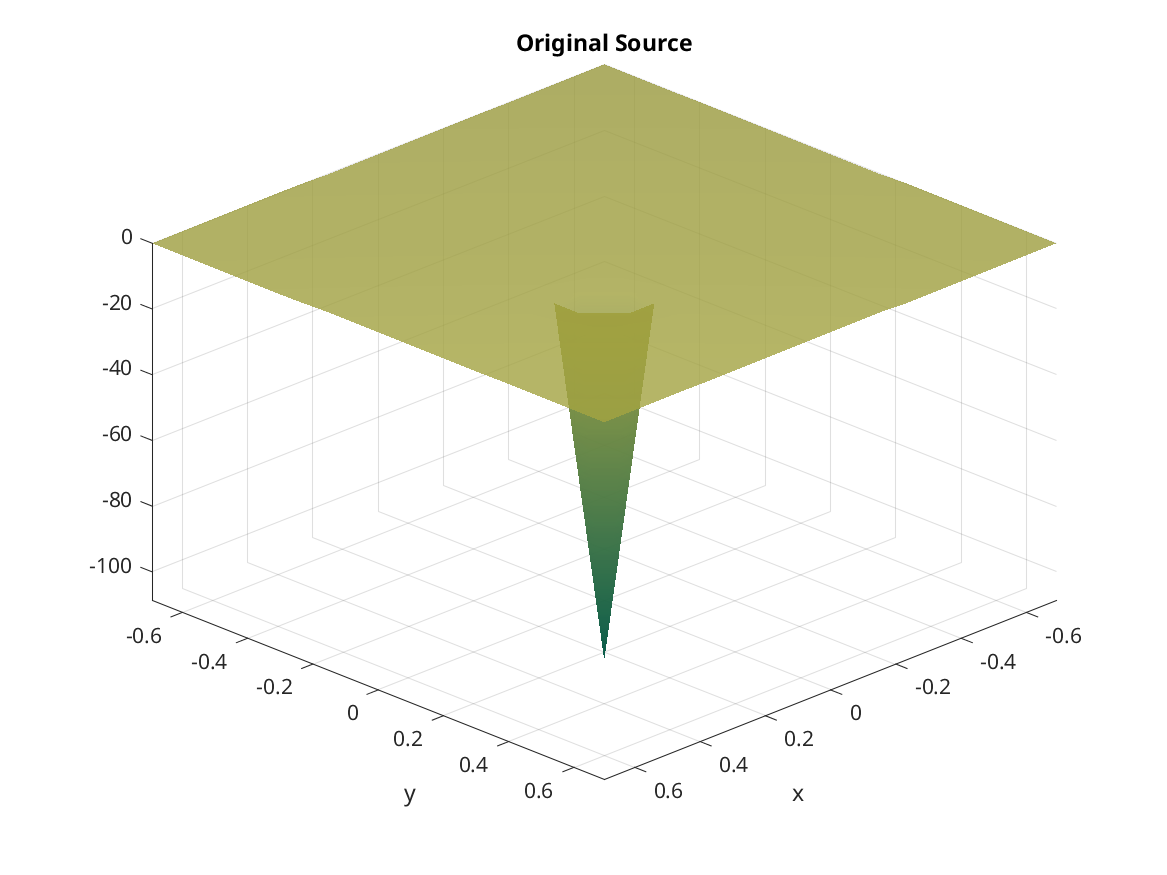}
        \caption{}
        \label{}
    \end{subfigure}
    \hfill
    \begin{subfigure}[b]{0.2\textwidth}
        \includegraphics[width=\textwidth]{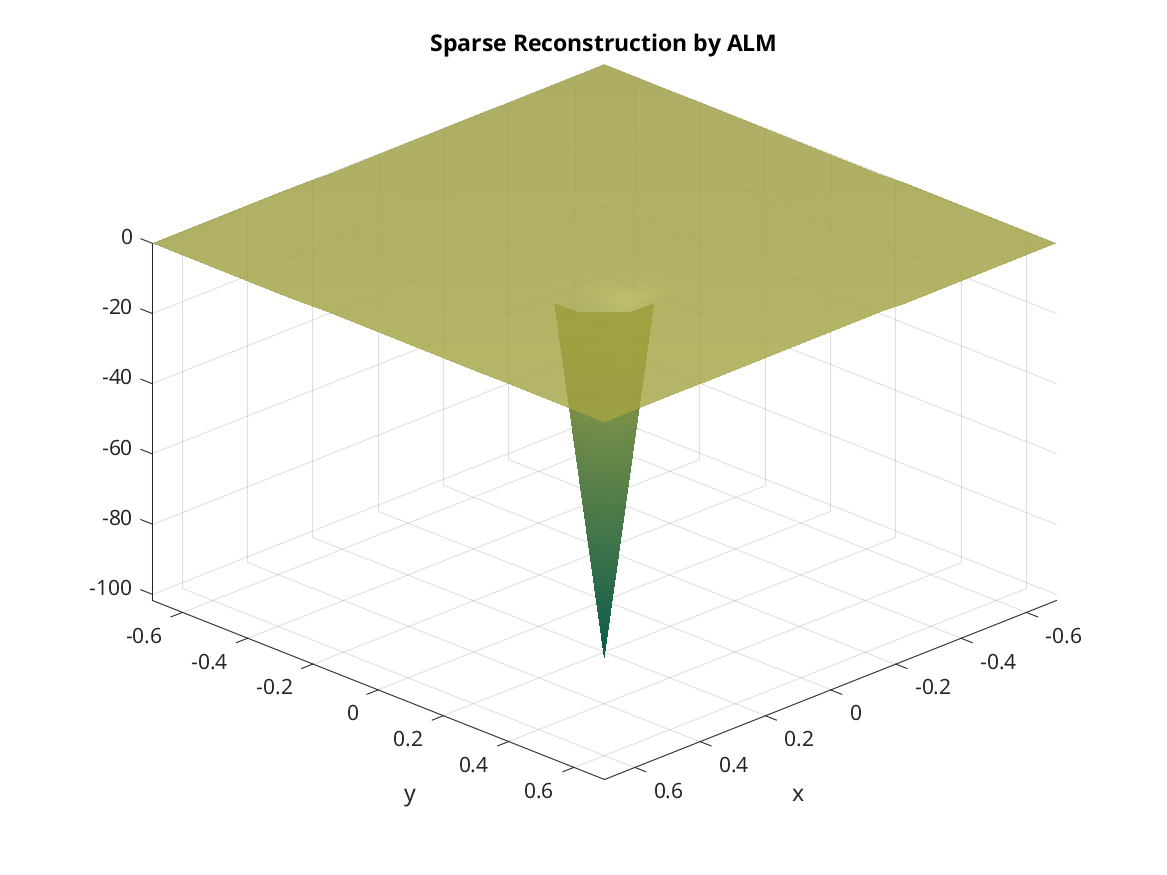}
        \caption{}
        \label{}
    \end{subfigure}
    \hfill
    \begin{subfigure}[b]{0.2\textwidth}
        \includegraphics[width=\textwidth]{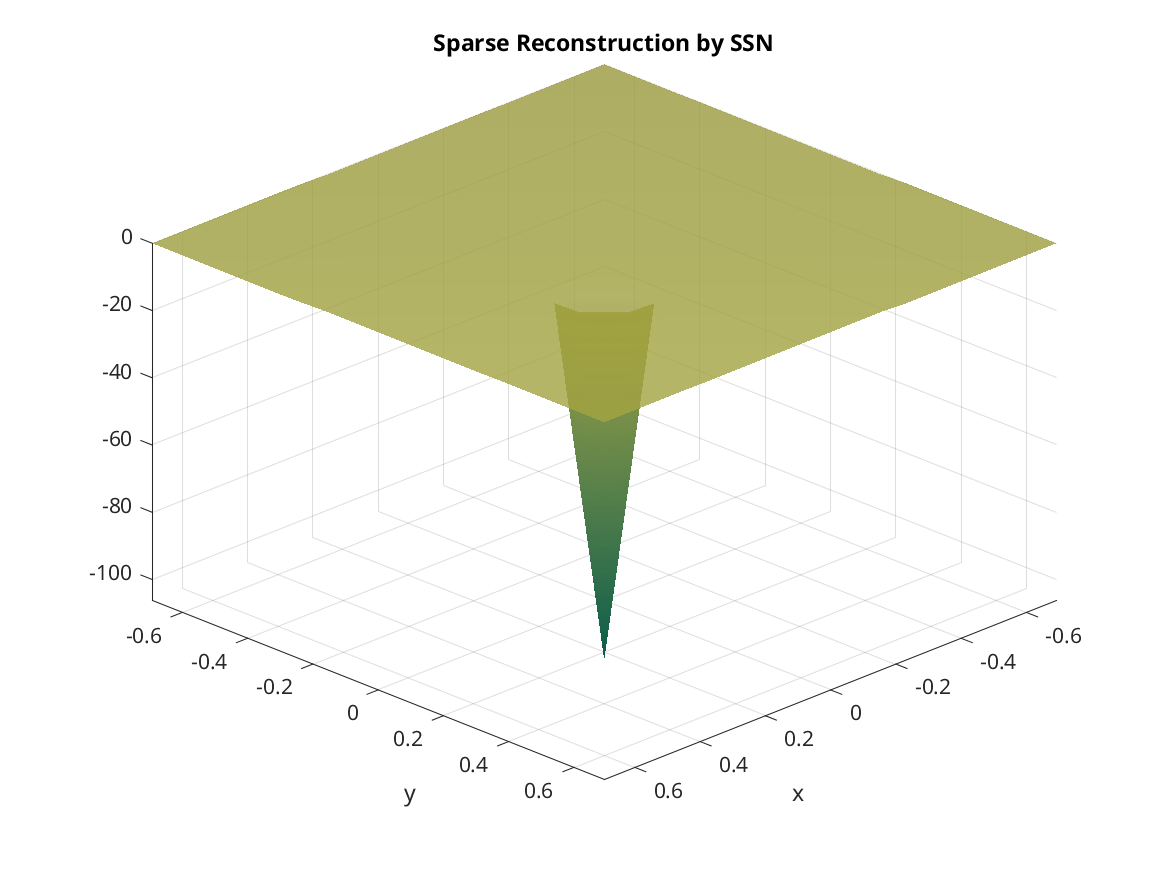}
        \caption{}
        \label{}
    \end{subfigure}
    \hfill
    \begin{subfigure}[b]{0.2\textwidth}
        \includegraphics[width=\textwidth]{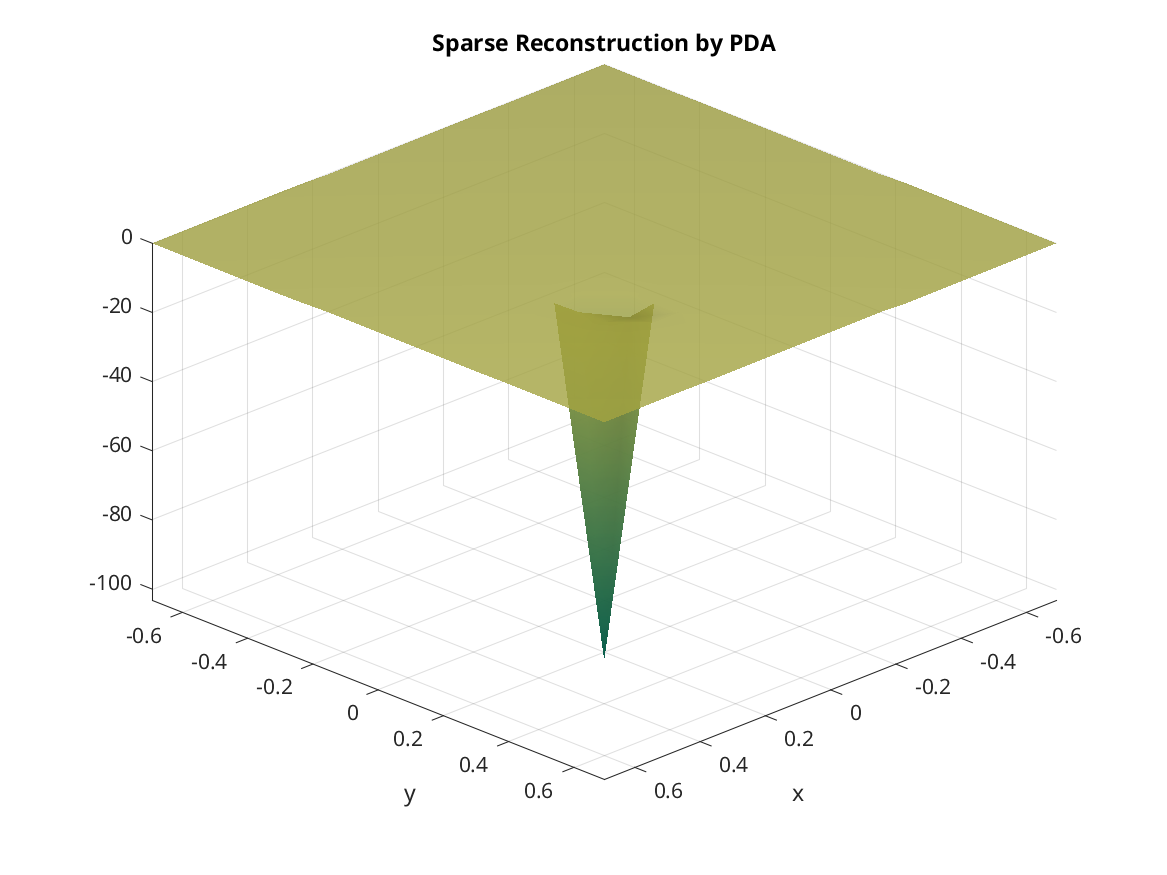}
        \caption{}
        \label{}
    \end{subfigure}
    
    \begin{subfigure}[b]{0.2\textwidth}
        \includegraphics[width=\textwidth]{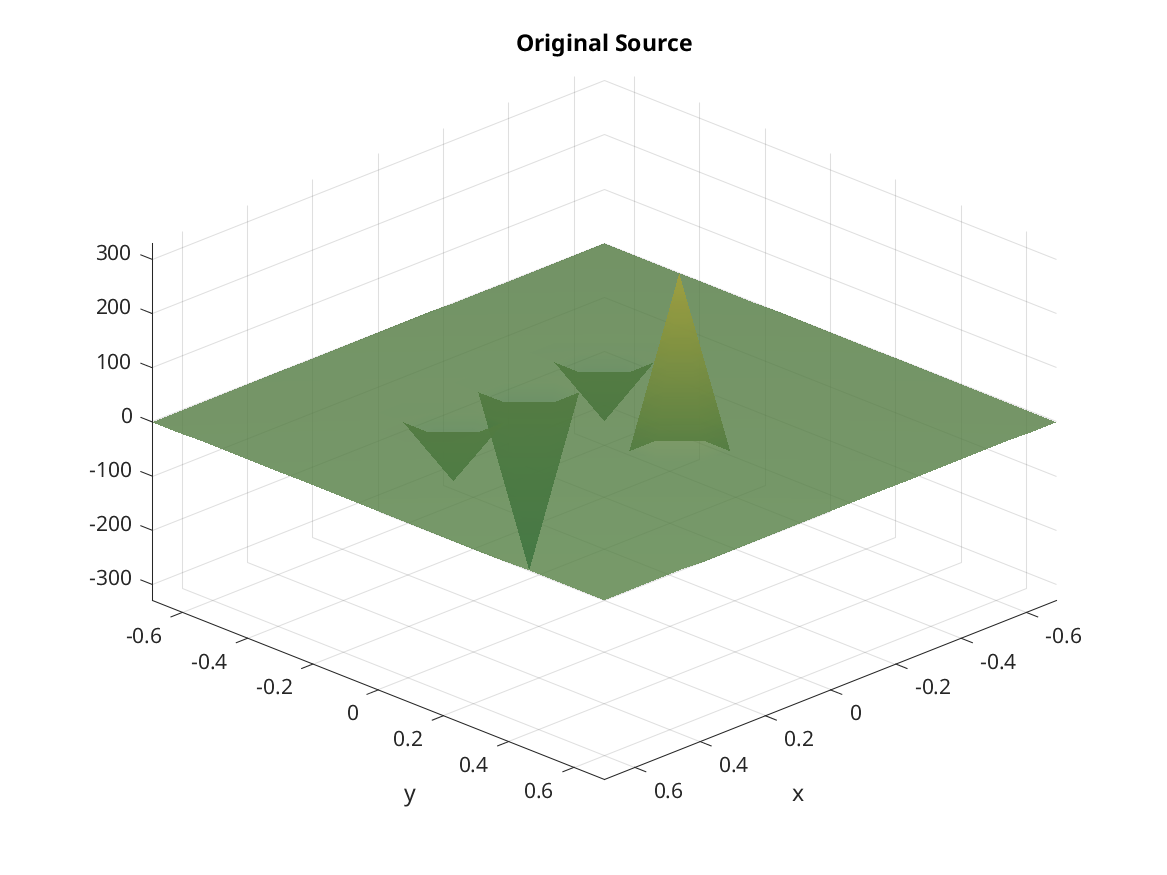}
        \caption{}
        \label{}
    \end{subfigure}
    \hfill
    \begin{subfigure}[b]{0.2\textwidth}
        \includegraphics[width=\textwidth]{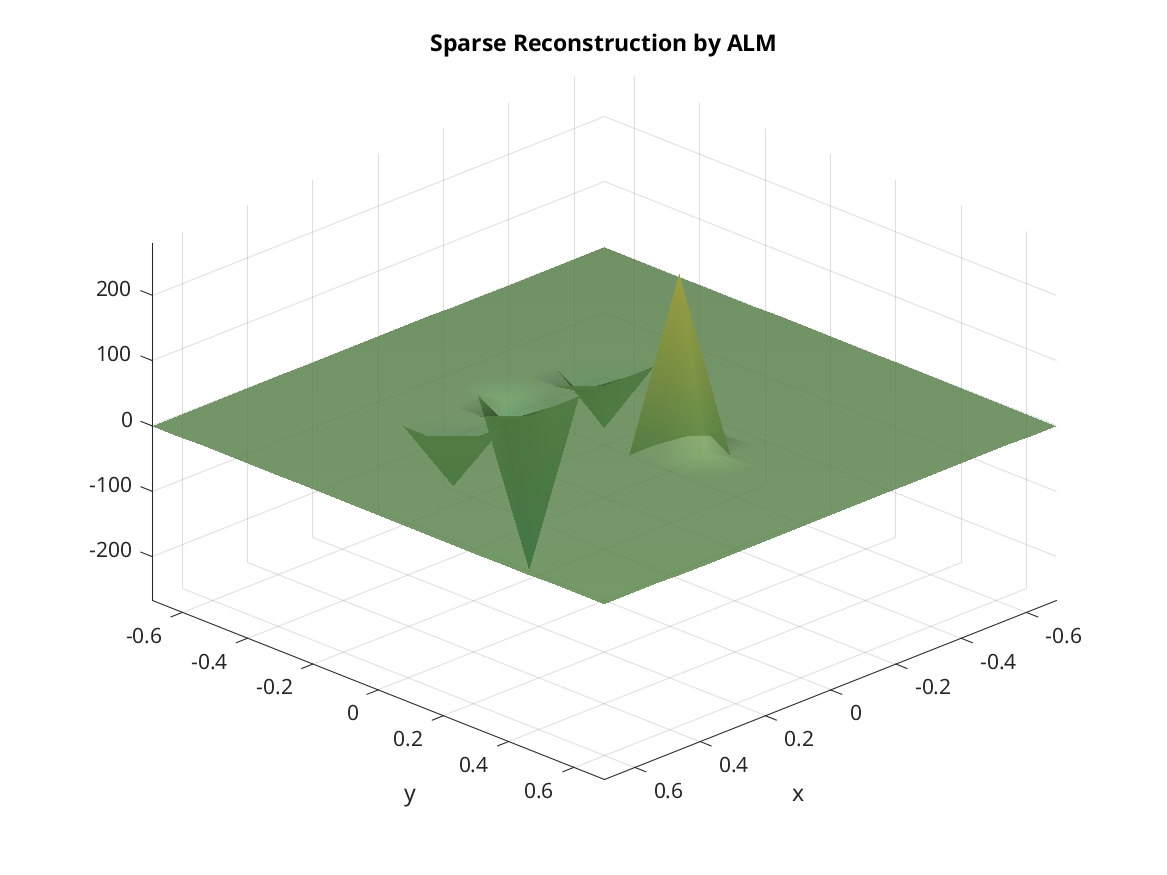}
        \caption{}
        \label{fig:1d}
    \end{subfigure}
    \hfill
    \begin{subfigure}[b]{0.2\textwidth}
        \includegraphics[width=\textwidth]{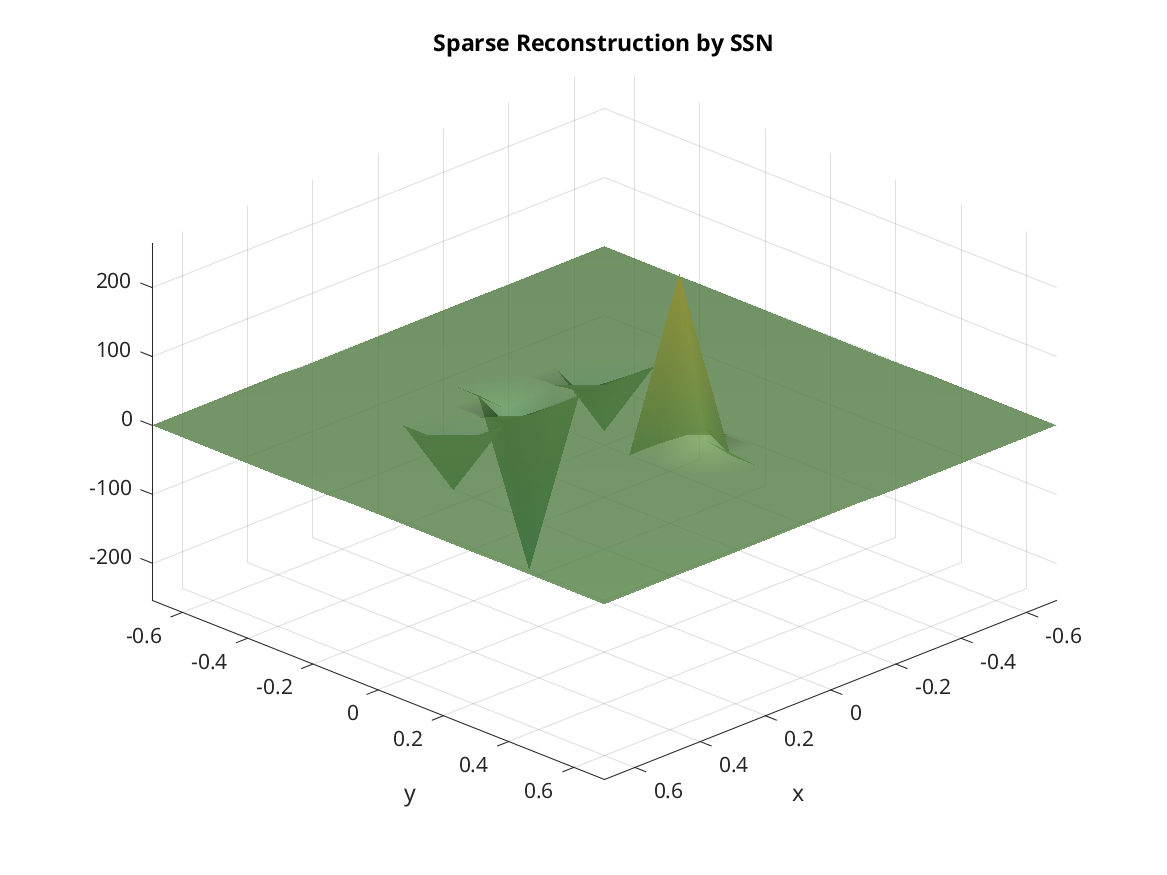}
        \caption{}
        \label{}
    \end{subfigure}
    \hfill
    \begin{subfigure}[b]{0.2\textwidth}
        \includegraphics[width=\textwidth]{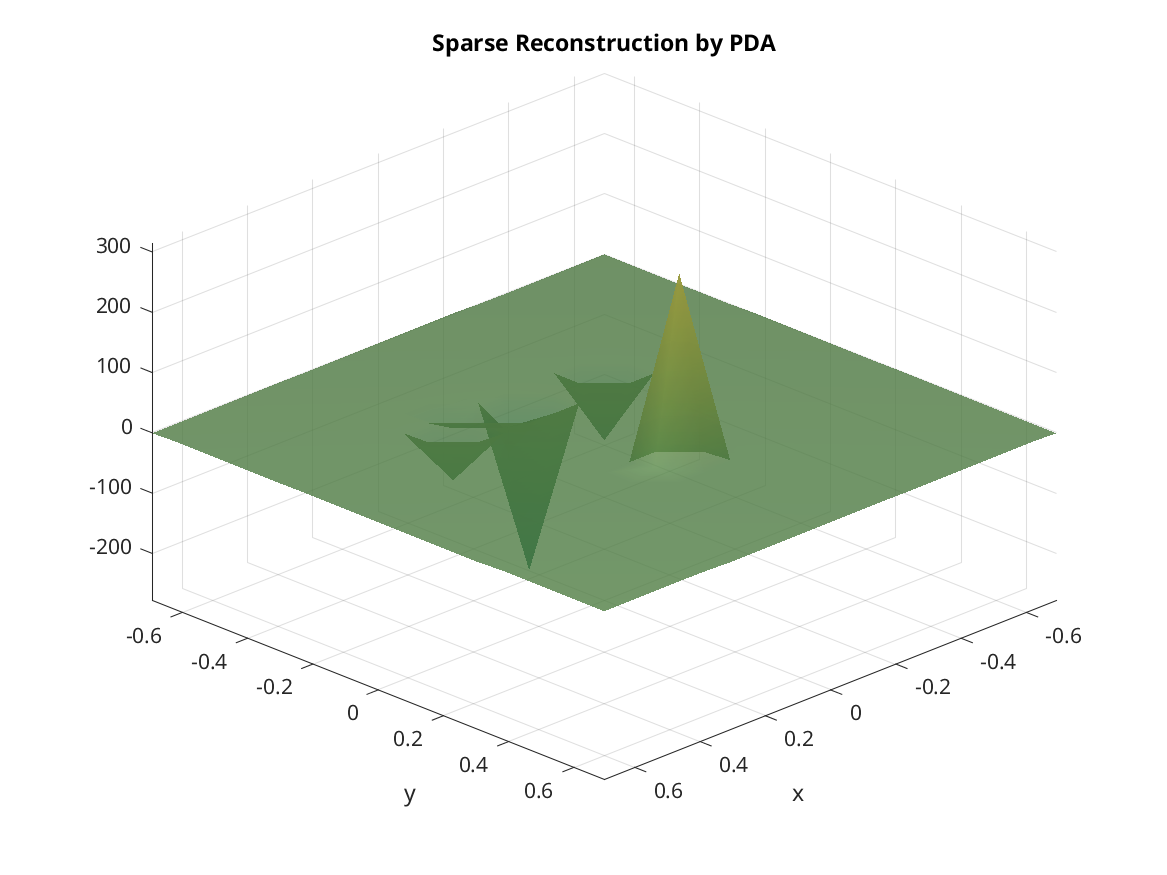}
        \caption{}
        \label{}
    \end{subfigure}
    
    \begin{subfigure}[b]{0.2\textwidth}
        \includegraphics[width=\textwidth]{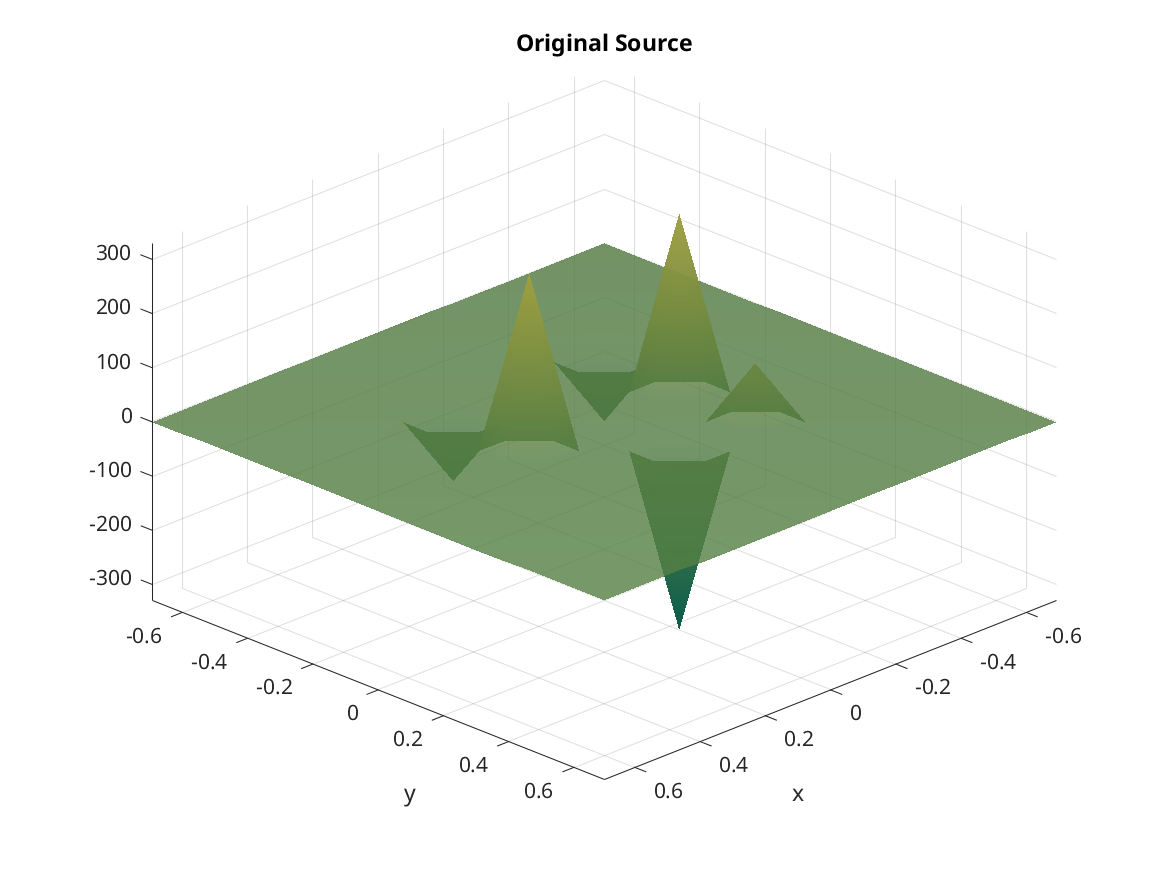}
        \caption{}
        \label{}
    \end{subfigure}
    \hfill
    \begin{subfigure}[b]{0.2\textwidth}
        \includegraphics[width=\textwidth]{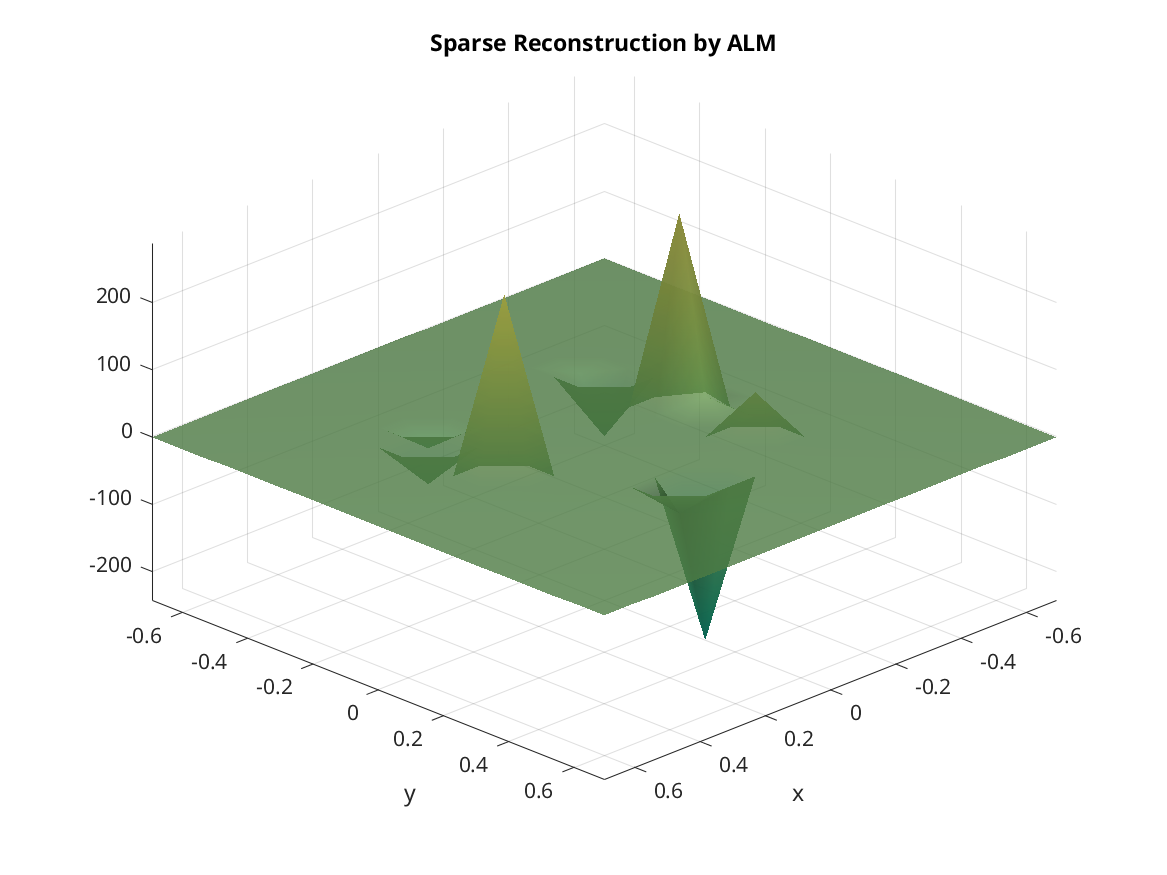}
        \caption{}
        \label{}
    \end{subfigure}
    \hfill
    \begin{subfigure}[b]{0.2\textwidth}
        \includegraphics[width=\textwidth]{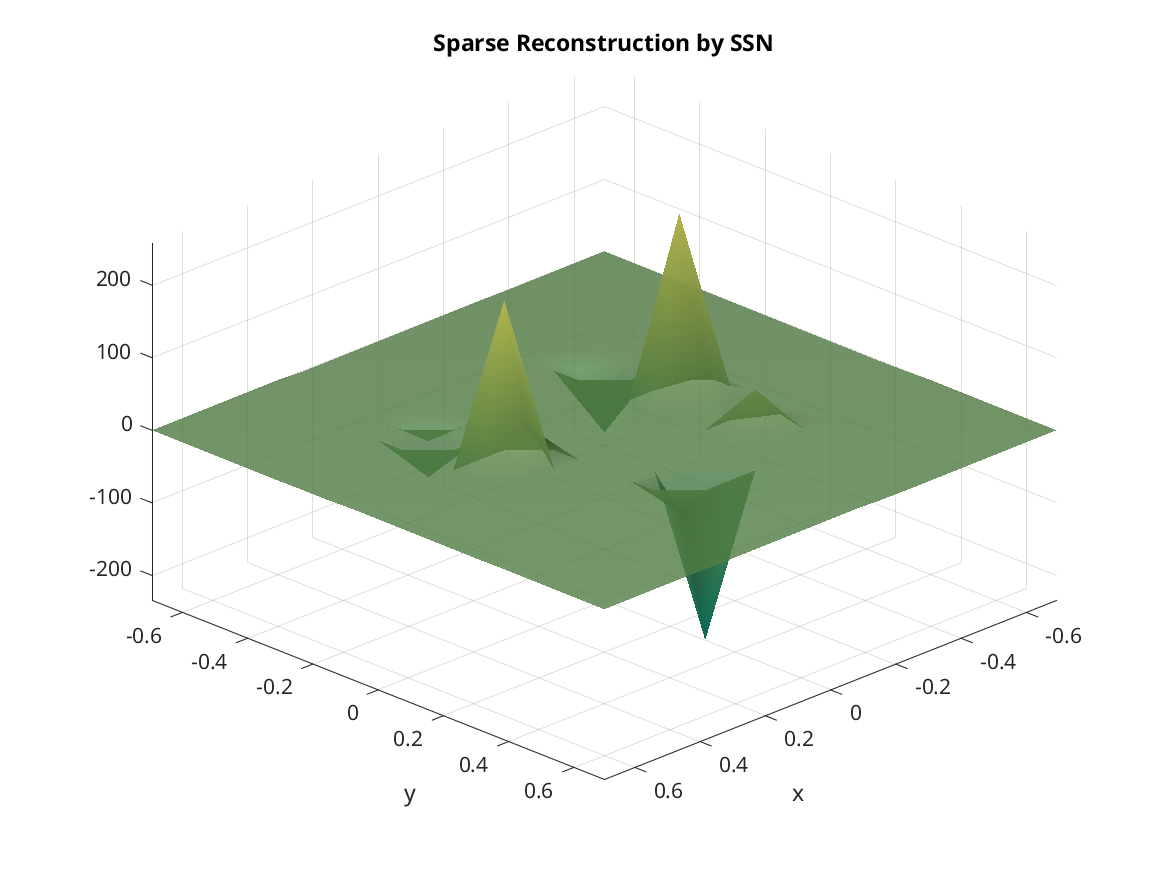}
        \caption{}
        \label{}
    \end{subfigure}
    \hfill
    \begin{subfigure}[b]{0.2\textwidth}
        \includegraphics[width=\textwidth]{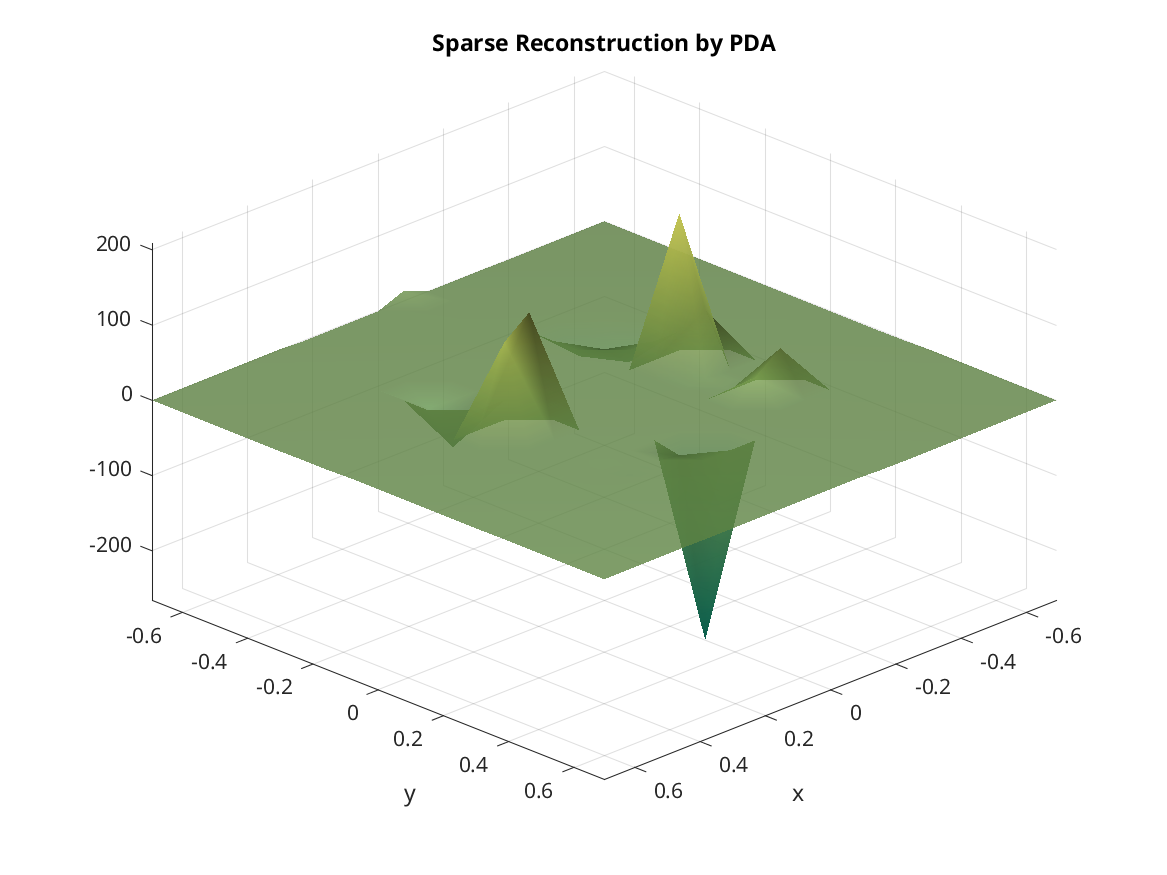}
        \caption{}
        \label{}
    \end{subfigure}
    
    \begin{subfigure}[b]{0.2\textwidth}
        \includegraphics[width=\textwidth]{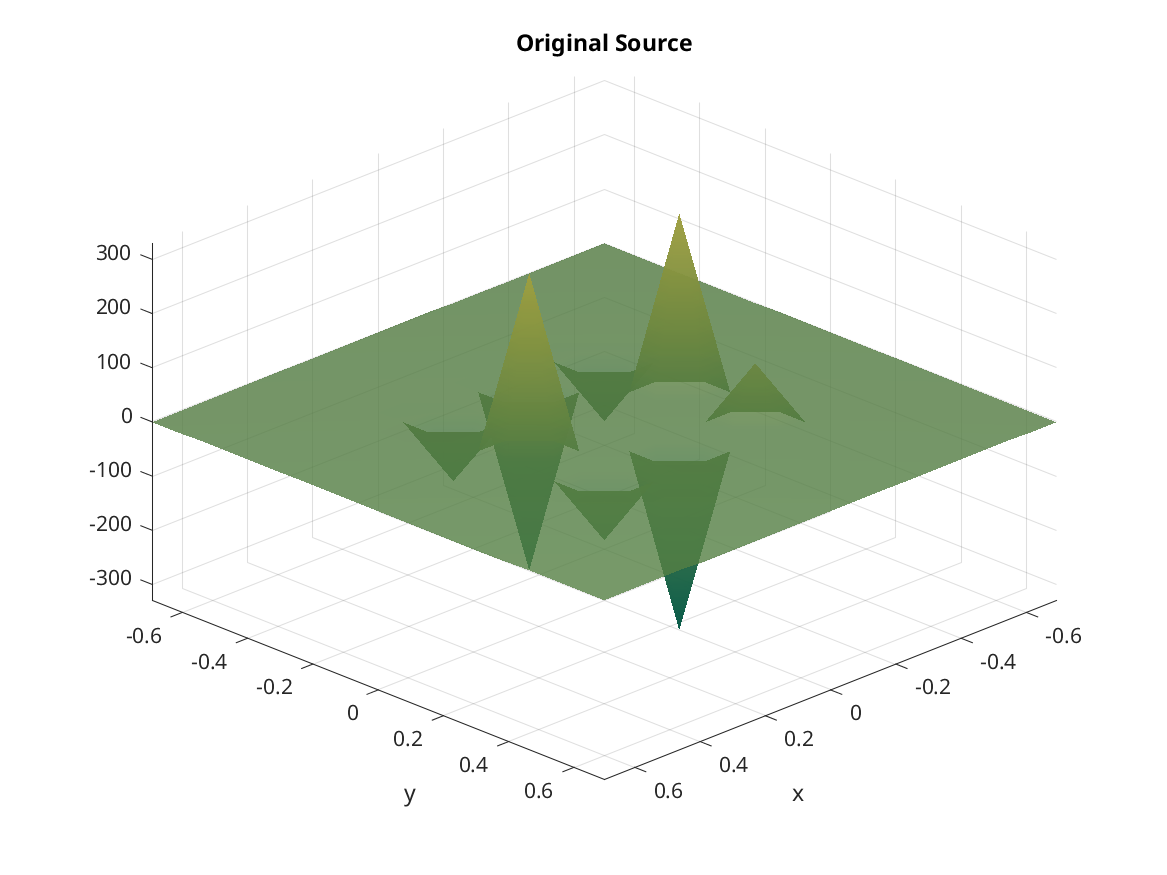}
        \caption{}
        \label{}
    \end{subfigure}
    \hfill
    \begin{subfigure}[b]{0.2\textwidth}
        \includegraphics[width=\textwidth]{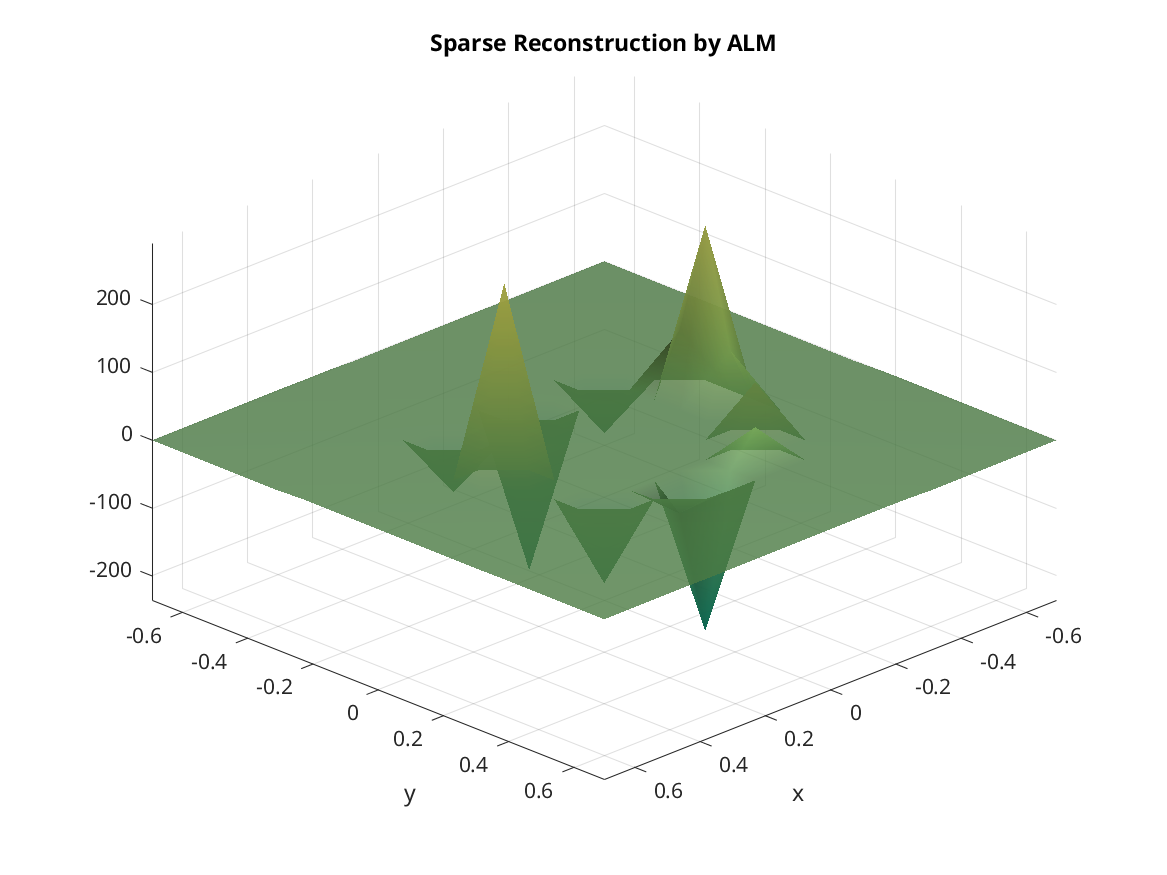}
        \caption{}
        \label{}
    \end{subfigure}
    \hfill
    \begin{subfigure}[b]{0.2\textwidth}
        \includegraphics[width=\textwidth]{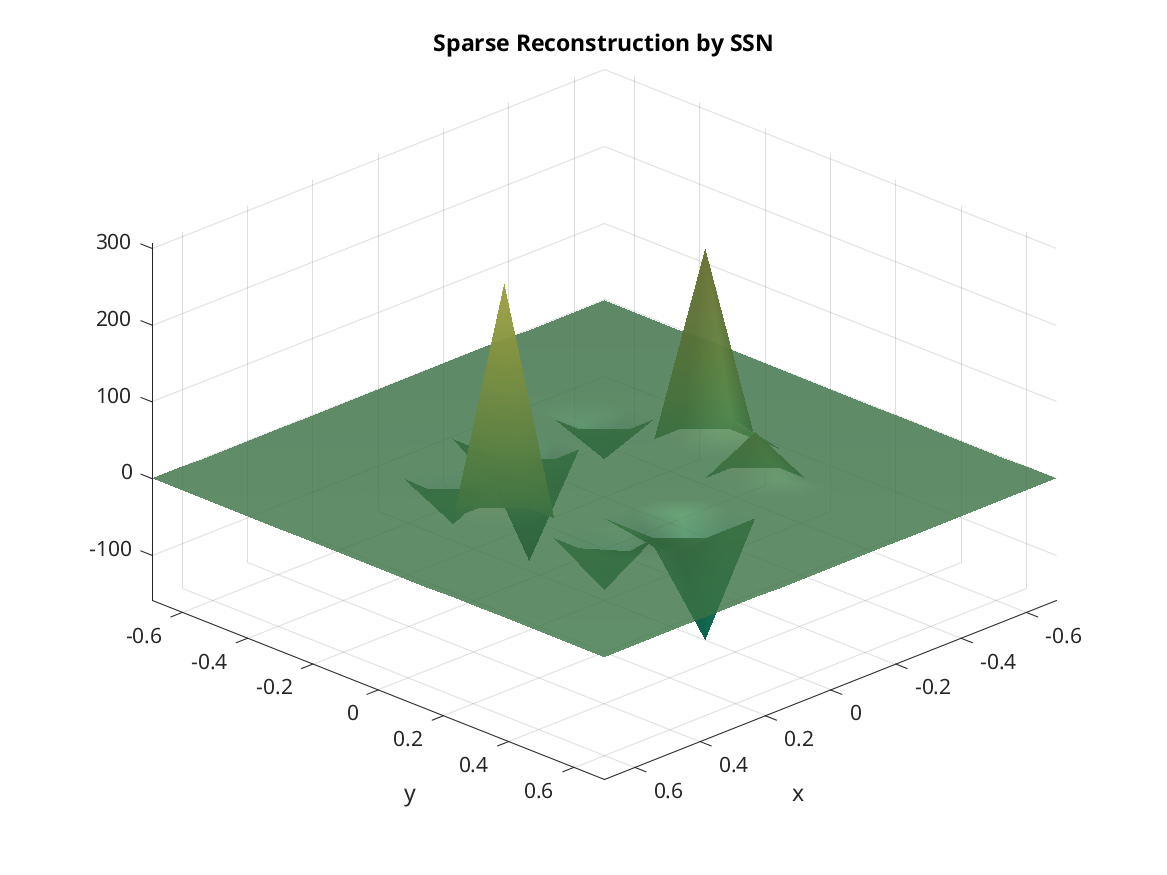}
        \caption{}
        \label{}
    \end{subfigure}
    \hfill
    \begin{subfigure}[b]{0.2\textwidth}
        \includegraphics[width=\textwidth]{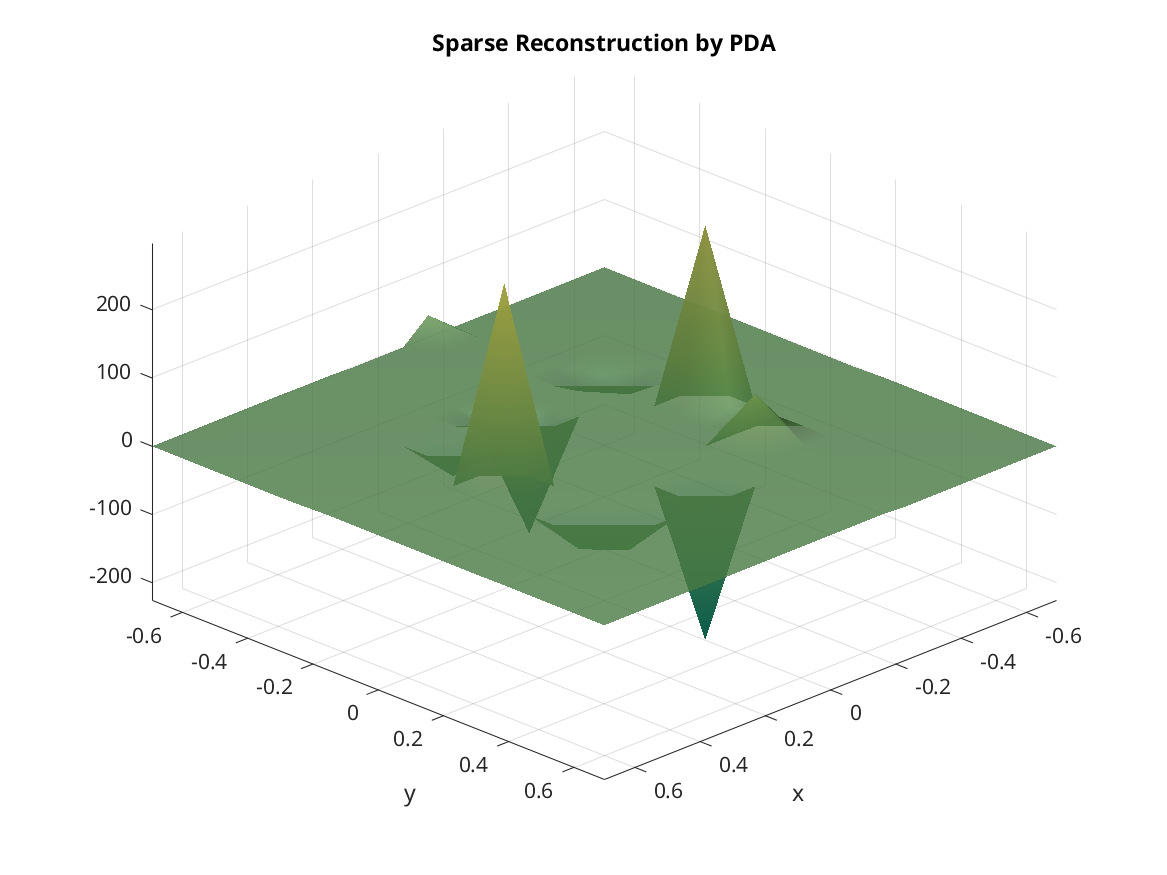}
        \caption{}
        \label{}
    \end{subfigure}
    \caption{Reconstruction of acoustic sources with multiple peaks in inhomogeneous media with $k=6$. The information of the acoustic sources and the corresponding reconstruction algorithms are the same as in Figure \ref{fig:big_figure1}. }
    \label{fig:big_figure2}
\end{figure}

Figures \ref{fig:big_figure1} and \ref{fig:big_figure2}  show reconstructions of point sources with \%1 relative Gaussian noise for homogeneous and inhomogeneous media, respectively. Table \ref{tab:results04} collects the running time and relative errors of three different algorithms. For our experiments, we observed no significant differences in running time between the homogeneous and inhomogeneous cases. We thus only present the running time for the homogeneous case.

\begin{figure}[htbp]
    \centering
    
    \begin{subfigure}[b]{0.2\textwidth}
        \includegraphics[width=\textwidth]{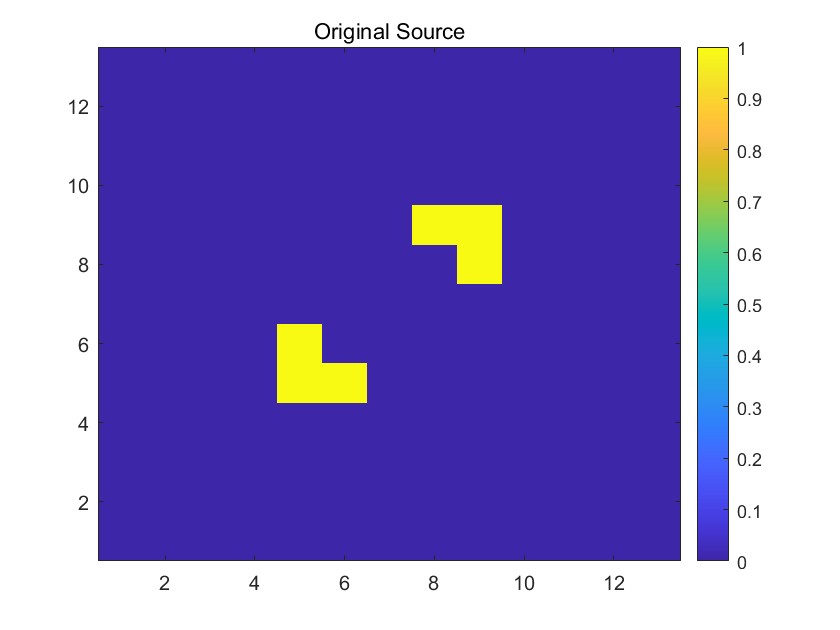}
        \caption{}
        \label{}
    \end{subfigure}
    \hfill
    \begin{subfigure}[b]{0.2\textwidth}
        \includegraphics[width=\textwidth]{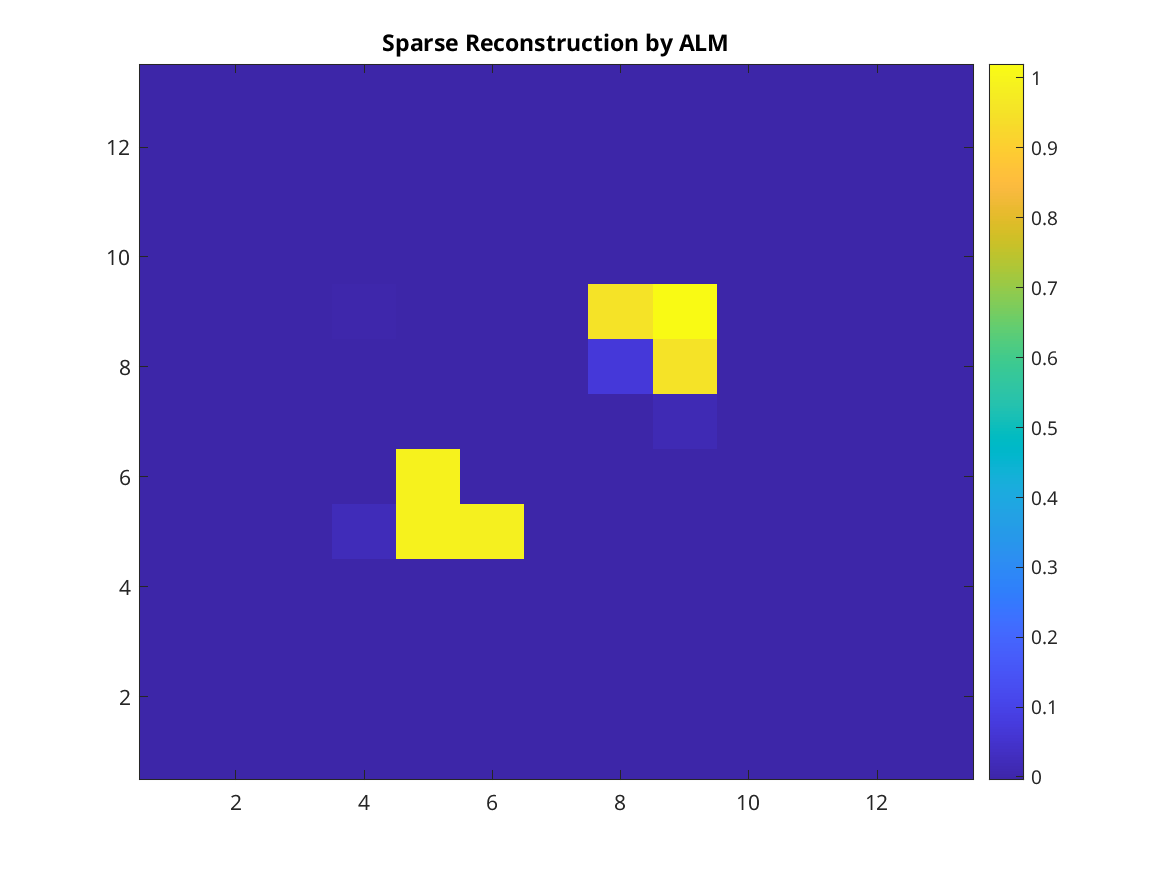}
        \caption{}
        \label{}
    \end{subfigure}
    \hfill
    \begin{subfigure}[b]{0.2\textwidth}
        \includegraphics[width=\textwidth]{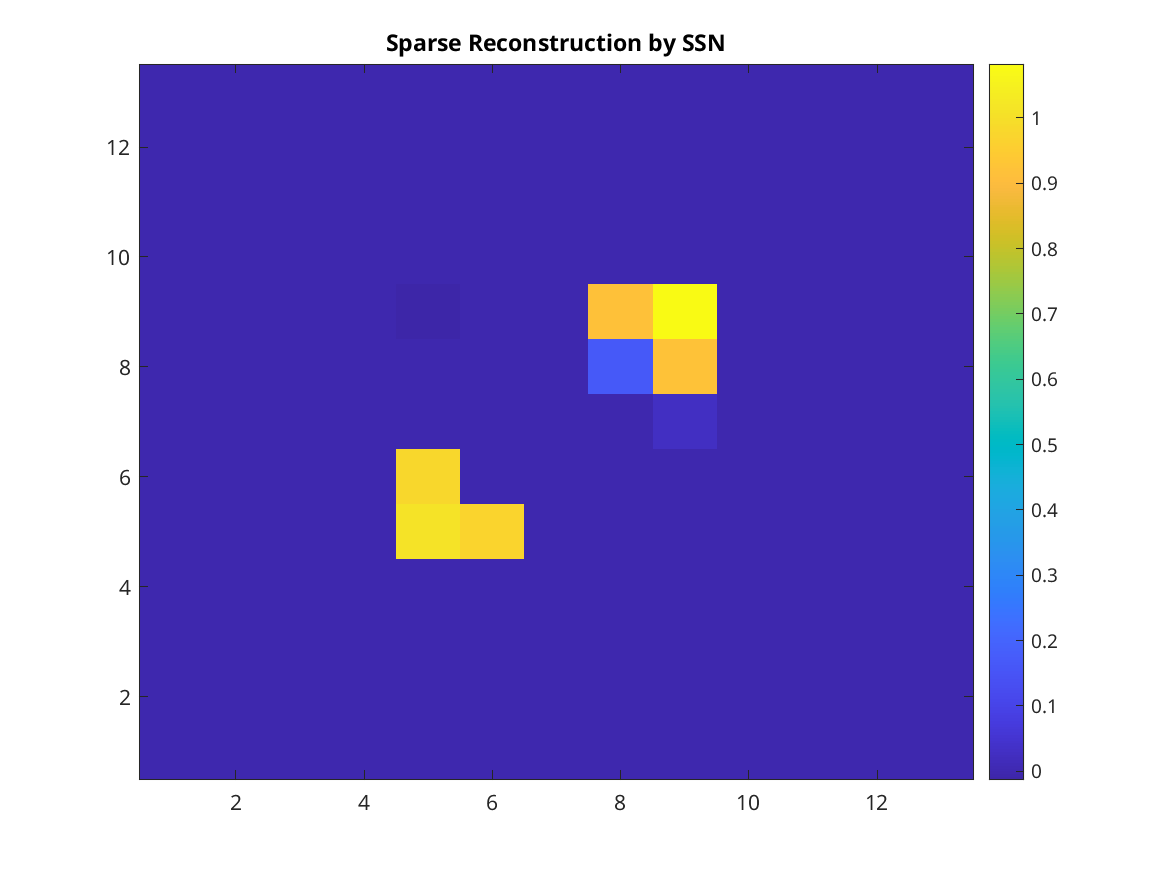}
        \caption{}
        \label{}
    \end{subfigure}
    \hfill
    \begin{subfigure}[b]{0.2\textwidth}
        \includegraphics[width=\textwidth]{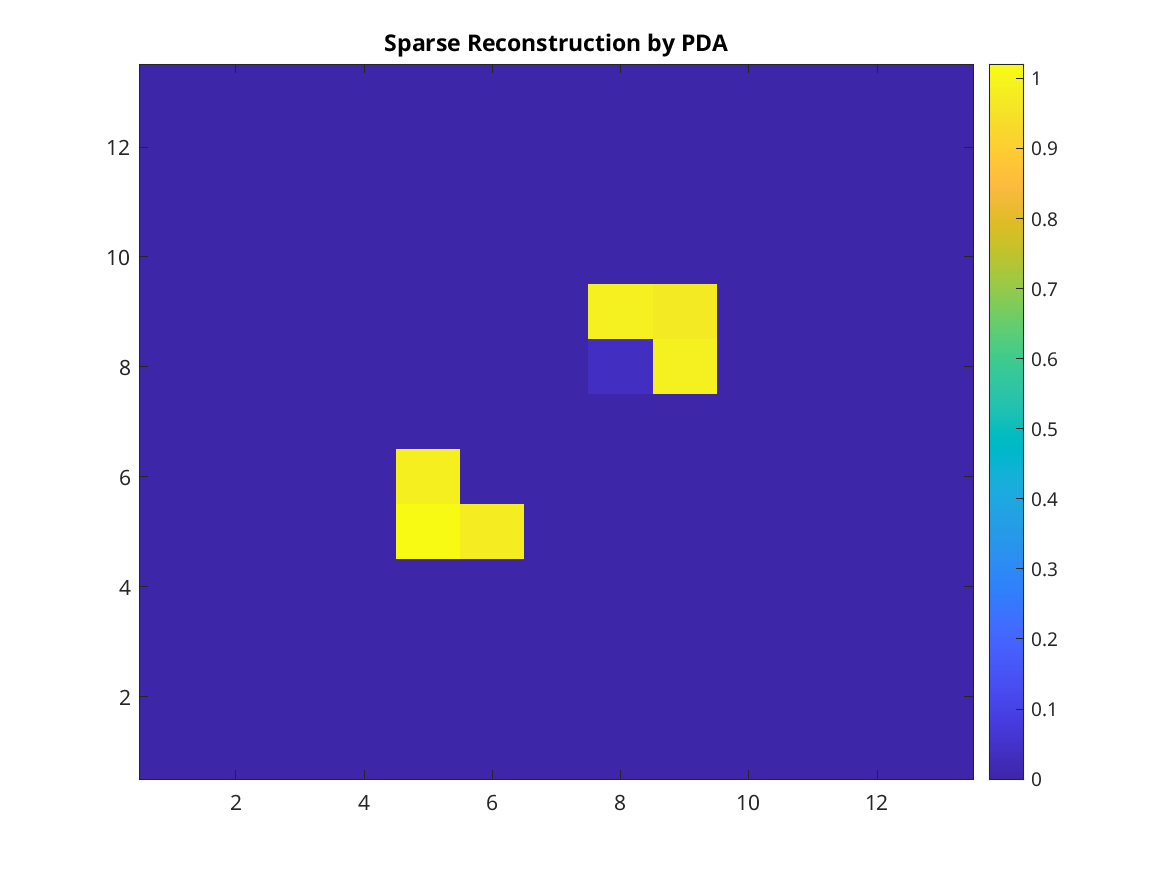}
        \caption{}
        \label{}
    \end{subfigure}
    
    \begin{subfigure}[b]{0.2\textwidth}
        \includegraphics[width=\textwidth]{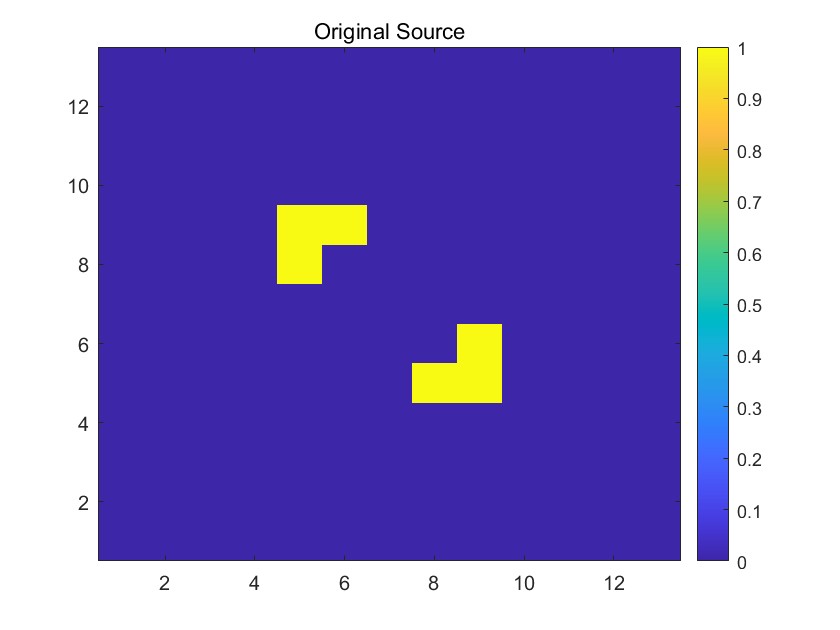}
        \caption{}
        \label{}
    \end{subfigure}
    \hfill
    \begin{subfigure}[b]{0.2\textwidth}
        \includegraphics[width=\textwidth]{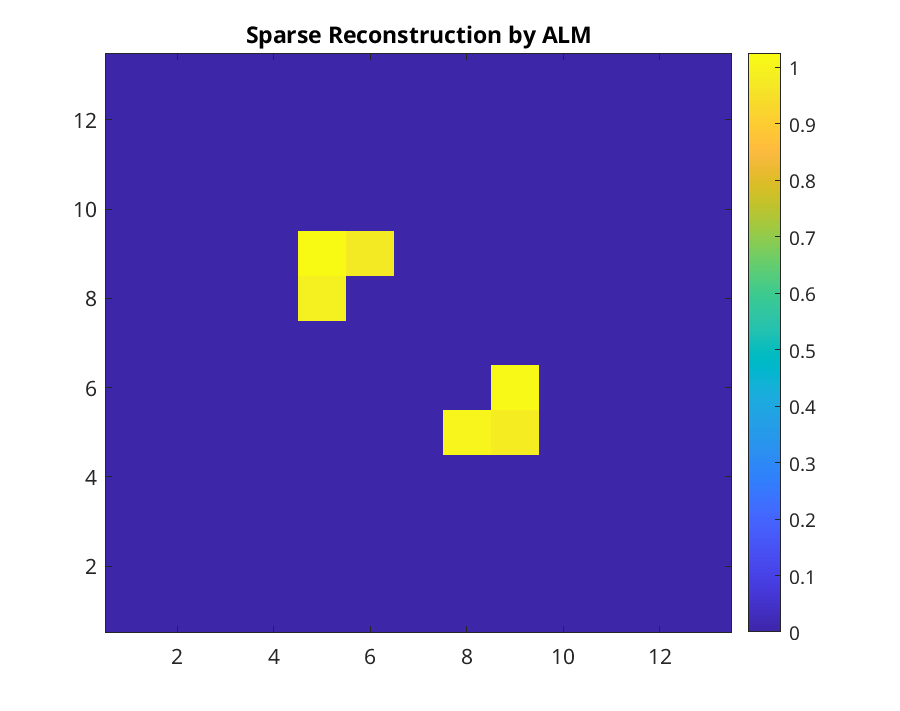}
        \caption{}
        \label{}
    \end{subfigure}
    \hfill
    \begin{subfigure}[b]{0.2\textwidth}
        \includegraphics[width=\textwidth]{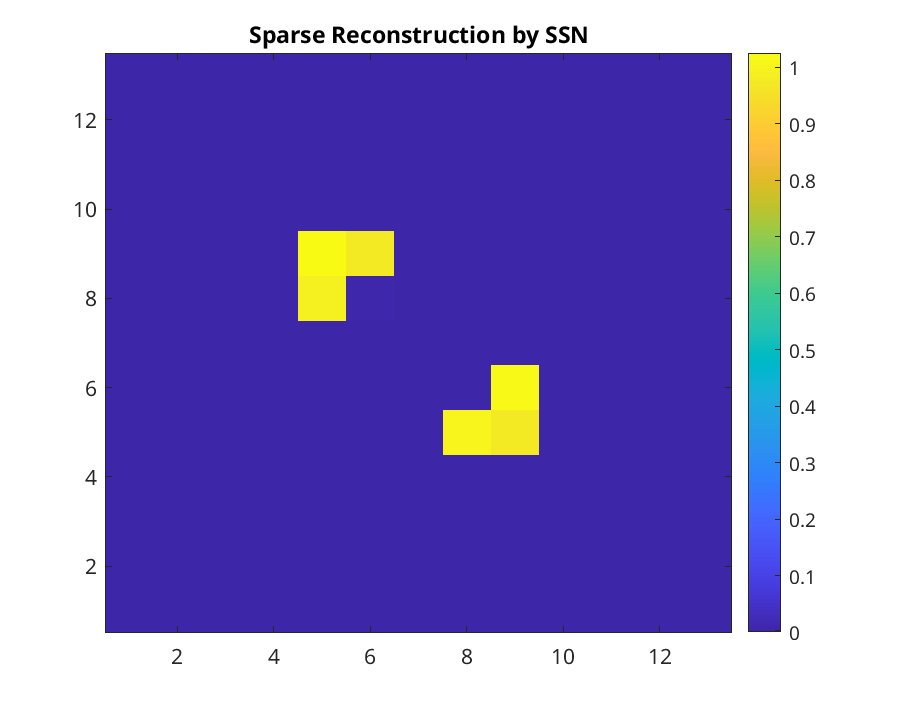}
        \caption{}
        \label{}
    \end{subfigure}
    \hfill
    \begin{subfigure}[b]{0.2\textwidth}
        \includegraphics[width=\textwidth]{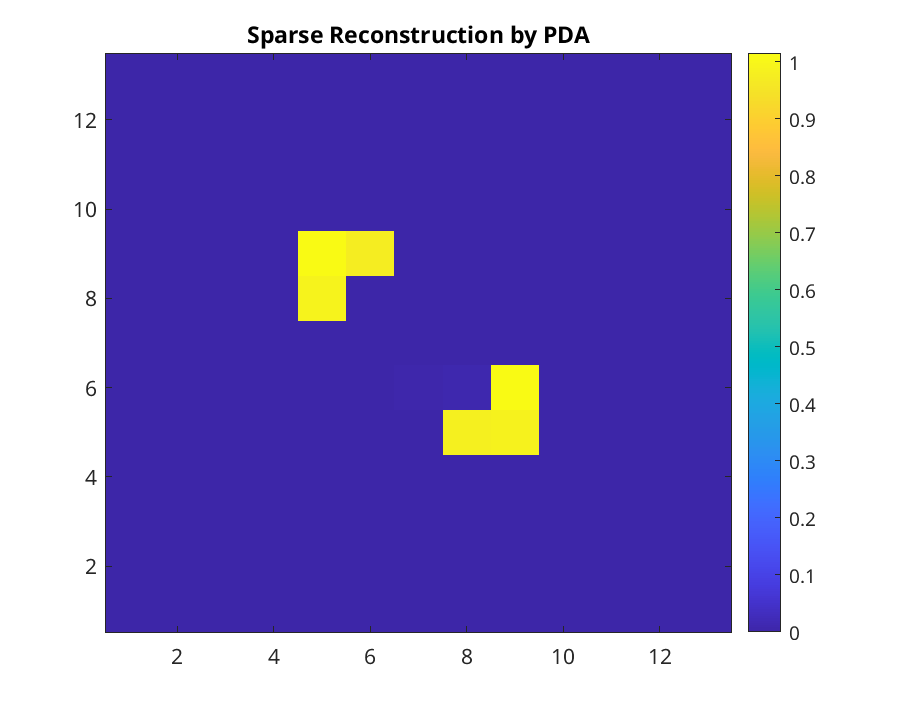}
        \caption{}
        \label{}
    \end{subfigure}
 \caption{Reconstructions of strip-shaped acoustic sources in homogeneous media with $k=4$, noise level 0.1\%, and $\alpha$ =2e-6,
$\alpha_0$= 2e-8.  The figures in the leftmost column are the original acoustic sources. The images in the second, the third from the left, and the rightmost columns are reconstructed results of ALM, SSN, and PDA, respectively.  The images in the first and second rows are the sources with ``skew diagonal" (the first row) and ``diagonal" (the second row) types of strips, respectively.}
    \label{fig:big_figure3}
\end{figure}

\begin{figure}[htbp]
    \centering
    
    \begin{subfigure}[b]{0.2\textwidth}
        \includegraphics[width=\textwidth]{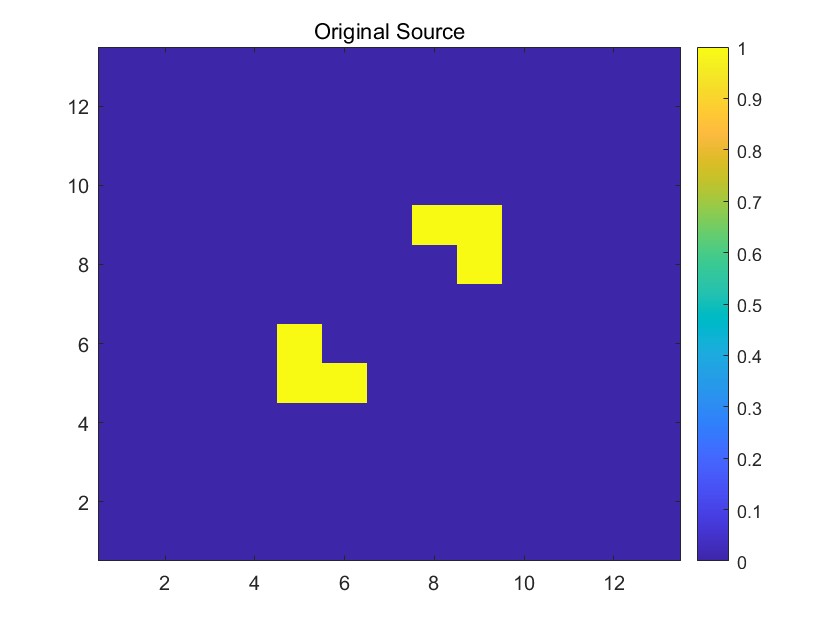}
        \caption{}
        \label{}
    \end{subfigure}
    \hfill
    \begin{subfigure}[b]{0.2\textwidth}
        \includegraphics[width=\textwidth]{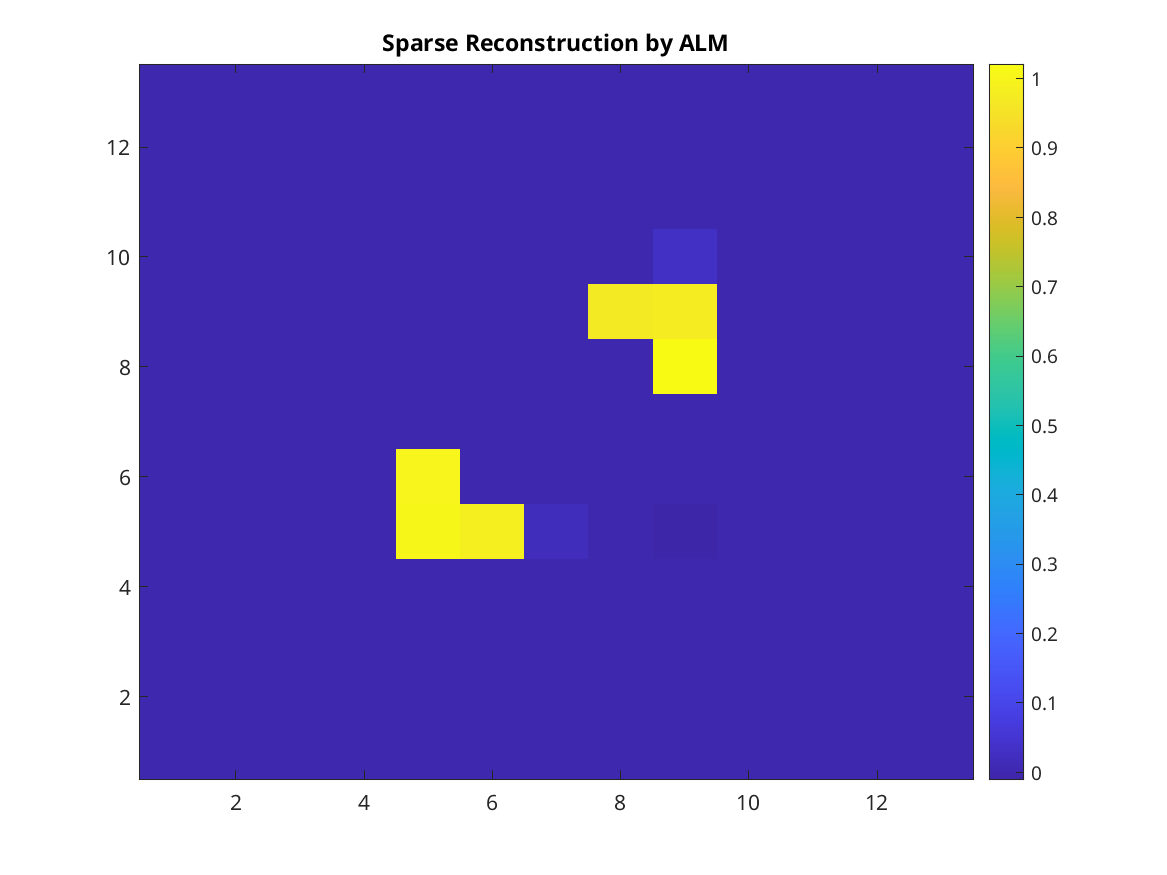}
        \caption{}
        \label{}
    \end{subfigure}
    \hfill
    \begin{subfigure}[b]{0.2\textwidth}
        \includegraphics[width=\textwidth]{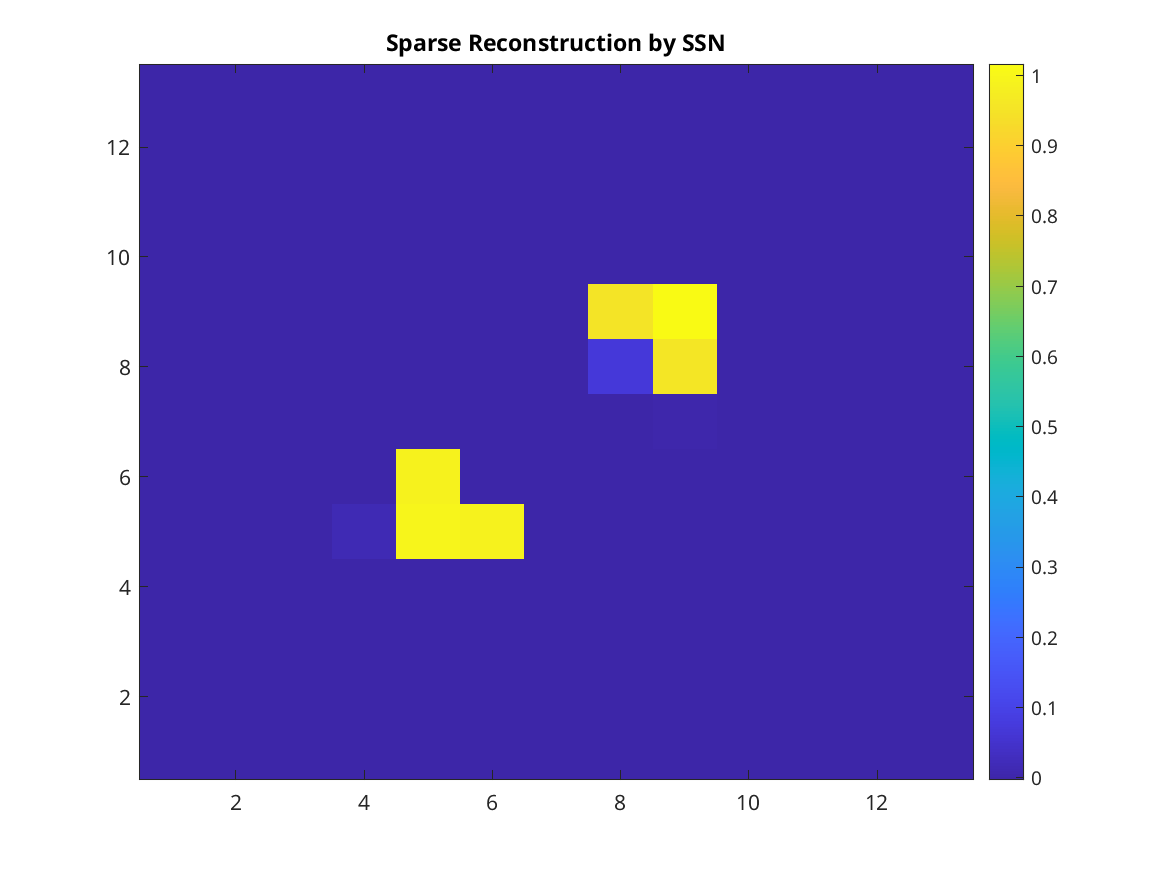}
        \caption{}
        \label{}
    \end{subfigure}
    \hfill
    \begin{subfigure}[b]{0.2\textwidth}
        \includegraphics[width=\textwidth]{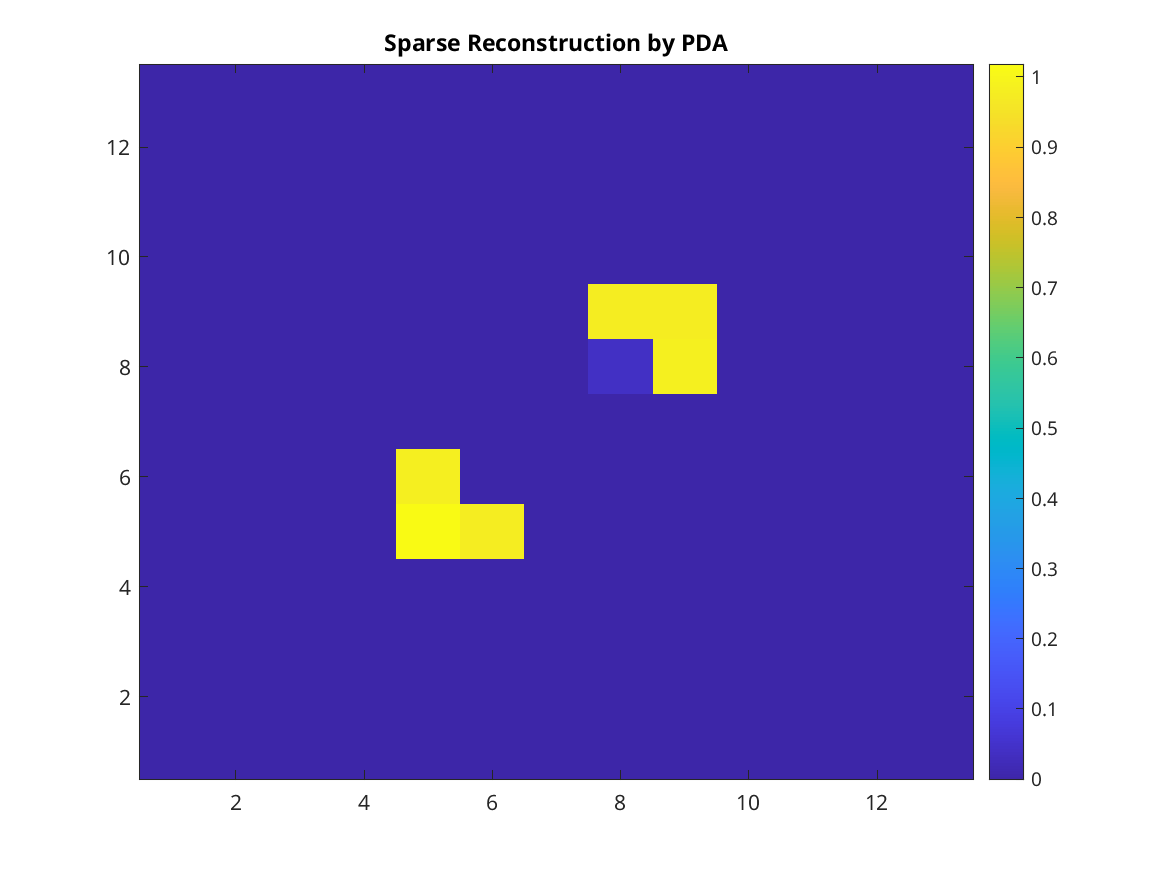}
        \caption{}
        \label{}
    \end{subfigure}
    
    \begin{subfigure}[b]{0.2\textwidth}
        \includegraphics[width=\textwidth]{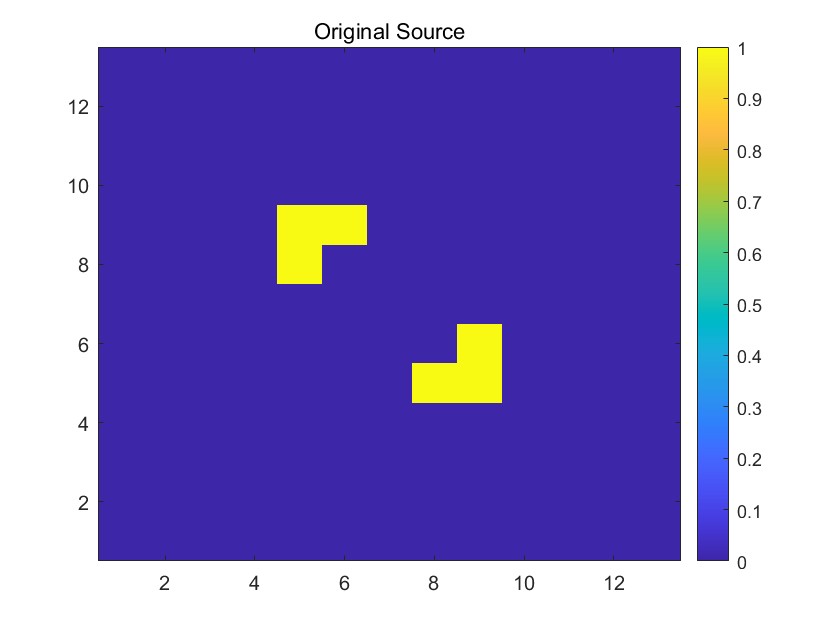}
        \caption{}
        \label{}
    \end{subfigure}
    \hfill
    \begin{subfigure}[b]{0.2\textwidth}
        \includegraphics[width=\textwidth]{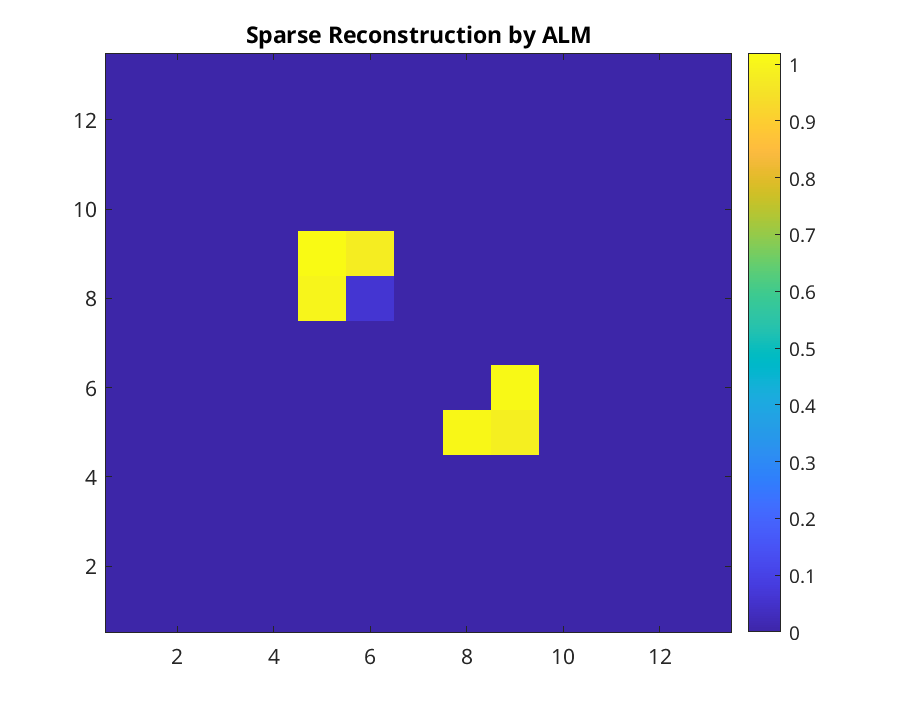}
        \caption{}
        \label{}
    \end{subfigure}
    \hfill
    \begin{subfigure}[b]{0.2\textwidth}
        \includegraphics[width=\textwidth]{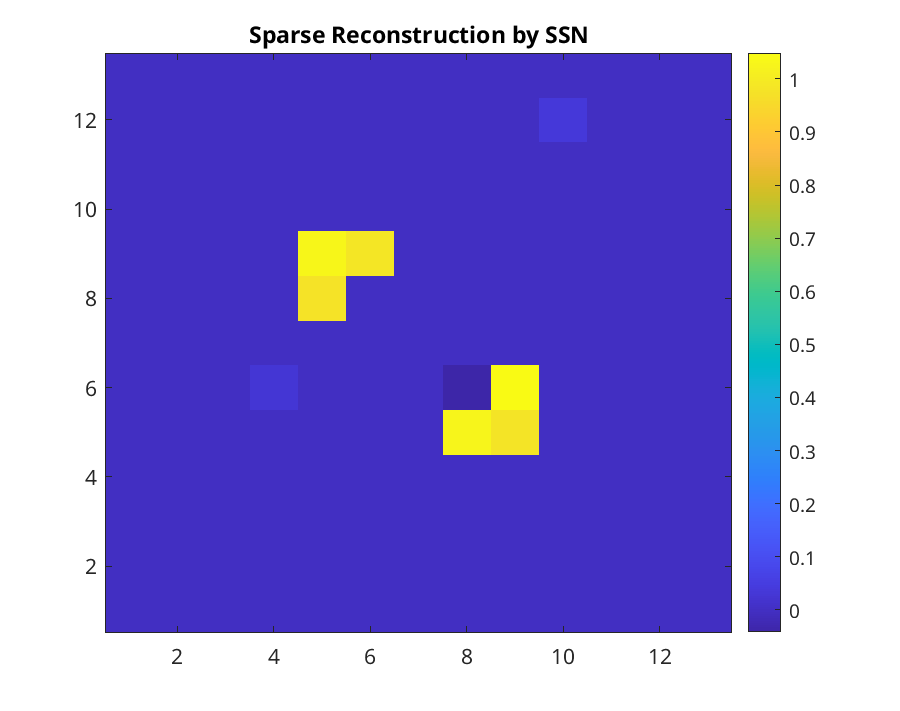}
        \caption{}
        \label{}
    \end{subfigure}
    \hfill
    \begin{subfigure}[b]{0.2\textwidth}
        \includegraphics[width=\textwidth]{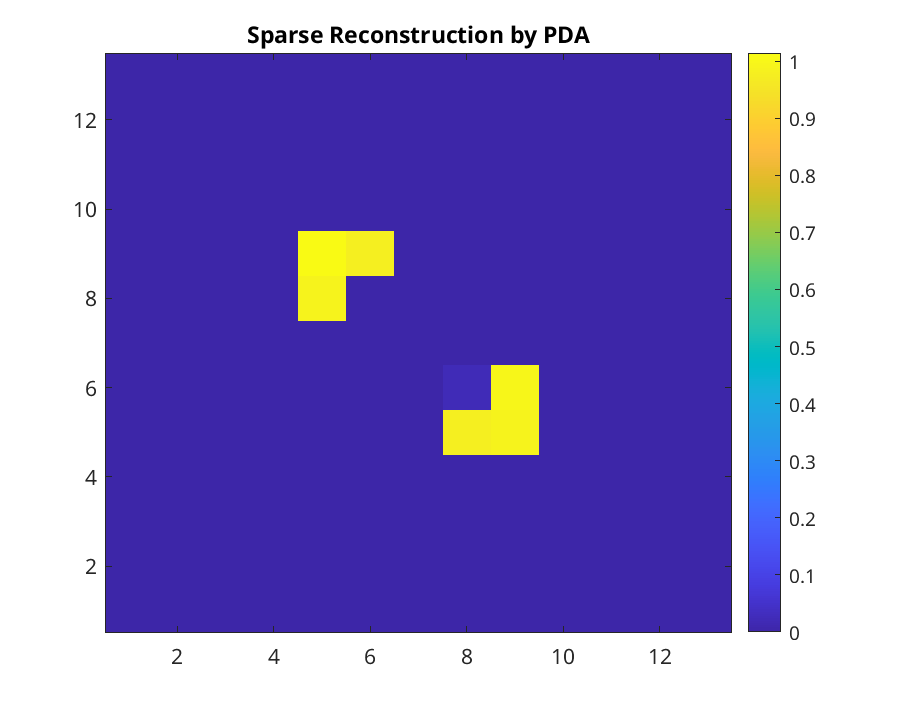}
        \caption{}
        \label{}
    \end{subfigure}
 \caption{Reconstruction of strip-shaped acoustic sources in inhomogeneous media with $k=4$, noise level 0.1\%, and $\alpha$ =2e-6,
$\alpha_0$= 2e-8. The information of the acoustic sources and the corresponding reconstruction algorithms are the same as in Figure \ref{fig:big_figure3}.}
    \label{fig:big_figure4}
\end{figure}

Figures \ref{fig:big_figure3} and \ref{fig:big_figure4} present the reconstructions of different strip sources with \%0.1 relative Gaussian noise for homogeneous and inhomogeneous media, respectively. Table \ref{tab:results:strips} collects the corresponding running times only for the homogeneous cases (since the running times for the inhomogeneous cases are nearly the same), and relative errors for the different algorithms in both homogeneous and inhomogeneous media.

Figures \ref{fig:big_figure:noise001:3d:homo} and \ref{fig:big_figure:inhomo:3d:001} present the reconstructions of different 3-dimensional sources with \%1 relative Gaussian noise  for homogeneous and inhomogeneous media, respectively. 
Table \ref{tab:results:3d} collects the running time and relative errors of ALM and SSN. It can be seen that  ALM can be 10 times faster than SSN. We did not compare with PDA, as we found it too slow. We computed the solution using $-\lambda^k$ since $\mu^*=-\lambda^*$ by \eqref{eq"lambda}, and numerically we found that $\|\mu^{K_{\text{max}}}+\lambda^{K_{\text{max}}}\|$ is generally less than $10^{-5}$. Here  ${K_{\text{max}}}$ is the maximal outer iteration number and $\mu^{K_{\text{max}}}$ is computed by \eqref{eq:mu}. 

Iterative reconstruction methods are known to be more sensitive to noise than sampling methods. A good strategy is to use a sampling method  \cite{Liu1,Liu2, Liu3} to find a good initial value, and continue the construction with iterative methods or sampling methods for more accurate reconstructions.

\begin{figure}[htbp]
    \centering
    
    
    \begin{subfigure}[b]{0.3\textwidth}
        \includegraphics[width=\textwidth]{twoballorig.jpg}
        \caption{}
        \label{fig:1d}
    \end{subfigure}
    \hfill
    \begin{subfigure}[b]{0.3\textwidth}
        \includegraphics[width=\textwidth]{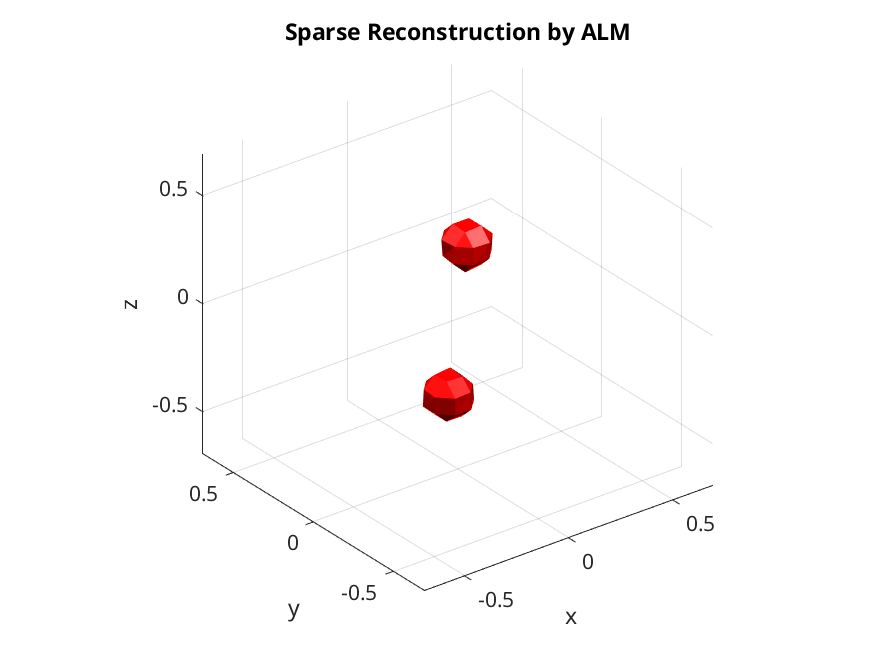}
        \caption{}
        \label{fig:1e}
    \end{subfigure}
    \hfill
    \begin{subfigure}[b]{0.3\textwidth}
        \includegraphics[width=\textwidth]{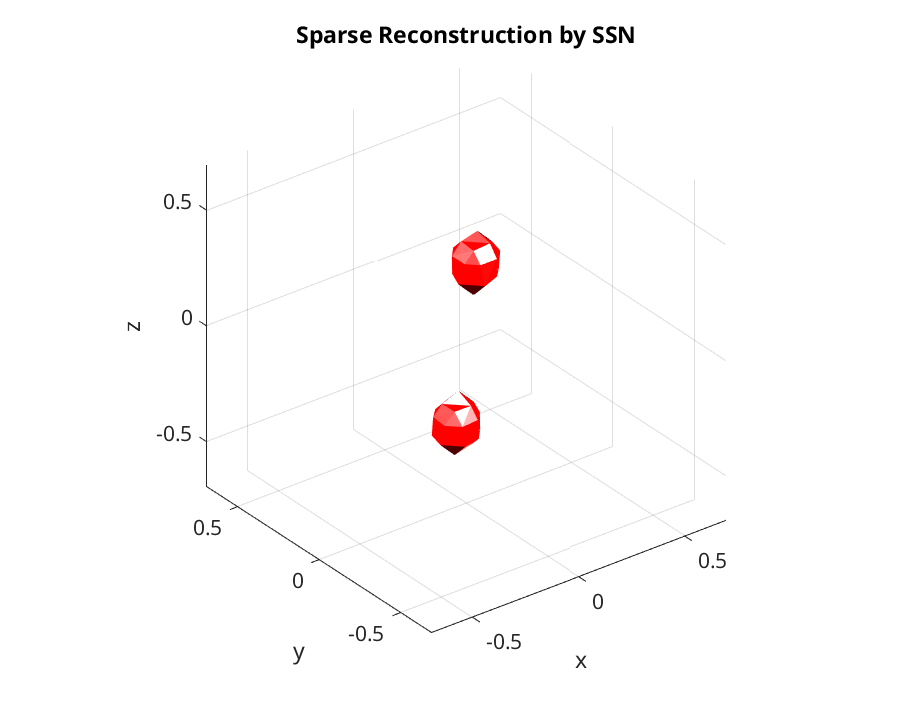}
        \caption{}
        \label{fig:1f}
    \end{subfigure}
    
    \begin{subfigure}[b]{0.3\textwidth}
        \includegraphics[width=\textwidth]{Orignalrightup.png}
        \caption{}
        \label{fig:1g}
    \end{subfigure}
    \hfill
    \begin{subfigure}[b]{0.3\textwidth}
        \includegraphics[width=\textwidth]{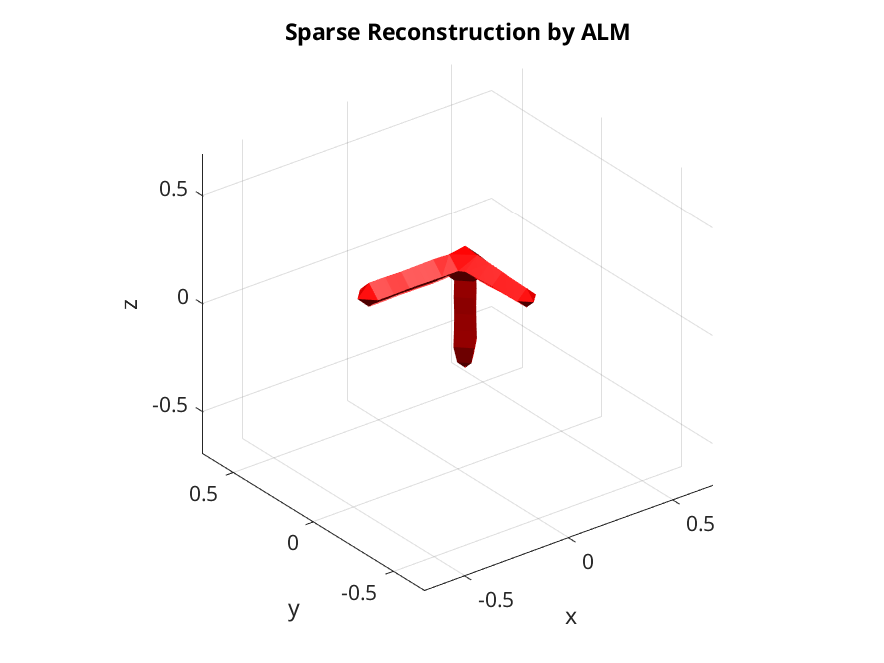}
        \caption{}
        \label{fig:1h}
    \end{subfigure}
    \hfill
    \begin{subfigure}[b]{0.3\textwidth}
        \includegraphics[width=\textwidth]{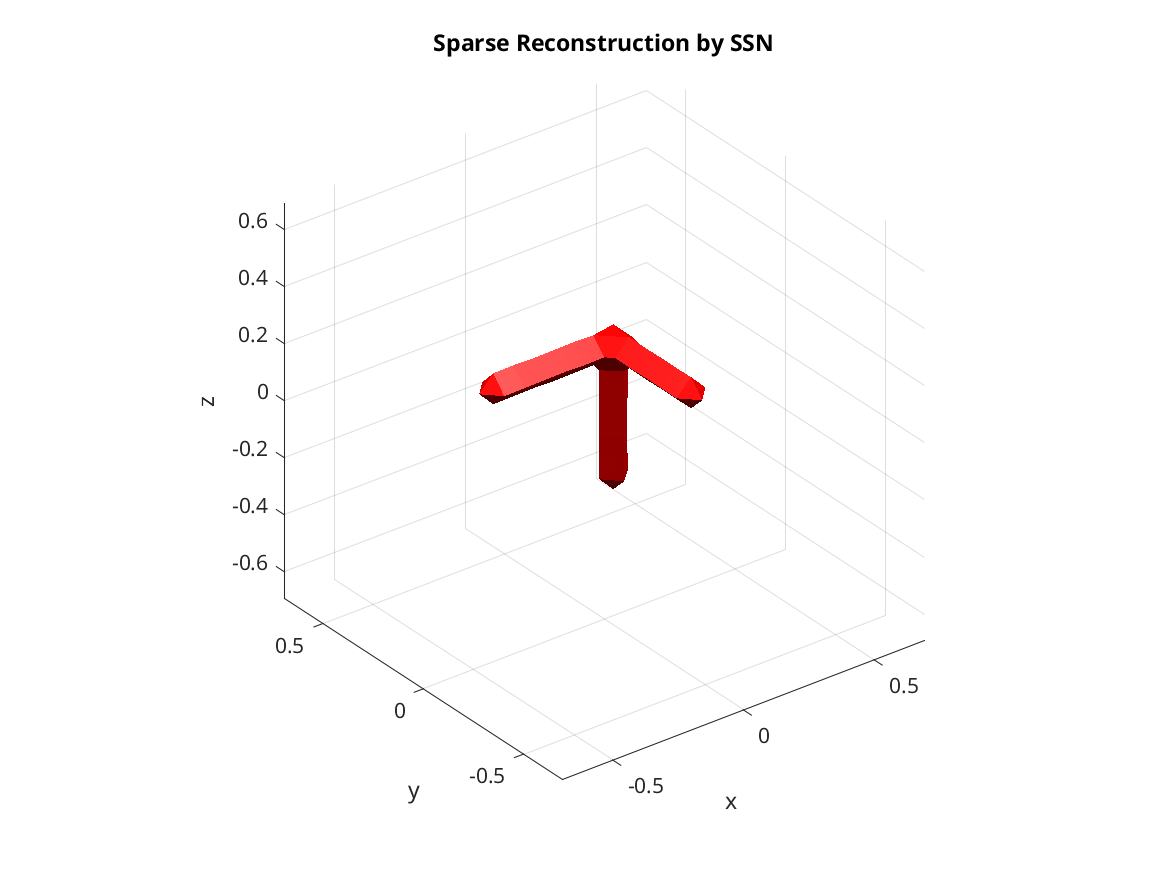}
        \caption{}
        \label{fig:1i}
    \end{subfigure}
    
    \begin{subfigure}[b]{0.3\textwidth}
        \includegraphics[width=\textwidth]{Origleftdown.png}
        \caption{}
        \label{fig:1j}
    \end{subfigure}
    \hfill
    \begin{subfigure}[b]{0.3\textwidth}
        \includegraphics[width=\textwidth]{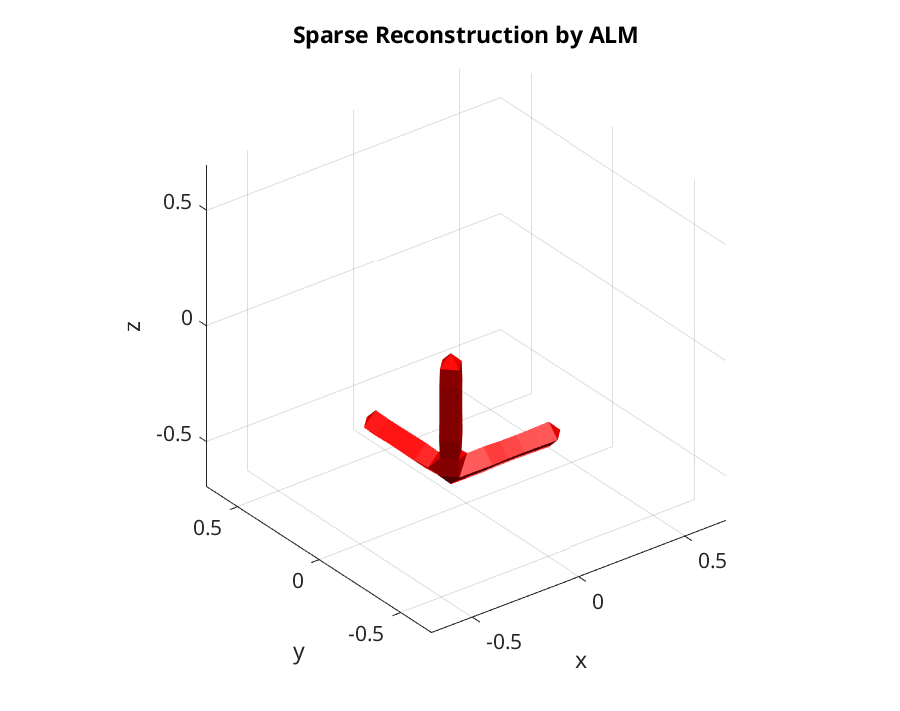}
        \caption{}
        \label{fig:1k}
    \end{subfigure}
    \hfill
    \begin{subfigure}[b]{0.3\textwidth}
        \includegraphics[width=\textwidth]{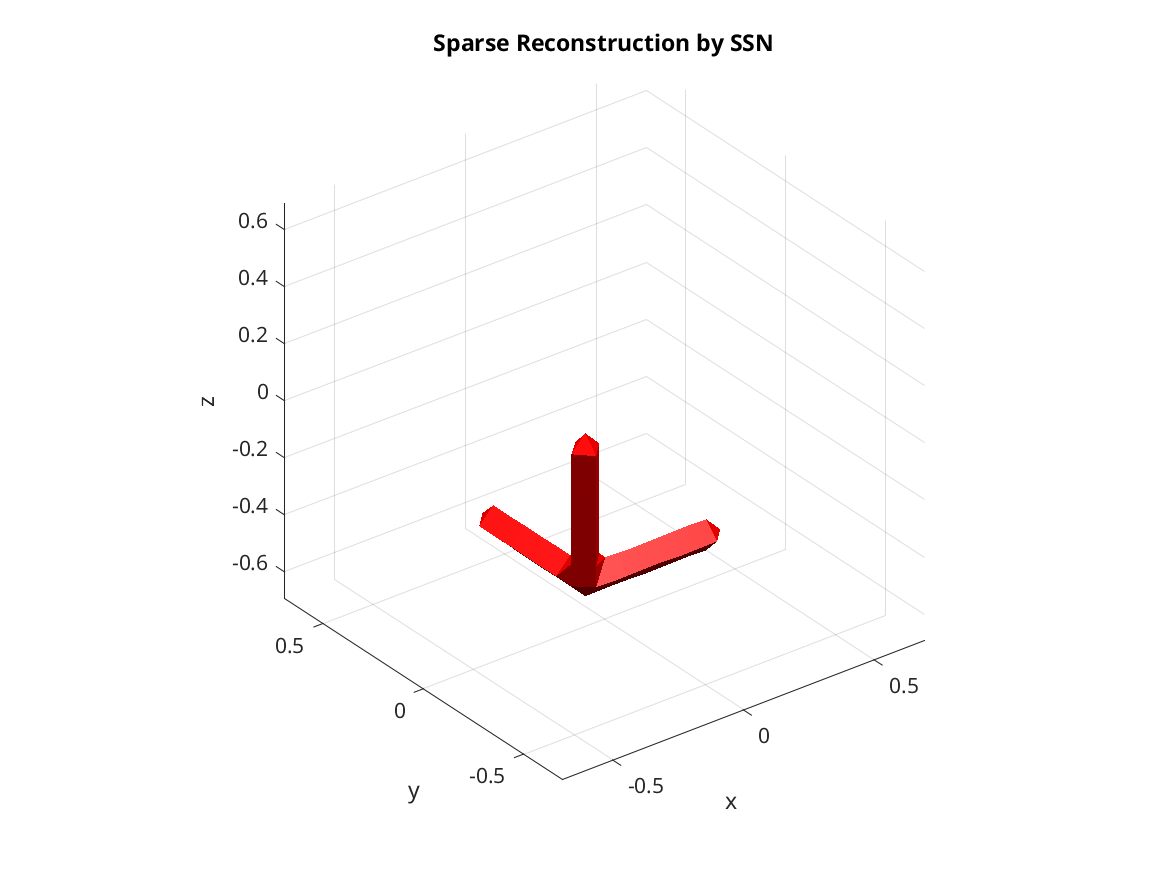}
        \caption{}
        \label{fig:1l}
    \end{subfigure}

%
    \caption{Reconstruction of 3-dimensional acoustic sources in homogeneous media with $k=6$, noise level 0.01\%, and $\alpha$ = 1e-7 together with
$\alpha_0$= 1e-9. The figures in the leftmost column are the original acoustic sources. The images in the middle and the rightmost columns are reconstructed results of ALM and SSN, respectively.  The images in the first,  second, and third rows are the sources with the two balls, the right-up tripod, and the left-down tripod, respectively. }
    \label{fig:big_figure:noise004:3d:homo}
\end{figure}

\begin{figure}[htbp]
    \centering

    \begin{subfigure}[b]{0.3\textwidth}
        \includegraphics[width=\textwidth]{twoballorig.jpg}
        \caption{}
        \label{fig:1d}
    \end{subfigure}
    \hfill
    \begin{subfigure}[b]{0.3\textwidth}
        \includegraphics[width=\textwidth]{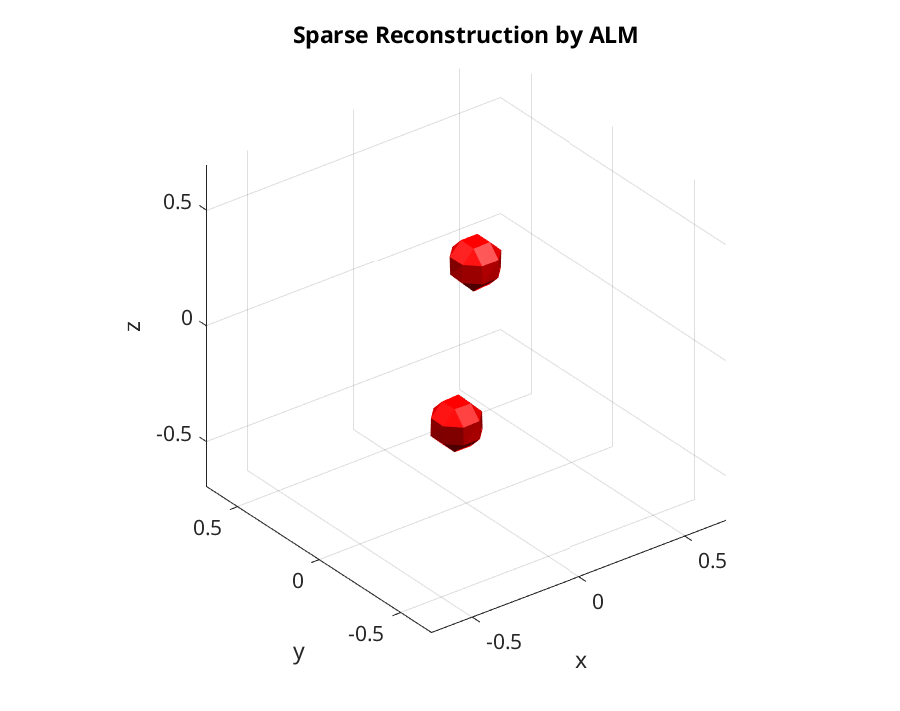}
        \caption{}
        \label{fig:1e}
    \end{subfigure}
    \hfill
    \begin{subfigure}[b]{0.3\textwidth}
        \includegraphics[width=\textwidth]{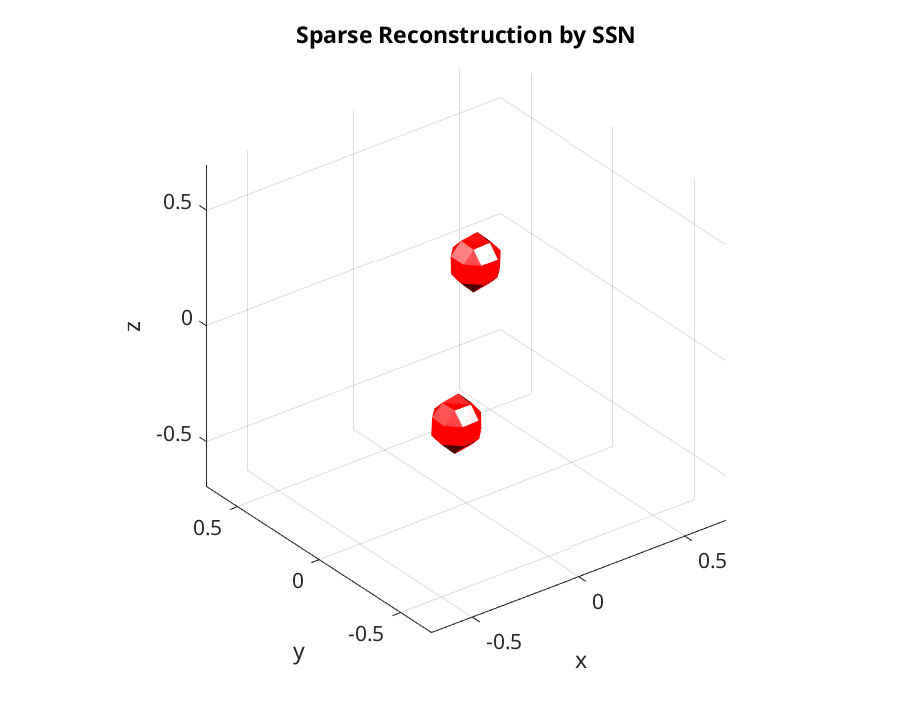}
        \caption{}
        \label{fig:1f}
    \end{subfigure}
    
    \begin{subfigure}[b]{0.3\textwidth}
        \includegraphics[width=\textwidth]{Orignalrightup.png}
        \caption{}
        \label{fig:1g}
    \end{subfigure}
    \hfill
    \begin{subfigure}[b]{0.3\textwidth}
        \includegraphics[width=\textwidth]{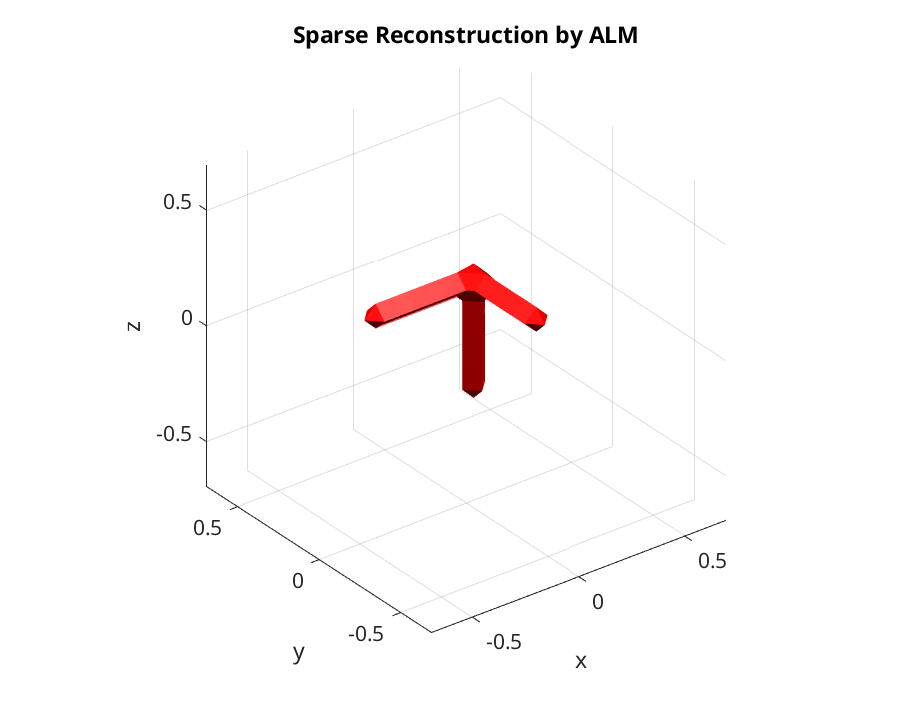}
        \caption{}
        \label{fig:1h}
    \end{subfigure}
    \hfill
    \begin{subfigure}[b]{0.3\textwidth}
        \includegraphics[width=\textwidth]{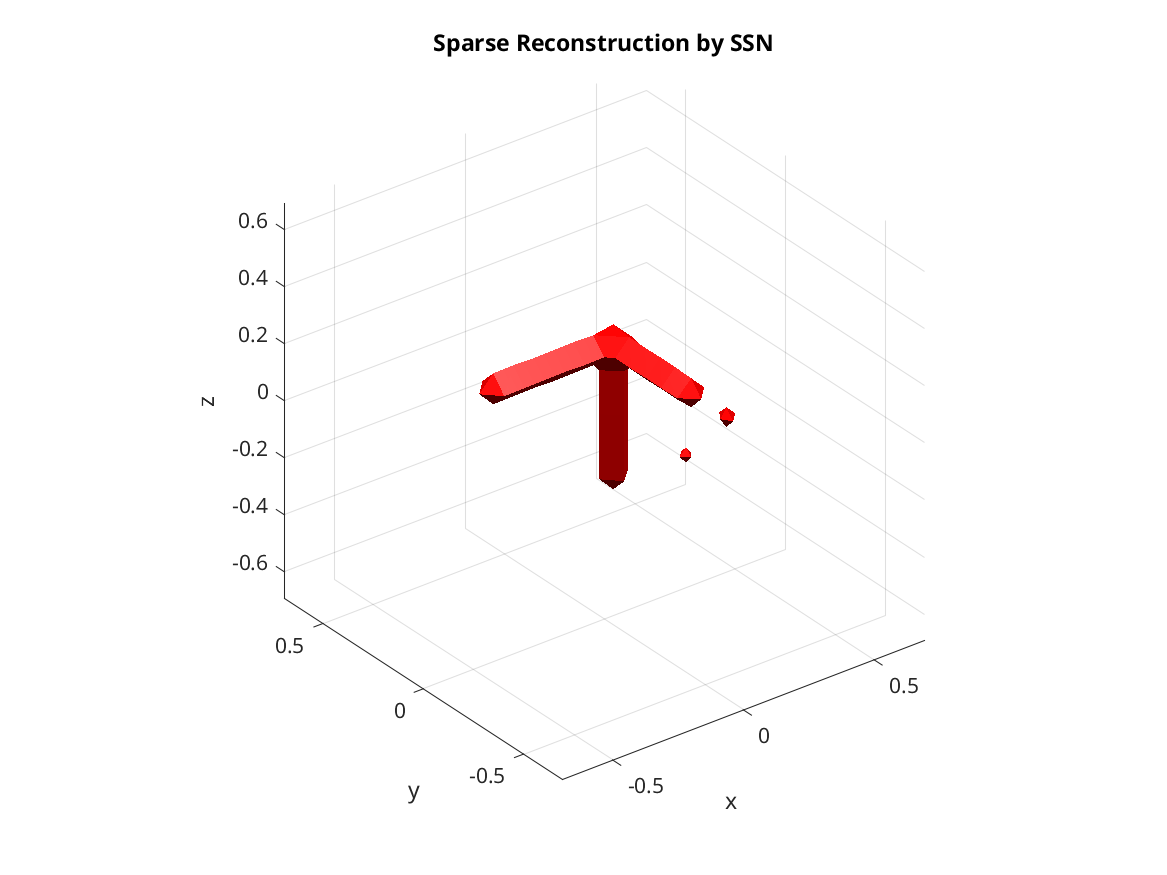}
        \caption{}
        \label{fig:1i}
    \end{subfigure}
    
    \begin{subfigure}[b]{0.3\textwidth}
        \includegraphics[width=\textwidth]{Origleftdown.png}
        \caption{}
        \label{fig:1j}
    \end{subfigure}
    \hfill
    \begin{subfigure}[b]{0.3\textwidth}
        \includegraphics[width=\textwidth]{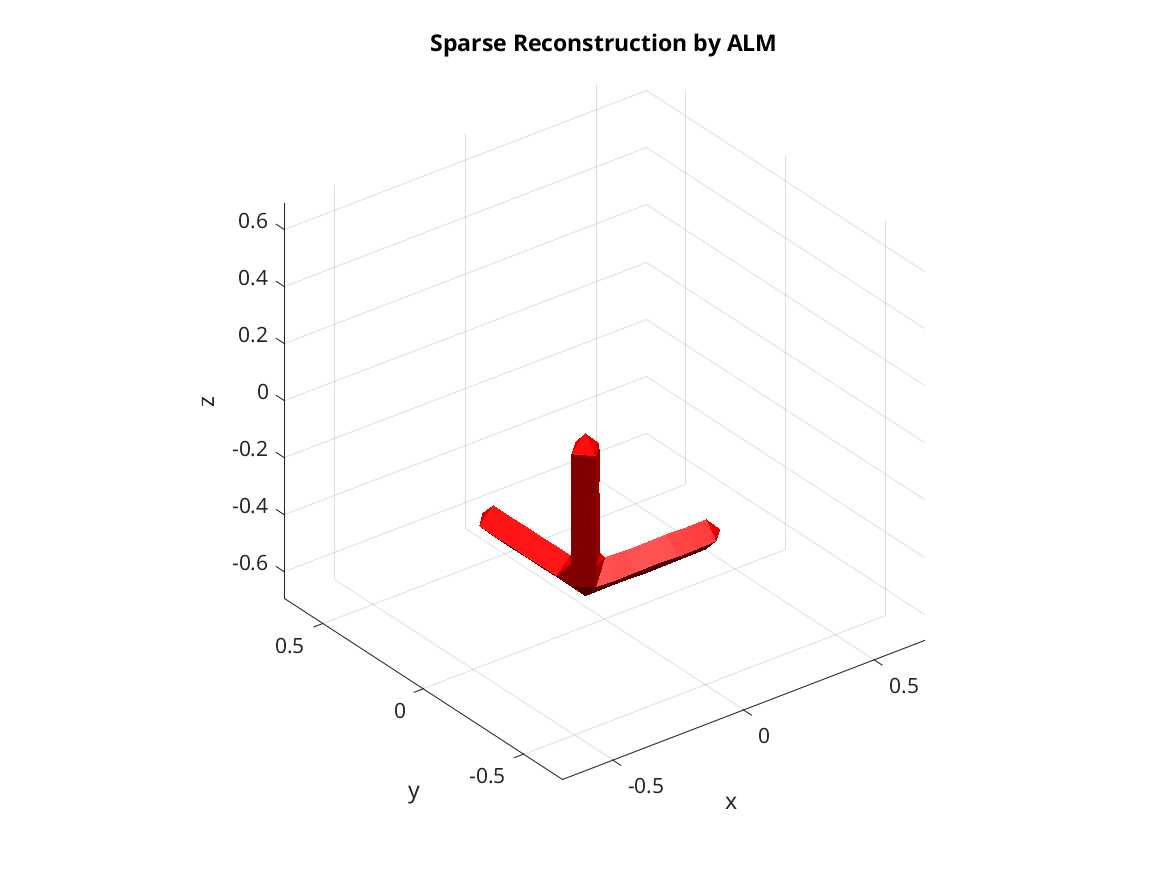}
        \caption{}
        \label{fig:1k}
    \end{subfigure}
    \hfill
    \begin{subfigure}[b]{0.3\textwidth}
        \includegraphics[width=\textwidth]{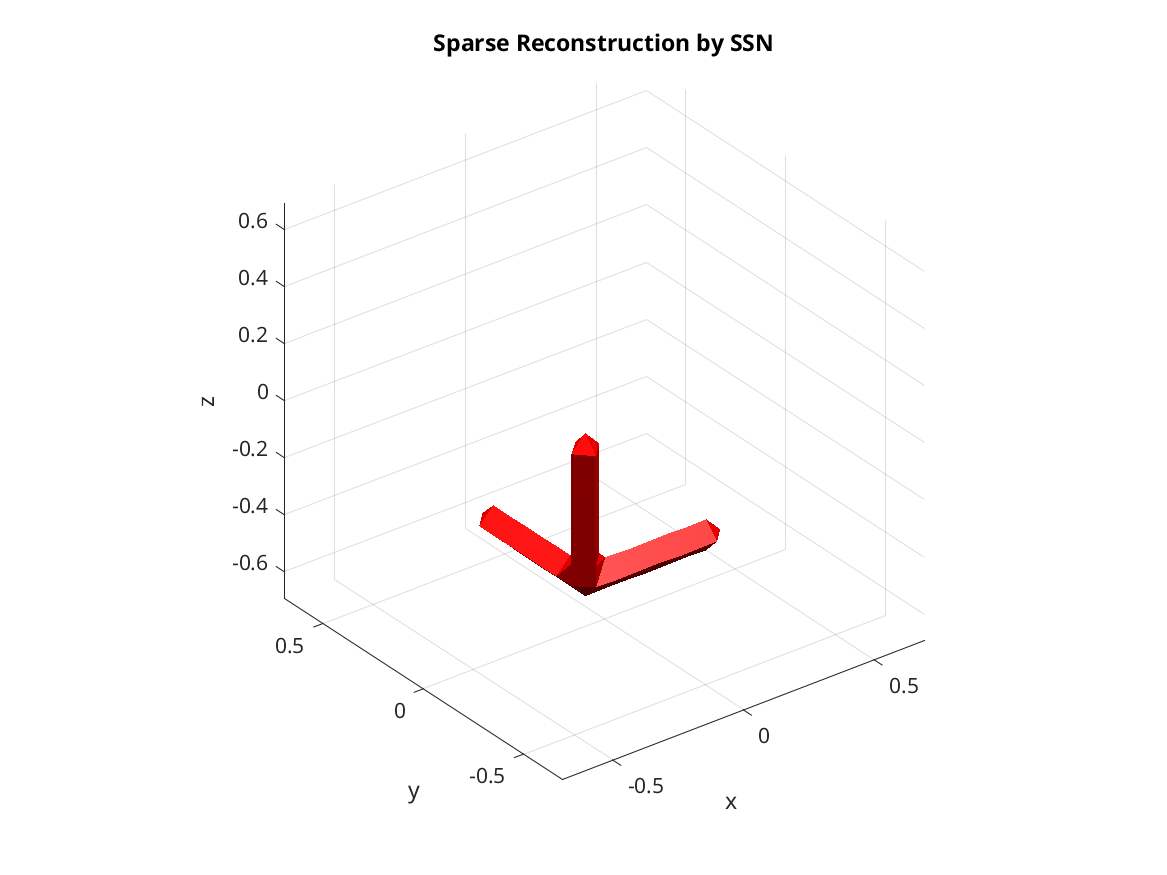}
        \caption{}
        \label{fig:1l}
    \end{subfigure}

%
    \caption{Reconstruction of a 3-dimensional acoustic sources in inhomogeneous media with $k=6$, noise level 0.01\%, and $\alpha$ = 1e-7 together with
$\alpha_0$= 1e-9. The information of the acoustic sources and the corresponding reconstruction algorithms are the same as in Figure \ref{fig:big_figure:noise004:3d:homo}.}
    \label{fig:big_figure:noise004:3d:inhomo}
\end{figure}

Figures \ref{fig:big_figure:noise004:3d:homo} and  \ref{fig:big_figure:noise004:3d:inhomo} present the reconstructions of different 3-dimensional sources with \%0.01 relative Gaussian noise  for homogeneous and inhomogeneous media, respectively. Table \ref{tab:results04:3d} collects relative errors of ALM and SSN. Since the running times are similar to those in the \%1 relative Gaussian noise case, we omit the running time in Table \ref{tab:results04:3d}.

\begin{table}[htbp]
\centering
\caption{Comparisons of the running time and relative errors for ALM, SSN, and PDA algorithms with noise level 1\%. This table corresponds to Figures \ref{fig:big_figure1} and  \ref{fig:big_figure2}. The running times, i.e., ``Times (s)" are only presented for the homogeneous cases as in Figure \ref{fig:big_figure1}.  The ``N-error"  and ``N-error (in)" are relative errors for the corresponding homogeneous and inhomogeneous cases. The notations ``one", ``four", ``six", and ``eight" correspond to the sources with one, four, six, and eight peaks and are the same as in Figures \ref{fig:big_figure1} and \ref{fig:big_figure2}.  (ALM and SSN with $\alpha$ = 9e-4,
$\alpha_0$= 1e-7, and PDA with $\alpha$ = 9e-5,
$\alpha_0$= 1e-12). }
\label{tab:results04}
\begin{tabular}{llccc}
\toprule
\textbf{Methods} & \textbf{Sources} & \textbf{Time (s)} &  \textbf{N-Error} & \textbf{N-Error (in)} \\
\midrule
\multirow{4}{*}{ALM}
& one    & 8.77     & 6.30e-02 & 4.98e-02 \\
& four    & 9.67    & 2.08e-01 & 7.69e-01 \\
& six     & 11.19    & 1.01e-00 & 1.00e-00\\
& eight    & 9.13     & 1.09e-00 & 1.07e-00 \\
\midrule
\multirow{4}{*}{SSN}
& one    & 8.62     & 1.93e-02 & 2.27e-02 \\
& four    & 10.64   & 4.32e-01  & 1.37e-01\\
& six     & 8.09    & 9.57e-01 & 9.50e-01 \\
& ten    &10.11     & 1.10e-00 & 1.07e-00 \\
\midrule
\multirow{4}{*}{PDA}
& one    & 17.12     & 6.99e-02 & 5.11e-2 \\
& four    & 15.66    & 2.49e-01 & 1.47e-01 \\
& six     & 16.71    & 9.09e-01 & 8.64e-01 \\
& eight    & 16.31     & 1.19e-00 & 1.12e-00 \\

\bottomrule
\end{tabular}
\end{table}

\begin{table}[htbp]
\centering
\caption{
Comparisons of the running time and relative errors for ALM, SSN, and PDA algorithms with noise level 0.1\%, and $\alpha$ = 2e-6, together with
$\alpha_0$= 2e-8.
This table corresponds to Figures \ref{fig:big_figure3} and \ref{fig:big_figure4}. The running times, i.e., ``Times (s)" are only presented for the homogeneous cases as in Figure \ref{fig:big_figure3}.  The ``N-error"  and ``N-error (in)" are relative errors for the corresponding homogeneous and inhomogeneous cases. The notations
 ``skew" and ``diag" represent the ``skew diagonal" and ``diagonal" strips in Figures \ref{fig:big_figure3} and \ref{fig:big_figure4}. }
\label{tab:results:strips}
\begin{tabular}{llccc}
\toprule
\textbf{Methods} & \textbf{Sources} &  \textbf{Time (s)} & \textbf{N-Error}  & \textbf{N-Error (in)}\\
\midrule
\multirow{2}{*}{ALM}
& skew     & 16.96   & 4.19e-02  &2.30e-02 \\
& diag       & 17.23   & 1.75e-02  & 4.18e-02 \\
\midrule
\multirow{2}{*}{SSN}
& skew    & 11.07   & 9.18e-02  & 4.04e-02 \\
& diag       & 11.72   &  1.79e-02  & 3.81e-02 \\
\midrule
\multirow{2}{*}{PDA}
& skew      & 18.15   & 2.40e-02 & 2.49e-02 \\
& diag    & 18.82   &  1.51e-02  & 1.70e-02\\
\bottomrule
\end{tabular}
\end{table}

\begin{table}[htbp]
\centering
\caption{Experimental results for ALM and SSN methods under different scenarios and treatments, noise level 1\%, and $\alpha$ = 5e-7 together with
$\alpha_0$= 2e-8.   Here ``N-Error" is the relative error for both homogeneous and inhomogeneous cases, corresponding to ``homo" and ``inhomo" in the table. The notations ``two balls", ``two tripods", ``right up", and ``left down" represent the acoustic sources in the first to the fourth rows in Figures \ref{fig:big_figure:noise001:3d:homo} and \ref{fig:big_figure:inhomo:3d:001}, respectively.  }
\label{tab:results:3d}
\begin{tabular}{llccc}
\toprule
\textbf{Methods} & \textbf{Sources} & \textbf{Medium} & \textbf{Time (s)} & \textbf{N-Error} \\
\midrule
\multirow{8}{*}{ALM}
& two balls    & homo   & 79.86   & 3.97e-01 \\
& two balls    & inhomo & 87.71   & 3.99e-01 \\
& two tripods    & homo   & 82.29   &  4.26e-01 \\
& two tripods     & inhomo & 88.84   & 4.96e-01 \\
& right up     & homo   & 79.49   & 3.14e-01 \\
& right up     & inhomo & 87.01   & 2.78e-01 \\
& left down    & homo   & 78.36   & 3.03e-01 \\
& left down    & inhomo & 88.69   & 2.55e-01 \\
\midrule
\multirow{8}{*}{SSN}
& two balls    & homo   & 1224.96 & 3.60e-01 \\
& two balls    & inhomo & 1304.16 & 3.26e-01 \\
& two tripods     & homo   & 1195.87 & 4.47e-01 \\
& two tripods    & inhomo & 1149.60 & 2.90e-01 \\
& right up     & homo   & 1139.21 & 4.35e-01 \\
& right up     & inhomo & 1296.09 & 3.55e-01 \\
& left down    & homo   & 1228.42 & 4.01e-01 \\
& left down    & inhomo & 1301.16 & 3.16e-01 \\
\bottomrule
\end{tabular}
\end{table}

\begin{table}[htbp]
\centering
\caption{Relative reconstruction errors for ALM and SSN methods with noise level 0.01\%, and $\alpha$ = 1e-7 together with
$\alpha_0$= 1e-9.  The ``N-error"  and ``N-error (in)" are relative errors for the corresponding homogeneous and inhomogeneous cases. The notations ``two balls", ``right up", and ``left down" represent the acoustic sources in the first to the third rows in Figures \ref{fig:big_figure:noise004:3d:homo} and \ref{fig:big_figure:noise004:3d:inhomo}, respectively. The running times are nearly the same as the 1\% noise level cases as in Table  \ref{tab:results:3d}, and we omit them here. }\label{tab:results04:3d}
\begin{tabular}{llcc}
\toprule
\textbf{Methods} & \textbf{Sources}  &  \textbf{N-Error} & \textbf{N-Error (in)} \\
\midrule
\multirow{3}{*}{ALM}
& two balls     & 1.27e-01 & 6.75e-02 \\
& right up       & 1.40e-01 & 1.43e-02 \\
& left down     & 7.33e-02 & 2.09e-02\\
\midrule
\multirow{3}{*}{SSN}
& two balls      & 3.72e-01   & 2.12e-01 \\
& right up        & 1.35e-01 & 2.29e-01\\
& left down     & 4.07e-02 & 5.36e-03 \\
\bottomrule
\end{tabular}
\end{table}

Let us focus on the mapping property of $\mathcal{V}_b$ as in Algorithm \ref{alg:alm} to illustrate the advantage of the proposed ALM based on adjoint space.
For all the 3D reconstructions in Figures \ref{fig:big_figure:noise001:3d:homo} and \ref{fig:big_figure:inhomo:3d:001},  we have 247 measured boundary data compared to the 6589 unknown discrete sources. After realization of the complex measurements and the source as in \eqref{eq:vectorize:matrixization}, we obtain the discrete  $\mathcal{V}_b$ that maps the source to the boundary measurement 
\begin{equation}
    \mathcal{V}_b: \mathbb{R}^{13178} \rightarrow \mathbb{R}^{494}.
\end{equation}
It means we only need to solve a Newton update with a linear equation with a matrix with $494\times 494$ dimensions in the adjoint space, while in the space of the acoustic source, the dimensions for the matrix associated with  a Newton update become $13178 \times 13178$. It can explain the high efficiency of the proposed ALM compared to the semismooth Newton methods that run directly on the space of the acoustic sources, as in Table \ref{tab:results:3d}.

\begin{figure}
     
    \begin{subfigure}[b]{0.2\textwidth}
        \includegraphics[width=\textwidth]{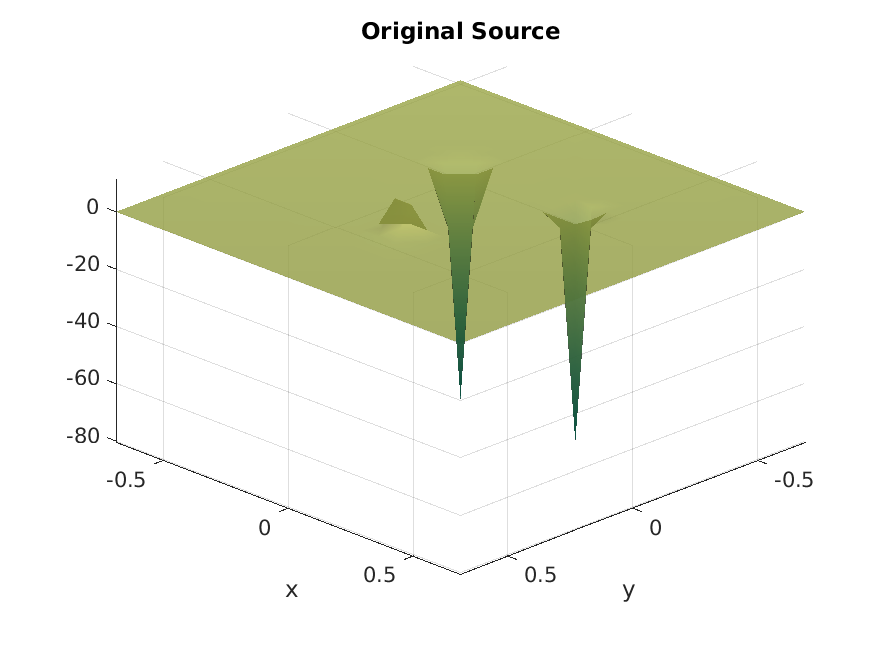}
        \caption{}
        \label{}
    \end{subfigure}
    \hfill
    \begin{subfigure}[b]{0.2\textwidth}
        \includegraphics[width=\textwidth]{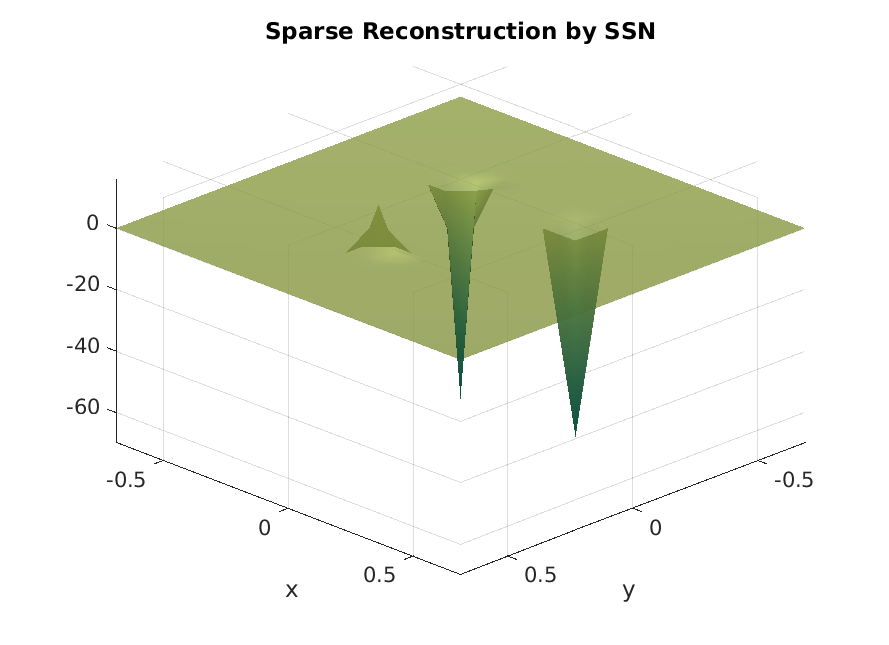}
        \caption{}
        \label{}
    \end{subfigure}
    \hfill
    \begin{subfigure}[b]{0.2\textwidth}
        \includegraphics[width=\textwidth]{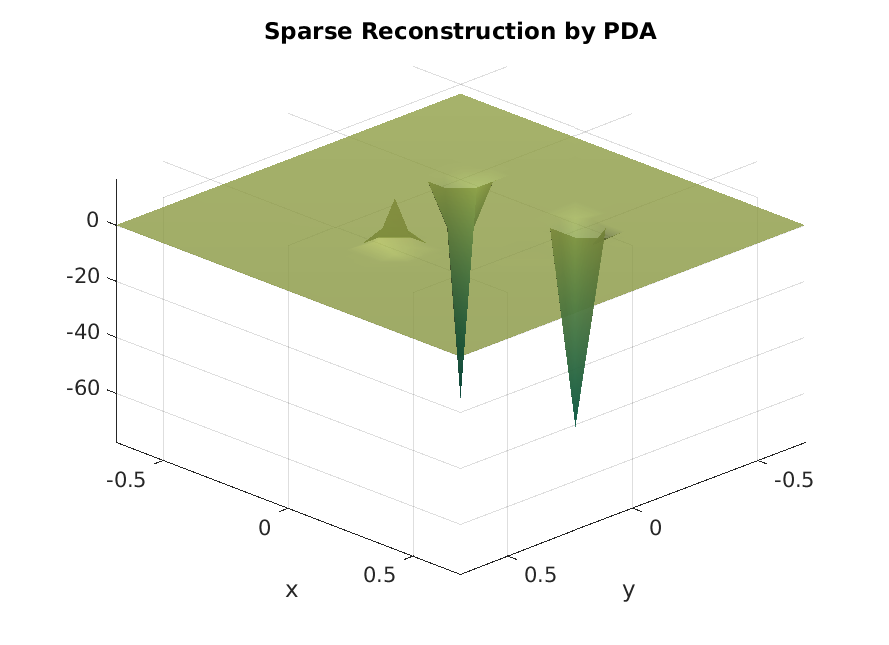}
        \caption{}
        \label{}
    \end{subfigure}
    \hfill
    \begin{subfigure}[b]{0.2\textwidth}
        \includegraphics[width=\textwidth]{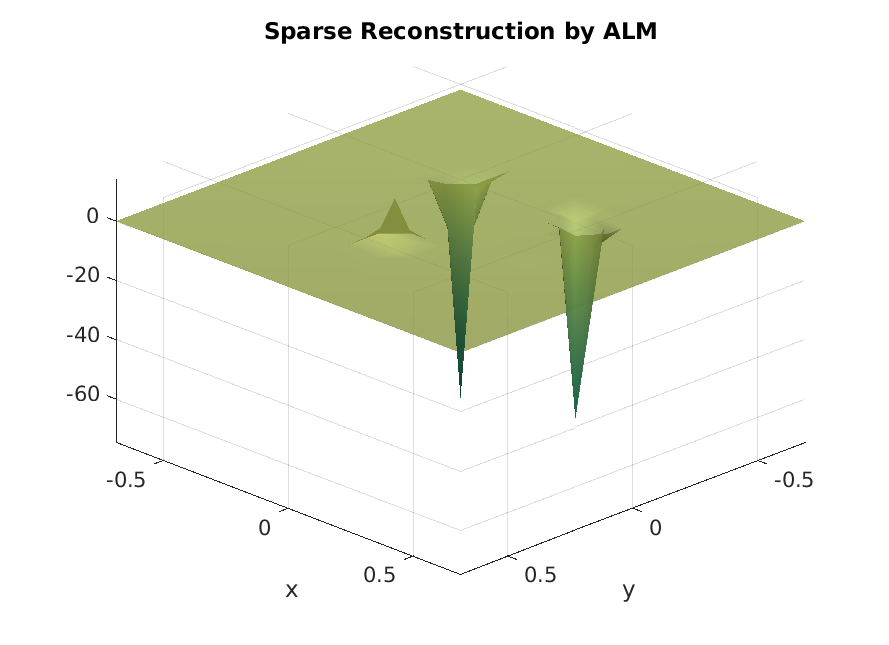}
        \caption{}
        \label{}
    \end{subfigure} \\  
    \caption{Reconstruction of an acoustic source with three peaks in homogeneous media by multiple frequencies. Here, Figure a is the original  source. Figures b, c, and d are the reconstructed figures using the SSN, PDA, and ALM algorithms, respectively.}
    \label{fig:multi:point:homo}
\end{figure}

\begin{figure}
      \begin{subfigure}[b]{0.2\textwidth}
        \includegraphics[width=\textwidth]{2d_original_source.png}
        \caption{}
        \label{}
    \end{subfigure}
    \hfill
    \begin{subfigure}[b]{0.2\textwidth}
        \includegraphics[width=\textwidth]{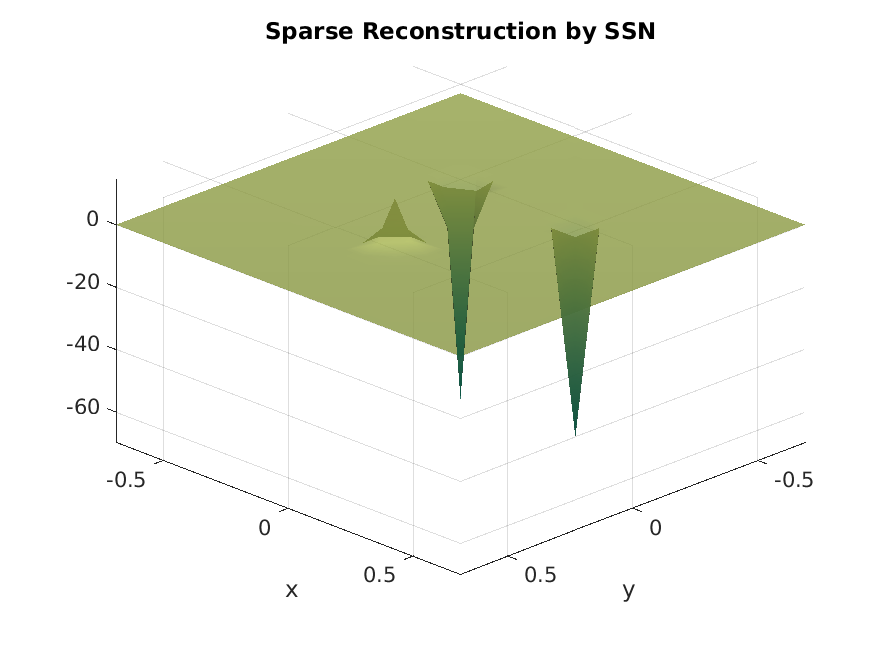}
        \caption{}
        \label{}
    \end{subfigure}
    \hfill
    \begin{subfigure}[b]{0.2\textwidth}
        \includegraphics[width=\textwidth]{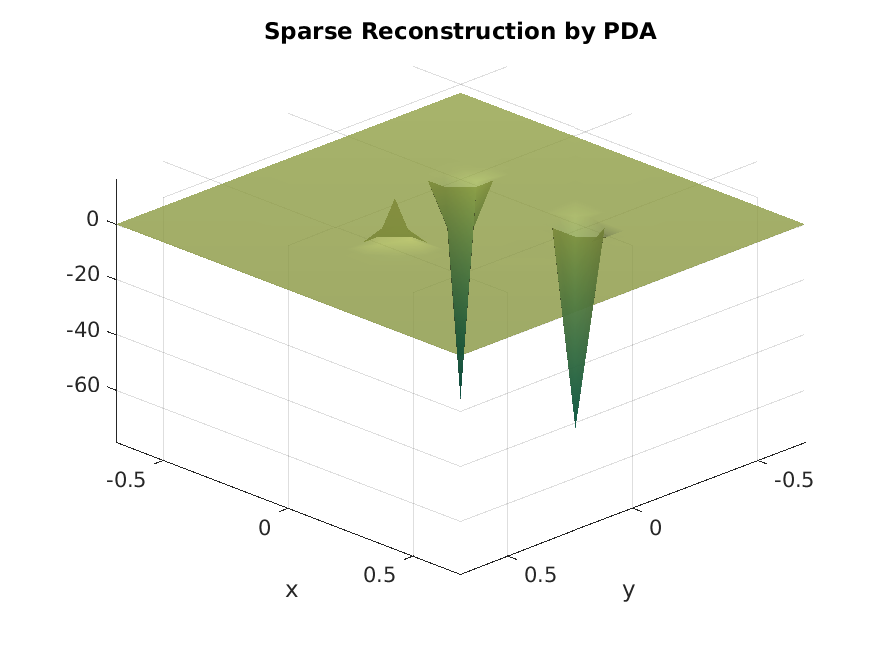}
        \caption{}
        \label{}
    \end{subfigure}
    \hfill
    \begin{subfigure}[b]{0.2\textwidth}
        \includegraphics[width=\textwidth]{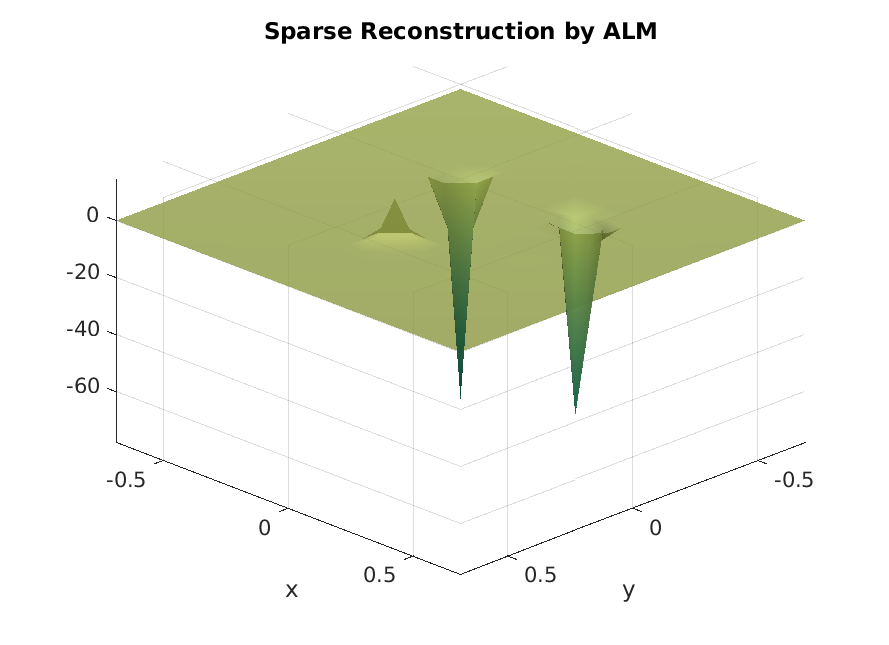}
        \caption{}
        \label{}
    \end{subfigure}
    \caption{Reconstruction of an acoustic source with three peaks in inhomogeneous media by multiple frequencies. Here, Figure a is the original acoustic source. Figures b, c, and d are the reconstructed solutions with algorithms SSN, PDA, and ALM, respectively. The inhomogeneous medium is the same as in Table \ref{tab:results04}. }
    \label{fig:multi:point:inhomo}
\end{figure}

\begin{figure}   
    \begin{subfigure}[b]{0.3\textwidth}
        \includegraphics[width=\textwidth]{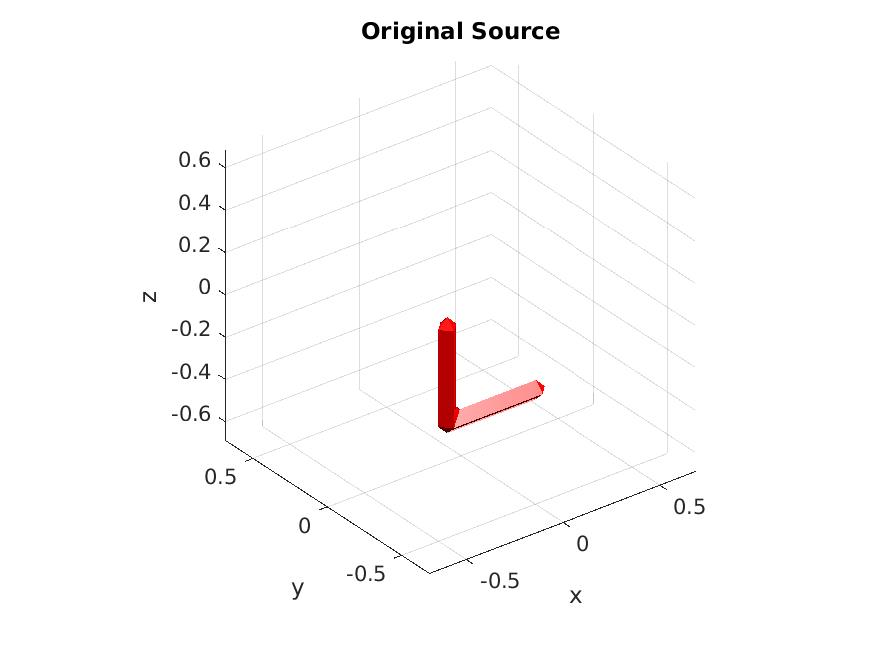}
        \caption{}
        \label{}
    \end{subfigure}
    \begin{subfigure}[b]{0.3\textwidth}
        \includegraphics[width=\textwidth]{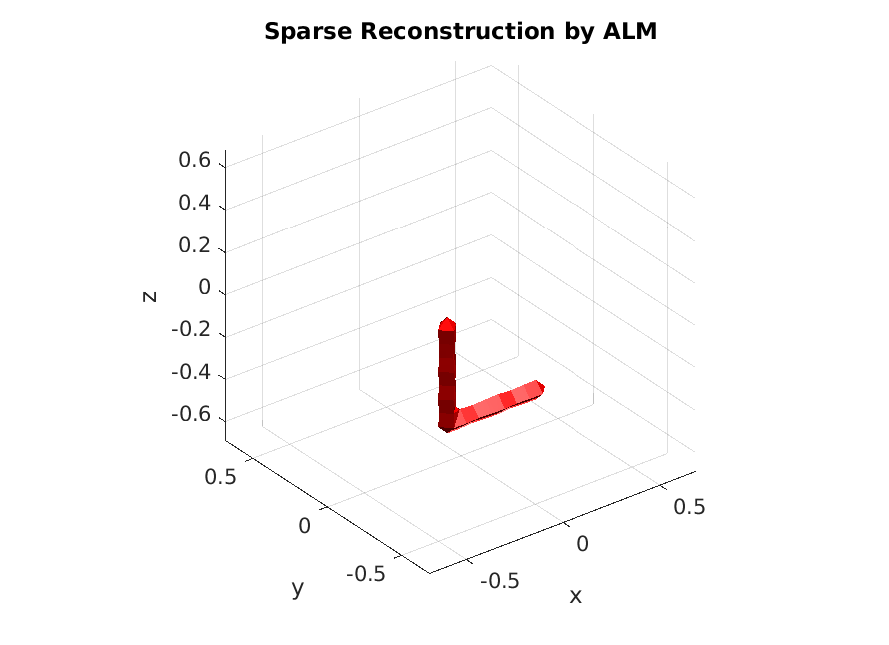}
               \caption{}
        \label{}
    \end{subfigure}
    \begin{subfigure}[b]{0.3\textwidth}
        \includegraphics[width=\textwidth]{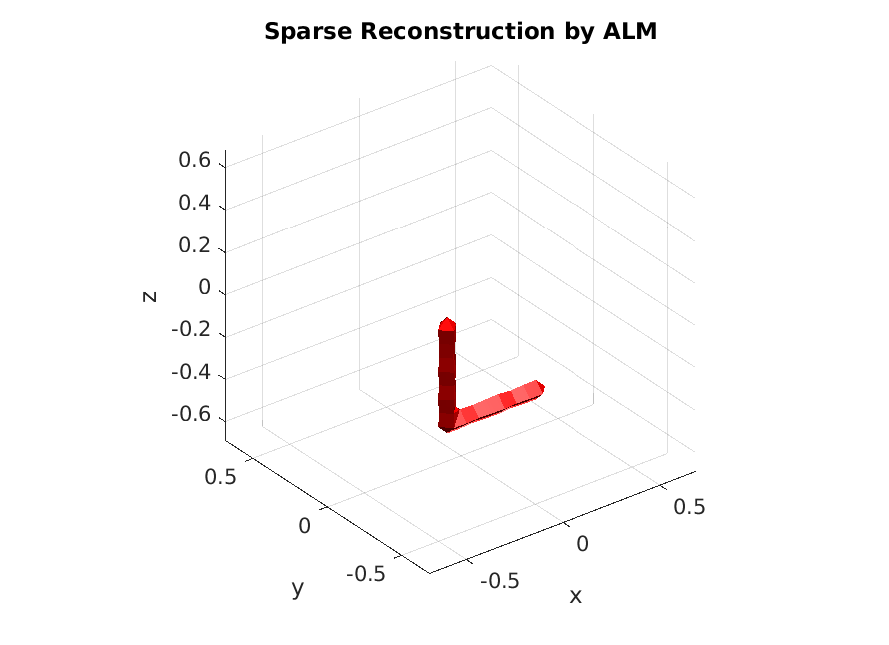}
             \caption{}
        \label{}
    \end{subfigure}
   

    \caption{Reconstruction of a 3D acoustic source. Figure a is the original acoustic source. Figures b and c are the reconstructions in homogeneous and inhomogeneous media with multiple frequencies, respectively. The parameter for the inhomogeneous medium is the same as in Table \ref{tab:results:3d}.}
    \label{fig:multi:rect}

\end{figure}

Now, we turn to the reconstructions with multiple frequencies cases, as shown in Table \ref{tab:multi}, the continuation strategy starting from the lowest frequency \cite{Bao3} can indeed bring an effective reconstruction. We set $k_{\min}=2$, $k_{\max}=6$, and the discretization of the grids is determined by the highest wavenumber. For all the compared algorithms including ALM, SSN, and PDA, the reconstructed solutions by the lower wavenumber as the initial value for the reconstruction of the next higher wavenumber \cite{Bao3,Bao2, Bao4}. For reconstructed acoustic sources, see Figures \ref{fig:multi:point:homo} and \ref{fig:multi:point:inhomo} for the two-dimensional cases and Figure \ref{fig:multi:rect} for the three-dimensional cases.

\begin{table}[]
\centering
\caption{Comparisons of the running time and relative errors for ALM, SSN, and PDA algorithms by multiple frequencies and with noise level 1\%. The running times, i.e., ``Times (s)" is the same as in Table \ref{tab:results:3d}. The ``N-error($k$)" ($k=2$, $k=4$, or $k=6$)  represents relative errors for the corresponding homogeneous (``homo" in the table) and inhomogeneous (``inhomo" in the table) cases with the frequency $k$.  The regularization parameters are as follows. For the 2D (two-dimensional) cases, ALM, i.e., ALM(2D), is with $\alpha$ = 9e-5,
$\alpha_0$= 1e-7. PDA, i.e., PDA(2D), is with $\alpha$ = 9e-5,
$\alpha_0$= 1e-12, and  SSN, i.e., SSN(2D), is with $\alpha$ = 9e-4,
$\alpha_0$= 1e-7. For 3D (three-dimensional) case, ALM, i.e., ALM(3D), is with $\alpha$ = 5e-7,
$\alpha_0$= 2e-8. The original acoustic sources and the reconstructed sources can be found in Figures \ref{fig:multi:point:homo} and \ref{fig:multi:point:inhomo} for the 2D cases and Figure \ref{fig:multi:rect} for the 3D cases.}\label{tab:multi}
\begin{tabular}{ccccccc}
\toprule
Methods& Medium& N-Error ($2$)& N-Error ($4$)& N-Error ($6$)&Time (s)\\
\midrule
SSN(2D)& homo& 3.83e-01& 2.81e-01&2.43e-01&14.88\\
SSN(2D)& inhomo& 3.81e-01& 2.61e-01&1.64e-01&24.86\\
PDA(2D)& homo& 6.98e-01& 2.45e-01&1.39e-01&89.97\\
PDA(2D)& inhomo& 6.64e-01&1.70e-01 &1.07e-01& 79.75\\
ALM(2D)& homo&6.68e-01&3.54e-01&1.84e-01& 20.39\\
ALM(2D)& inhomo&6.51e-01&2.59e-01& 1.51e-01&23.83\\
ALM(3D)& homo&  1.43e-01& 1.25e-01& 1.09e-01&165.42\\
ALM(3D)& inhomo& 1.46e-01& 1.19e-01 & 1.07e-01& 218.01\\

\bottomrule
\end{tabular}
\end{table}

\section{Conclusions}\label{sec:con}
We propose a semismooth Newton-based augmented Lagrangian method in the adjoint space of the measurements. As shown in the three-dimensional numerical examples, reconstruction can be greatly accelerated by more than a factor of 10 when the number of measurements is much smaller than the number of unknown acoustic sources.
 The multi-frequency scattering data \cite{AHLS, Bao1,Bao2,EV} is also considered. It would be interesting to study the partial measurements \cite{KSU, LLSL}.
We believe this framework can also benefit inverse medium problems and inverse electromagnetic wave scattering problems. The stability of the inverse source problem \cite{Ilu, LZZ} is also very interesting, especially for the proposed algorithm with a certain noise level.

\section*{Acknowledgements}

Nirui Tan and Hongpeng Sun acknowledge the support of the National Natural Science Foundation of China under grant No. \,12271521, National Key R\&D Program of China (2022ZD0116800), and Beijing Natural Science Foundation No. Z210001.

\bibliographystyle{plain}

\end{document}